\documentclass[mathpazo]{cicp}
\usepackage[margin=2.5cm]{geometry}
\usepackage{setspace}
\usepackage{lmodern}
\usepackage{amsmath,amssymb}
\usepackage{amsfonts,amsthm}
\usepackage{lineno}
\usepackage{graphicx}
\usepackage{bm, bbm}
\usepackage{algorithm}
\usepackage{physics}
\usepackage{color}
\usepackage{pxfonts}
\usepackage[footnotesize,bf]{caption}
\usepackage{hhline}
\usepackage[T1]{fontenc}
\usepackage{placeins}

\usepackage{float}
\usepackage{subfig}
\usepackage{tabularx}

\usepackage[cp1252]{inputenc}

\newcommand{\D}{D}
\newcommand{\E}{\mathbb{E}}

\newcommand{\overbar}[1]{\mkern 1.5mu\overline{\mkern-1.5mu#1\mkern-1.5mu}\mkern 1.5mu}

\renewcommand{\arraystretch}{1.5}
\graphicspath{{figure/}}

\begin{document}
\title{Bi-orthogonal fPINN: A physics-informed neural network  method for solving time-dependent stochastic fractional PDEs}

 \author[Ma L et.~al.]{Lei Ma\affil{1}, Rongxin Li\affil{1}, Fanhai Zeng\affil{2}\comma\corrauth,
       Ling Guo\affil{1}\comma${}^*$, and George Em Karniadakis\affil{3}}
 \address{\affilnum{1}\ Department of Mathematics, Shanghai Normal University, Shanghai 200234, China. \\
           \affilnum{2}\ School of Mathematics, Shandong University,  Shandong 250100, China.\\
           \affilnum{2}\ Division of Applied Mathematics, Brown University, Providence, RI 02906, USA,
            and School of Engineering, Brown University, Providence, RI 02906, USA.}
 \emails{{\tt fanhai\_zeng@sdu.edu.cn} (F.~Zeng), {\tt lguo@shnu.edu.cn} (L.~Guo),
          {\tt george\_karniadakis@Brown.edu} (G.~Karniadakis)}

\begin{abstract}
Fractional partial differential equations (FPDEs) can effectively represent anomalous transport and nonlocal interactions. However, inherent uncertainties arise naturally in real applications due to random forcing or unknown material properties. Mathematical models considering nonlocal interactions with uncertainty quantification can be formulated as stochastic fractional partial differential equations (SFPDEs). There are many challenges in solving SFPDEs numerically, especially for long-time integration since such problems are high-dimensional and nonlocal.
Here, we combine the bi-orthogonal (BO) method for representing stochastic processes with physics-informed neural networks (PINNs) for solving partial differential equations to formulate the bi-orthogonal PINN method (BO-fPINN) for solving time-dependent SFPDEs. Specifically,
we introduce a deep neural network for the stochastic solution of the time-dependent SFPDEs, and include the BO constraints in the loss function following a weak formulation.
Since automatic differentiation is not
currently  applicable to fractional derivatives, we employ discretization on a  grid to to compute the fractional derivatives of the neural network output. The weak formulation loss function of the BO-fPINN method can overcome some drawbacks of the BO methods and thus can be used to solve SFPDEs with eigenvalue crossings. Moreover, the BO-fPINN method can be used for inverse SFPDEs with the same framework and same computational complexity as for forward problems. We demonstrate the effectiveness of the BO-fPINN method for different benchmark problems. Specifically, we first consider an SFPDE with eigenvalue crossing and  obtain good results while the original BO method fails.  We then solve several forward and inverse problems governed by SFPDEs, including problems with noisy initial conditions. We study the effect of the fractional order as well as the number of the BO modes on the accuracy of the BO-fPINN method. The results demonstrate the flexibility and efficiency of the proposed method, especially for inverse problems. We also present a simple example of transfer learning (for the fractional order) that can help in accelerating the training of BO-fPINN for SFPDEs.
Taken together, the simulation results show that the BO-fPINN method can be employed to effectively solve time-dependent SFPDEs and may provide a reliable computational strategy for real applications exhibiting anomalous transport.
\end{abstract}


\maketitle

\section{Introduction}
Fractional partial differential equations (FPDEs) can effectively represent anomalous transport and nonlocal interactions in  engineering and medical fields, e.g., porous media, viscoelasticity, non-Newtonian fluid mechanics, soft tissue mechanics, etc.  \cite{benson2000application,mainardi2010fractional,song2016fractional,yu2016fractional,lischke2020fractional}. However, simulating real-world  applications requires modelers to consider many uncertain factors, such as material properties, random  forcing terms, experimental measurement errors, and the complexity of geometric regions with random roughness. These uncertain factors may have an important impact on the system evolution, especially in long-term  forecasting, and hence quantifying uncertainty is very important. Following a probability framework, the uncertainty is usually modeled as a random field \cite{tao2015recent}, and therefore modeling nonlocal interactions with uncertainty requires the formulation of fractional partial differential equation with random inputs. Although there have been some achievements in the numerical solution of stochastic fractional partial differential equations (SFPDEs) \cite{liu2018quasi, liu2018fourier, yokoyama2016regularity}, the design of reliable algorithms that can tackle high-dimensions and long-term integration is still an open challenge in the context of solving efficiently time-dependent stochastic fractional partial differential equations (SFPDEs).

Monte Carlo (MC) and Quasi Monte Carlo (QMC) methods are common methods for solving differential equations with random inputs, with the statistical values (e.g.,  mean and variance) of the random solutions obtained by numerically solving a set of corresponding deterministic differential equations at sample points. However, such  methods have slow convergence speed and tax computational resources heavily for large complex systems \cite{kuo2012quasi, kuo2016application}.
In recent years, Polynomial Chaos (PC) has been widely used to solve partial differential equations with random inputs, which can be regarded as a spectral approximation method in random space. The basic idea is to first employ a Karhunen Lo\`{e}ve expansion of the random field to reduce the dimension of the infinite-dimensional random inputs and represent them  with a finite-dimensional series, and subsequently construct a polynomial surrogate model of the random solution. The polynomial expansion coefficients are solved by the stochastic Galerkin method or the stochastic collocation method, and then the statistics of interest can be readily computed \cite{xiu2002wiener, xiu2010numerical}. For more applications and introduction of the PC method, please refer to \cite{tao2015recent, guo2020constructing}. Although the PC method is very effective, with the increase of the random input dimension, the number of basis functions increases exponentially leading to the well known problem of the curse-of-dimensionality. In the standard PC method, the basis function of the random space does not evolve with time, which brings difficulties to the solution of time evolution nonlinear problems.

In order to solve time-dependent stochastic partial differential equations more efficiently, the authors of \cite{sapsis2009dynamically} first proposed a generalized Karhunen-Lo\`{e}ve expansion that approximates second-order random fields. This method is different from the PC method in that both the physical space and random space basis functions evolve with time. To derive equations for the mean value, the evolutionary system of spatial basis functions and stochastic basis functions, the authors of \cite{sapsis2011dynamically,sapsis2012dynamical} developed the dynamically orthogonal (DO) condition by constraining the spatial basis functions. In follow up work,  the authors of \cite{cheng2013dynamically1, cheng2013dynamically2} developed bi-orthogonal (BO) conditions for constraining the spatial basis functions and random basis functions. The error analysis of the DO method was first given in \cite{musharbash2015error}. Choi et al. \cite{choi2014equivalence} proved the formal equivalence of the DO and BO methods, and demonstrated that the two methods differ only in terms of computability and numerical artifacts.
In this paper, we focus on the BO method, which can effectively solve the long-term stochastic system evolution problem. However, BO also faces some problems  if the eigenvalues of the system have crossing or if the covariance matrix is singular, in which case catasrophic instabilities may appear spontaneously \cite{cheng2013dynamically2, choi2014equivalence}. To this end, it is necessary to further combine BO with other methods, such as the MC and PC methods, to improve the BO method, e.g., see details in  \cite{babaee2017robust}.

Recently, scientific machine learning has received extensive attention in numerically solving data-driven partial differential equations, mainly including learning methods based on Gaussian process regression \cite{graepel2003solving, sarkka2011linear, bilionis2016probabilistic, raissi2018numerical}, and approximation algorithms based on deep neural networks \cite{lagaris1998artificial, lagaris2000neural, khoo2021solving, raissi2017physics}. In view of the strong representation ability of deep neural networks to approximate nonlinear operators \cite{chen1993approximations, chen1995universal}, in this paper we focus on Physics-informed neural networks (PINNs)~\cite{raissi2017physics}. The basic idea of PINNs is to encode physical laws (e.g., expressed by partial differential equations as conservation laws) in the neural network, by modifying the loss function along with the competing terms of limited data mismatch, such as data corresponding to initial and boundary conditions or scattered measurements in the domain. Compared with the traditional numerical methods such as finite difference or finite element method, the PINN algorithm is a meshless method, so it has great advantages in solving inverse problems governed by nonlinear differential equations.
In addition, with the help of the automatic differentiation available in current deep learning frameworks, the code of the PINN algorithm is concise and easy to understand and  implement, so the algorithm becomes a powerful tool for numerically solving general  differential equations, see references \cite{Raissi2017PhysicsID, raissi2019deep, tartakovsky2020physics, kissas2020machine, gao2022failure, gao2023failure} and related literature for details.  A deep learning algorithm for solving fractional PDEs was first presented in \cite{pang2019fpinns}. Since the back-propagation algorithm of deep learning cannot be used to calculate fractional partial derivatives (due to the chain rule), it is necessary to adopt the traditional discrete scheme such as finite difference or finite element method to discretize  the derivative of the neural network surrogate model, and then calculate the loss function of the equation to  optimize the neural network parameters. The results in \cite{pang2019fpinns} demonstrated the effectiveness of the deep learning algorithm in solving fractional partial differential equations, especially for parameter estimation in  inverse problems.

On the other hand, deep learning has also been developed rapidly in numerically solving stochastic PDEs, and achieved some important results, including solving high-dimensional parabolic equations, backward stochastic differential equations, and uncertainty quantization; for details, see references \cite{han2017deep, raissi2018forward, zhu2018bayesian, tartakovsky2018learning, yan2019adaptive}. Zhang et al. proposed a deep learning-based method for solving partial differential equations with random inputs \cite{zhang2019quantifying}. Further, in order to solve time-dependent stochastic partial differential equations more efficiently, Zhang et al.  \cite{zhang2020learning} studied the dynamic-orthogonal (DO) and bi-orthogonal (BO) deep neural network surrogate models, which could overcome some difficulties encountered previously in solving  uncertainty quantification problems with such methods. Recently, Zhao et al. proposed a dynamic orthogonal and biorthogonal Galerkin spectral method for solving stochastic fractional nonlinear equations \cite{zhao2021spectral}.


The main objective of the present work is to combine the bi-orthogonal decomposition method with PINNs to establish an effective algorithm for solving the time-dependent SFPDEs. A discretization scheme for SFPDEs based on the bi-orthogonal decomposition method was proposed in \cite{zhao2021spectral}. However, a starting stage is required, that is, a hybrid algorithm based on QMC or PC needs to be employed for the problem whose initial conditions are deterministic. The above hybrid algorithm will be unstable when the eigenvalues of the system cross or  are very closely. To this end, herein we propose a new hybrid PINN method based on bi-orthogonal decomposition, which we call the BO-fPINN algorithm. Specifically, we construct a surrogate model of deep neural network based on the generalized KL expansion of the solution of the time-dependent SFPDEs, and then design the loss function through the BO conditions and the weak form of the equation.  In particular, since the back-propagation algorithm cannot compute the fractional derivative automatically, a discrete scheme based on the Gr\"{u}nwald-Letnikov (GL) formula is used for computing the fractional derivatives of the neural network. Numerical experiments show that the BO-fPINN method can solve problems with eigenvalue crossings. In addition, the same computational framework can be used to solve data-driven inverse problems associated with SFPDEs as well as problems with noisy initial data. Our diverse numerical examples demonstrate the efficiency and robustness of the BO-fPINN algorithm for a range of conditions, including different fractional orders. We also present a transfer learning example to demonstrate  that transfer learning can speed up the training of the BO-fPINN method for solving SFPDEs.

The paper is organized as follows. In Section 2, we set up the time-dependent SFPDEs. In Section 3, we build our BO-fPINN method after the introduction of the BO decomposition and its weak formulation, as well as the review of fPINN method. In Section 4, we demonstrate the accuracy and performance of our BO-fPINN method with several numerical experiments. We conclude the paper and provide a brief discussion in Section 5. In the appendices, we include definitions as well as additional results.

\section{Problem setup}

Let $(\Omega, \mathcal{F}, P)$ be a probability space, where $\Omega$ is the sample space, $\mathcal{F}$ is the $\sigma$-algebra of subsets of $\Omega$, and $P$ is a probability measure. We consider the following time-dependent SFPDE:
\begin{equation}\label{eqn:SFPDE}
\frac{\partial u(x,t;\omega)}{\partial t} =D_{|x_1|}^{\alpha}u(x,t;\omega) + D_{|x_2|}^{\beta}u(x,t;\omega)+ K(x;\omega)f(u) + g(x,t;\omega), \quad (x,t)\in D\times (0,T], \quad \omega\in \Omega,
\end{equation}
with initial and boundary conditions
\begin{equation}
\label{eqn:SFPDEIB}
\begin{aligned}
&u(x,0;\omega)  = u_0(x;\omega), \qquad x\in \mathcal{D},\omega\in \Omega, \\
 &u(x,t;\omega) = 0, \qquad x\in \partial \mathcal{D},
    \end{aligned}
\end{equation}
in which $\alpha, \beta\in (1,2)$,  $x=(x_1,x_2)$, $\mathcal{D}=(a,b)\times (c,d)$ with boundary denoted by $\partial \mathcal{D}$, $f(u)$ is a nonlinear operator. Here, $K(x;\omega)$ and $g(x,t;\omega)$ are the sources of randomness, which can be modeled as a set of random variables or as stochastic fields.   $D_{|x_1|}^{\alpha}$ and $D_{|x_2|}^{\beta}$ are the Riesz fractional operators defined by
\begin{equation}
\begin{aligned}
&D_{|x_1|}^{\alpha}u=-\frac{1}{2\cos(\pi \alpha/2)}(\sideset{}{_{a,{x_1}}^{\alpha}}{\mathop{D}}u
+\sideset{}{_{{x_1},b}^{\alpha}}{\mathop{D}}u), \\
 &D_{|x_2|}^{\beta}u=-\frac{1}{2\cos(\pi \beta/2)}(\sideset{}{_{c,{x_2}}^{\beta}}{\mathop{D}}u
+\sideset{}{_{{x_2},d}^{\beta}}{\mathop{D}}u),
    \end{aligned}
\end{equation}
where $\sideset{}{_{a,{x_1}}^{\alpha}}{\mathop{D}}$, $\sideset{}{_{{x_1},b}^{\alpha}}{\mathop{D}}$, $\sideset{}{_{c,{x_2}}^{\beta}}{\mathop{D}}$, and $\sideset{}{_{{x_2},d}^{\beta}}{\mathop{D}}$ are the left and right Riemann-Liouville fractional operators, the definitions of which are given in \ref{appendix:fractional-derivative-def}. We aim to evaluate the mean and standard deviation of the solution $u(x,t;\omega)$ and thus quantify the uncertainty propagation along with the random inputs.

\section{Methodology}
\subsection{The BO decomposition method and its weak formulation}\label{S:3-1}
We consider the following finite  approximation of a random field $u(x,t;\omega)$:
\begin{equation}\label{eqn:BOexpansion}
   u(x,t;\omega)\approx u_N(x,t;\omega) = \overbar{u}(x,t) + \sum_{i=1}^{N}u_i(x,t)Y_{i}(t;\omega),
\end{equation}
which is referred to as the bi-orthonormal decomposition and was first proposed in ~\cite{cheng2013dynamically1}. In the above expression, $\overbar{u}(x,t)$ is the deterministic mean field function, $u_i(x,t)$, $i=1,\cdots, N$ are  a set of orthonormal spatial basis functions:
\begin{align}\label{eqn:bo-condition1}
   \langle u_i(\cdot,t), u_j(\cdot,t)\rangle=\lambda_i(t)\delta_{ij},
\end{align}
where $\delta_{ij}$ is the Dirac¡¯s delta
function and $\lambda_i$s are the eigenvalues of the covariance kernel
$$C_{u(x,t)u(y,t)} = \E \left[ (u(x,t;\omega)-\overbar{u}(x,t))^{\mathrm{T}} (u(y,t;\omega)-\overbar{u}(y,t)) \right].$$
$Y_i(t;\omega)$, $i=1,\cdots, N$ are a set of orthonormal stochastic basis functions defined by
\begin{align}\label{eqn:bo-condition2}
   \mathbb{E}[Y_iY_j](t)=\delta_{ij},
\end{align}
which have zero mean i.e., $\E[Y_i(t;\omega)] = 0, i=1,\cdots, N$.

The conditions (\ref{eqn:bo-condition1}) and (\ref{eqn:bo-condition2}) guarantee both the spatial basis functions and stochastic coefficients are orthogonal in time. It is proved in ~\cite{cheng2013dynamically1} that these BO conditions lead to a set of independent and explicit evolution equations for all of the unknown quantities in (\ref{eqn:BOexpansion}). We state the BO evolution equations for the SFPDE (\ref{eqn:SFPDE}) without proof in \ref{appendix:QMC-BO}. To cope with the deterministic initial conditions of the SFPDE, a starting stage is required, that is, a
hybrid algorithm based on Quasi-Monte Carlo Methods needs to be employed, coined as QMC-BO method. For the implementation details of the QMC-BO algorithm, please see  Algorithm \ref{alg:QMC-BO} in \ref{appendix:QMC-BO}. The above hybrid algorithm will be unstable when the eigenvalues of the system cross or are very closely or perform poorly when the covariance matrix is singular or near singular. To overcome these difficulties, we employ the weak formulation of the BO methods, first given in \cite{zhang2020learning}, to  propose a new hybrid PINN method based on bi-orthogonal decomposition.

Now we derive the weak formulation of the original SFPDE. By applying the expectation  operator $\E[\cdot]$ on both sides of the SFPDE (\ref{eqn:SFPDE}) and replace $u$ by the finite expansion in \eqref{eqn:BOexpansion}. Noticing that $\E[Y_i]=0$, $i=1,\cdot,\cdot,\cdot,N$, we obtain the evolution equation for the mean of the solution:
    \begin{equation}\label{eqn:weakmean}
        \E \left[ \pdv{u}{t} \right] = \pdv{\overbar{u}}{t} = \E \left[ \mathcal{N}_x^{\alpha,\beta}[u(x,t;\omega)]+g(x,t;\omega) \right],
        \end{equation}
where, $\mathcal{N}_{|x|}^{\alpha,\beta}[u(x,t;\omega)]:=D_{|x_1|}^{\alpha}u(x,t;\omega) + D_{|x_2|}^{\beta}u(x,t;\omega)+ K(x;\omega)f(u)$ for notation simplification.

Applying the operator (inner product wrt $x$) $ \langle  \cdot, u_i \rangle$ on both sides of the SFPDE, we have:
    \begin{equation}\label{eqn:weakui}
        \left\langle \pdv{u}{t}, u_i \right\rangle = \left\langle \mathcal{N}_x^{\alpha,\beta}[u(x,t;\omega)]+g(x,t;\omega), u_i \right\rangle.
    \end{equation}

By multiplying $Y_j$ on both sides of the equation (\ref{eqn:SFPDE}) and applying the expectation operator $\E[\cdot Y_i]$ we obtain:
    \begin{equation}\label{eqn:weakyi}
        \E \left[ \pdv{u}{t} Y_i\right] = \E \left[ \left( \mathcal{N}_x^{\alpha,\beta}[u(x,t;\omega)]+g(x,t;\omega)\right) Y_i \right]
    \end{equation}

\subsection{fPINN for one-dimensional fractional PDEs}
In this section, we review the fractional physics-informed neural networks (fPINN) methods \cite{pang2019fpinns} by solving the following one dimensional deterministic FPDE:
\begin{equation}\label{eqn:prob}
\begin{aligned}
&\frac{\partial u(x,t)}{\partial t}-D_{|x|}^{\alpha}u(x,t) = g(x,t), \quad (x,t)\in (0,1)\times (0,T]\\
&u(x,0)  = v(x), \\
&u(0,t)  =u(1,t)=0.
    \end{aligned}
\end{equation}

\begin{figure}[!ht]
	\centering	\includegraphics[width=0.96\textwidth]{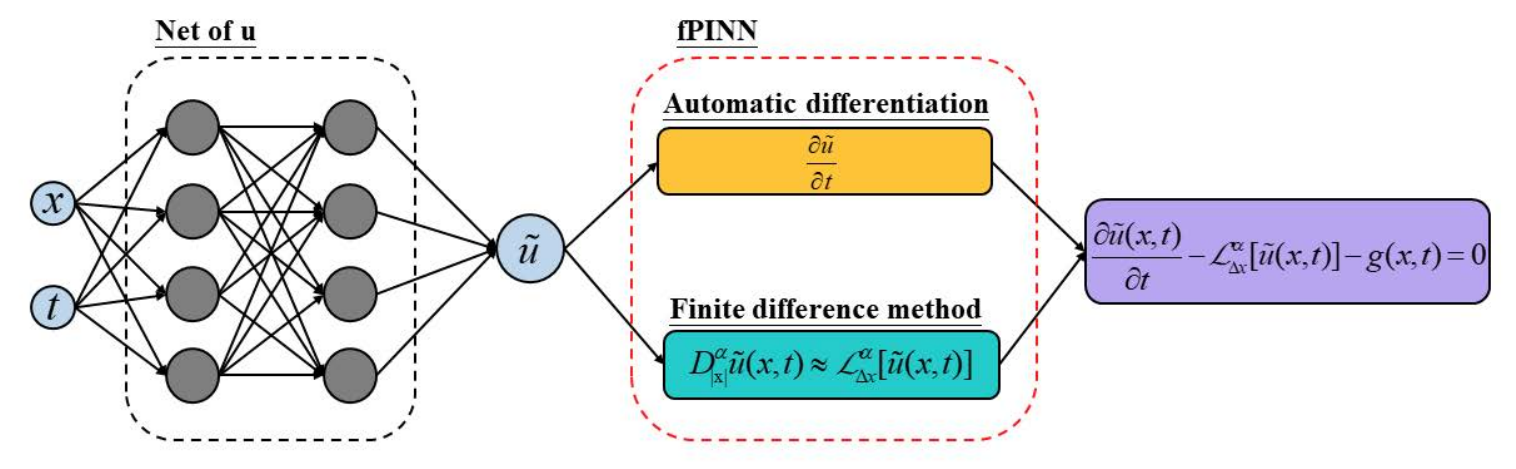}
	\caption{Schematic of the fPINN method for solving fractional partial differential equations.}\label{fig:fPINN}
\end{figure}

The approximate solution of $u(x,t)$ is constructed as $\tilde{u}(x,t)=x(1-x)u_{nn}(x,t;\theta)$ such that it satisfies the boundary conditions automatically. Here, $u_{nn}(x,t;\theta)$ is a deep neural network (DNN) surrogate parameterized by $\theta$, which takes the coordinate $x$ and $t$ as the input and outputs a vector that has the same dimension as $u$. The temporal derivative of $u_{nn}(x,t;\theta)$ can be calculated directly by automatic differentiation, however, we need to discretize the fractional spatial derivatives with the finite difference methods (FDM) firstly and then incorporate the discrete schemes into the loss function of PINNs. For simplicity, we consider the following shifted GL formula as in \cite{pang2019fpinns}
\begin{align*}
  D_{|x|}^{\alpha}\tilde{u}(x,t;\theta) & \approx
  \mathcal{L}_{\Delta x}^{\alpha}[\tilde{u}(x,t;\theta)]
   =  c_{\alpha} \left[ \frac{\alpha}{2}\delta_{\Delta x,1}^{\alpha}\tilde{u}(x,t;\theta)+(1-\frac{\alpha}{2}) \delta_{\Delta x,0}^{\alpha}\tilde{u}(x,t;\theta)  \right] ,
\end{align*}
where~$c_{\alpha}=-\frac{1}{2\cos(\frac{\pi\alpha}{2})}$, and the definition of $\delta_{\Delta x,p}^{\alpha}, p=0,1$ can be found in \ref{appendix:GL}.
Given the training dataset $\{(x_g^i,t_g^i, g(x_g^i,t_g^i))\}_{i=1}^{N_g}$ and $\{(x_v^i,v(x_v^i))\}_{i=1}^{N_v}$, we can optimize the parameters $\theta$ in the neural network by minimizing the loss function as follows:
\begin{equation}
    \mathcal{LOSS}_{fPINN}(\theta)=\frac{1}{N_v}\sum_{i=1}^{N_v}[\tilde{u}(x_v^i, 0;\theta)-v(x_v^i)]^2 + \frac{1}{N_g}\sum_{i=1}^{N_g}\bigg
    \{\frac{\partial \tilde{ u}(x_g^i,t_g^i;\theta)}{\partial t}-\mathcal{L}_{\Delta x}^{\alpha}[\tilde{\boldsymbol u}(x_g^i,t_g^i;\theta)]- g(x_g^i,t_g^i)\bigg\}^2.
  \end{equation}
Fig. \ref{fig:fPINN} shows a sketch of fPINN, while the workflow of solving fractional differential equation with fPINN can be summarized in Algorithm \ref{alg:PINN}. We use a toy example to show the performance of fPINN compared with the FDM scheme for fractional PDEs in \ref{appendix:pinn-gl}.

\begin{algorithm}[H]
\caption{fPINN for solving deterministic fractional PDEs}
\label{alg:PINN}
\textbf{Step 1:} Specify the training set:
  \begin{equation*}
    \text{training data for initial condition: }\{(x_v^i,v(x_v^i)\}_{i=1}^{N_v}, \quad \text{training data for equation: }\{(x_g^i,t_g^i, g(x_g^i,t_g^i))\}_{i=1}^{N_g}.
  \end{equation*}
\textbf{Step 2:} Construct a DNN $\tilde{u}(x,t;\theta)$ with random initialized parameters $\theta$.\\
\textbf{Step 3:} Specify a loss function by summing the mean squared error of both the $v$ observations and the equation:
\begin{equation*}
  \mathcal{LOSS}_{fPINN}(\theta)=\frac{1}{N_v}\sum_{i=1}^{N_v}[\tilde{u}(x_v^i, 0;\theta)-v(x_v^i)]^2 + \frac{1}{N_g}\sum_{i=1}^{N_g}\bigg
    \{\frac{\partial \tilde{ u}(x_g^i,t_g^i;\theta)}{\partial t}-\mathcal{L}_{\Delta x}^{\alpha}[\tilde{\boldsymbol u}(x_g^i,t_g^i;\theta)]- g(x_g^i,t_g^i)\bigg\}^2.
    \end{equation*}
\textbf{Step 4:} Train the DNN to find the optimal parameters $\theta^*$ by minimizing the loss function
  \begin{equation}
    \theta^*=\arg\min \mathcal{LOSS}_{fPINN}(\theta).
  \end{equation}
\end{algorithm}

\subsection{BO-fPINNs for solving time-dependent fractional PDEs with random inputs}
In this section we build up the framework of BO-fPINNs based on the weak formulation of the stochastic FPDE (\ref{eqn:weakmean})-(\ref{eqn:weakyi}). Suppose we can parameterize the original SFPDE into a PDE with finite dimensional random variables $\xi(\omega)$, then $u(x,t;\omega )$ can be written as $u(x,t;\xi)$ and we can rewrite the finite  approximation of a random field $u(x,t;\xi)$ (\ref{eqn:BOexpansion}) as
\begin{align}\label{eqn:BO-nn}
  u_N(x,t;\xi)=\overline{u}(x,t)+\sum_{i=1}^{N}a_i(t)u_i(x,t)Y_i(t;\xi),
\end{align}
while enforcing $\langle u_i,u_i\rangle=1$, and $E[Y_i^2]=1$. The time-dependent coefficients $a_i(t)$ are scaling factors and play the role of $\sqrt{\lambda_i}$ in the standard KL expansion. We construct four separate neural networks for each function in the above expansion. Specifically, $\overline{u}_{nn}(x,t)$ takes $x$ and $t$ as inputs and outputs $\mathbb{E}[u(x,t;\omega)]$; $A_{nn}(t)$ is the neural net that outputs an $N$-dimensional vector for approximating $a_i$ for $i=1,...,N$; $U_{nn}(x,t)$ is the neural net with $N$-dimensional vector outputs representing $u_i(x,t)$ for $i=1,...,N$; and $Y_{nn}(\xi, t)$ takes $\xi$ and $t$ as inputs and outputs an $N$-dimensional vector representing $Y_i(t;\xi)$ for $i=1,...,N$. By substituting them into (\ref{eqn:BO-nn}) we obtain a neural net surrogate for the SFPDE solution $u(x,t;\omega)$, and
then, we can obtain a surrogate of the solution:
\begin{align}\label{nn2}
  u_{nn}(x,t;\omega)=\overline{u}_{nn}+\sum_{i=1}^{N}A_{nn,i}U_{nn,i}Y_{nn,i}.
\end{align}

We notice that the integer derivatives of
the quantities of interest can be easily
obtained via the automatic differentiation algorithm, while for the fractional derivatives with respect to $x$, we can use the GL formula to discretize them. We assume that we have $n_x$ training points $\{x^k\}_{k=1}^{n_x}$  in the physical domain, $n_{\xi}$ training points
$\{\xi^l\}_{l=1}^{n_{\xi}}$ in the probabilistic space, and $n_t$ uniformly sampled points $\{t^s\}_{s=1}^{n_t}$ in the time domain $[0,T]$. Then, the loss function of the BO-fPINN method can be defined as a weighted summation of the weak formulation of SFPDE, initial/boundary conditions,
BO constraints and the additional regularization terms, i.e.,
\begin{align}\label{eqn:Loss}
  \mathcal{LOSS}(\theta)= \lambda_wMSE_w + \lambda_{IC}MSE_{IC}+
  \lambda_{BC}MSE_{BC}+\lambda_{BO}MSE_{BO}+\lambda_gMSE_g,
\end{align}
where $\lambda_w$,$\lambda_{IC}$,$\lambda_{BC}$, $\lambda_{BO}$ and $\lambda_g$ are objective function weights for balancing the various term in (\ref{eqn:Loss}). Consequently, $MSE_{w}$ is the loss function for the weak formulation,  $MSE_{IC}$ and $MSE_{BC}$ corresponds to the loss on the initial data/boundary conditions respectively,
$MSE_{BO}$ enforces the bi-orthogonal conditions, and $MSE_g$ is a penalty term, which may help speed up the training and alleviate over-fitting. The explicit forms of the the different loss parts are given below.

By employing a numerical quadrature rule we can calculate the integration in the weak formulation and obtain
\begin{align*}
  MSE_w=\frac{1}{n_xn_t}\sum_{k,s}(\epsilon_1^{ks})^2 +\frac{1}{n_tn_{\xi}}\sum_{s,l}(\epsilon_2^{sl})^2+\frac{1}{n_xn_t}\sum_{k,s}(\epsilon_3^{ks})^2,
\end{align*}
where
\begin{align*}
  \epsilon_1^{ks}:=&\mathbb{E}\bigg [\frac{\partial u_{nn}}{\partial t}-
  \mathcal{L}_{\Delta x}^{\alpha}u_{nn}(x^k,t^s;\xi)-g(x^k,t^s;\xi)\bigg ], \\
  \epsilon_2^{sl}:=&\bigg\langle \frac{\partial u_{nn}}{\partial t}-\mathcal{L}_{\Delta x}^{\alpha}u_{nn}(x,t^s;\xi_c^l)-g(x,t^s;\xi^ll), U_{nn,i}(x,t^s)\bigg\rangle, \\
  \epsilon_3^{ks}:=&E\bigg [(\frac{\partial u_{nn}}{\partial t}-
  \mathcal{L}_{\Delta x}^{\alpha}u_{nn}(x^k,t^s;\xi)-g(x^k,t^s;\xi))Y_{nn,i}(t^s,\xi)\bigg],
\end{align*}
and the fractional operator in $\mathcal{L}_{\Delta x}^{\alpha}$ is discretized by FDM.

The loss corresponding to the initial condition is calculated by
\begin{align*}
  MSE_{IC}= &\frac{1}{n_x}\sum_{k=1}^{n_x}(\overline{u}_{nn}(x^k,t_0)-\overline{u}(x^k,t_0))^2+\frac{1}{Nn_x}\sum_{i=1}^{N}\sum_{k=1}^{n_x}(U_{nn,i}(x^k,t_0)-u_i(x^k,t_0))^2\\
  &+\frac{1}{N}\sum_{i=1}^{N}(A_{nn,i}(t_0)-a_i(t_0))^2+\frac{1}{Nn_{\xi}}\sum_{i=1}^{N}\sum_{l=1}^{n_{\xi}}(Y_{nn,i}(t_0;\xi^l)-Y_i(t_0;\xi^l))^2,
\end{align*}
where
\begin{equation}
\begin{aligned}
&\overline{u}(x,t_{0}) = \mathbb{E}[u(x,t_{0};\xi)],\\
&u_{i}(x,t_{0}) = v_{i}(x),\\
&a_{i}(t_{0}) = \sqrt{\mathbb{E}\bigg[\langle u(x,t_{0};\xi) - \mathbb{E}[u(x,t_{0};\xi)], v_{i} \rangle^2 \bigg]},\\
&Y_{i}(t_{0};\xi) = \frac{1}{a_{i}(t_{0})} \langle u(x,t_{0};\xi) - \mathbb{E}[u(x,t_{0};\xi)], v_{i} \rangle.
\end{aligned}
\end{equation}
For the deterministic initial condition, we can set $u_i(x,t_0)$ to be orthonormal bases satisfying the boundary conditions, $Y_i(t_0;\xi)$ to be the generalized polynomial chaos base of $\xi$ with unit variance, and $a_i(t_0)$ to be $0$ \cite{zhang2020learning}.

By taking the weak formulation of \eqref{eqn:SFPDEIB} in the random space, we can obtain the loss corresponding to the boundary condition

\begin{align*}
  MSE_{BC}= \frac{1}{n_t}\sum_{s=1}^{n_t}(\overline{u}_{nn}(x_b,t^s))^2+\frac{1}{Nn_t}\sum_{i=1}^{N}\sum_{s=1}^{n_t}\bigg (A_{nn,i}(t^s)U_{nn,i}(x_b,t^s)\bigg )^2.
\end{align*}

By enforcing the BO condition, we can get the loss function for the BO constraints~\cite{zhang2020learning}:

\begin{align*}
  MSE_{BO}&=\frac{1}{Nn_t}\sum_{i=1}^{N}\sum_{s=1}^{n_t}(\mathbb{E}[Y_i(t^s;\xi)])^2
  +\frac{1}{N^2n_t}\sum_{i=1}^{N}\sum_{j=1}^{N}\sum_{s=1}^{n_t}\bigg(\bigg\langle \frac{dU_{nn,i}(x,t^s)}{dt}, U_{nn,j}(x,t^s)\bigg\rangle \\
  &+\bigg\langle \frac{dU_{nn,j}(x,t^s)}{dt}, U_{nn,i}\bigg\rangle\bigg)^2 + \frac{1}{Nn_t}\sum_{i=1}^{N}\sum_{j=1}^{N}\sum_{s=1}^{n_t}\bigg(\mathbb{E}\bigg [Y_{nn,i}(t^s;\xi)\frac{dY_{nn,j}(t^s;\xi)}{dt}\bigg] \\
  &+ \mathbb{E}\bigg[Y_{nn,j}(t^s;\xi)\frac{dY_{nn,i}(t^s;\xi)}{dt}\bigg ]\bigg)^2.
\end{align*}
Finally, the loss from the original equation is
\begin{align*}
 MSE_g=\frac{1}{n_xn_tn_{\xi}}\sum_{k=1}^{n_x}\sum_{s=1}^{n_t}\sum_{l=1}^{n_{\xi}}\left(\frac{\partial u_{nn}}{\partial t}
 -\mathcal{L}_{\Delta x}^{\alpha}u_{nn}(x^k,t^s;\xi^l)-g(x^k,t^s;\xi^l)\right)^2.
\end{align*}

\begin{figure}[!ht]
	\centering	\includegraphics[width=0.9\textwidth]{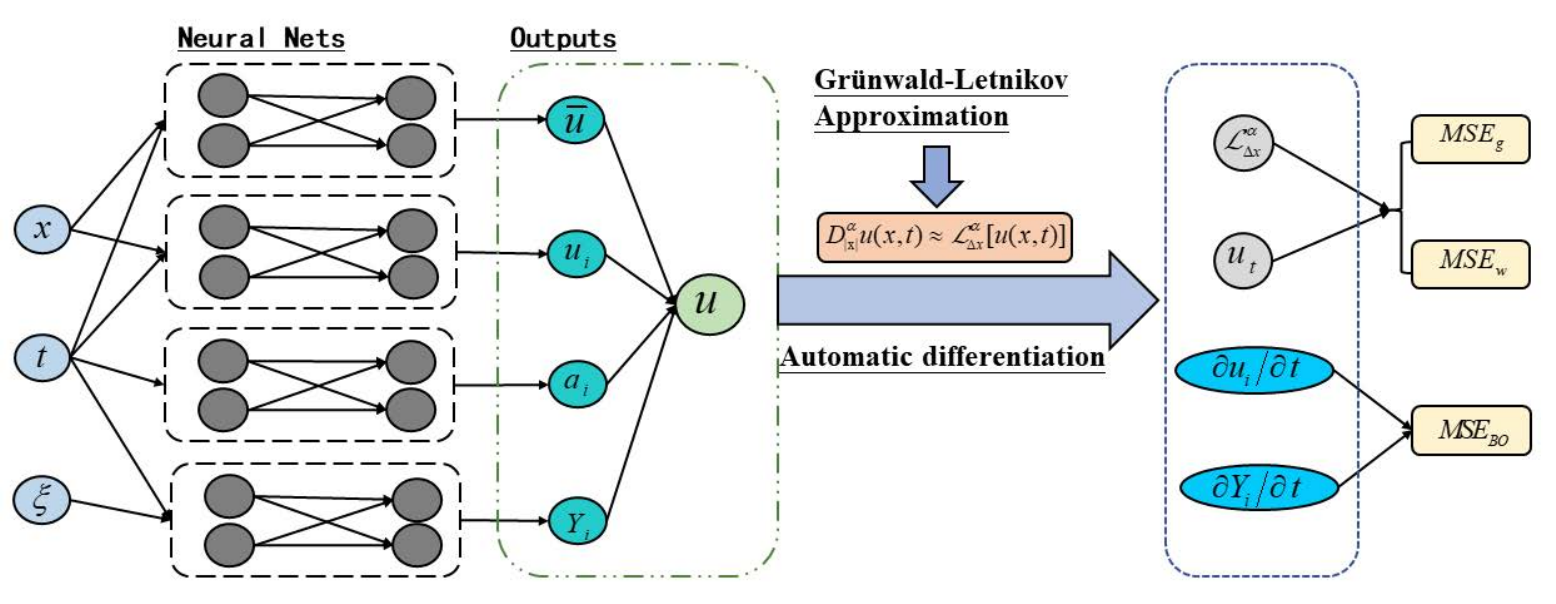}
	\caption{Schematic of the BO-fPINN method for solving time-dependent stochastic fractional partial differential equations. }\label{fig:fPINN-BO}
\end{figure}

A sketch of the computation
graph for the BO-fPINN is shown in Figure \ref{fig:fPINN-BO}, while the algorithm is summarized in Algorithm \ref{alg:fPINN-BO}. We remark that the BO constraints are written into the loss functions instead of deriving explicit expressions as in the classical BO method. Thus, the new BO-fPINN method can overcome the bottleneck of the BO method, which will be verified by the numerical examples in the next section.

\begin{algorithm}
\setlength{\baselineskip}{20pt}
\caption{~BO-fPINN~ method for solving time-dependent stochastic fractional PDEs}\label{alg:fPINN-BO}
\textbf{Step 1:} Given $n_x$ equidistantly distributed training points $\{x^k\}_{k=1}^{n_x}$ in physical domain, $n_t$ uniformly distributed training points $\{t^s\}_{k=1}^{n_t}$ in the time domain, $n_{\xi}$ training points $\{\xi^l\}_{l=1}^{n_{\xi}}$ in the stochastic  space.\\
\textbf{Step 2:} Build four neural networks:~ $$\overline{u}_{nn}(x,t),~A_{nn}(t),~U_{nn}(x,t),~Y_{nn}(t;\xi).$$
\textbf{Step 3:} Compute loss function \eqref{eqn:Loss} by using GL formulas and quadrature rules. \\
\textbf{Step 4:} Train the neural networks by minimizing the loss function.\\
\textbf{Step 5:} Predict the SFPDE solution using \eqref{nn2}.
\end{algorithm}

\section{Results}
In this section, we test the proposed BO-fPINN method for solving different time-dependent SFPDEs. We first test our BO-fPINN method with a benchmark case that is especially designed to have exact solutions for the BO
representations, where the eigenvalues have crossings in the time domain. To demonstrate the advantage of the BO-fPINN method over
the standard methods, we then solve a nonlinear stochastic fractional diffusion-reaction equation with high-dimensional random inputs. An inverse problem is also considered to demonstrate the new capability of the proposed BO-fPINN method. We also show that transfer learning can speed up the
training of the BO-fPINN method for solving SFPDEs. Finally, a 2D nonlinear fractional diffusion-reaction equation driven by random forcing is considered to further show the effectiveness of the proposed method. For all
test cases, we use deep feed-forward neural networks for ~$\overline{u}_{nn}(x,t),~A_{nn}(t),~U_{nn}(x,t)$~ and ~$Y_{nn}(t;\xi)$,~and the sizes of the hidden layers of neural networks are shown in Table \ref{tab:NN}. In the following examples, the predictive accuracy of the trained model is assessed using the $L_2$ error and the relative $L_2$ error, defined by $\|f_{NN}-f_{exact}\|_2$ and $\|f_{NN}-f_{exact}\|_2 / \|f_{exact}\|_2$, respectively, for any function $f$.


\subsection{1D nonlinear stochastic fractional reaction-diffusion equation with manufactured solution}\label{sec:RD}
First, we demonstrate that the BO-fPINN method is robust for time-dependent SFPDEs for the eigenvalue crossing phenomenon, where the standard BO method would fail \cite{choi2014equivalence,babaee2017robust}.

We consider the following nonlinear reaction-diffusion equation
\begin{equation}\label{prob:1d-cross}
      u_t-D_{|x|}^\alpha u-u(1-u^2)=g(x,t;\xi_1,\xi_2), \quad  x\in \mathcal{D}=[0,1],\quad t\in [0,\pi],
\end{equation}
 where $\xi_1(\omega)$ and $\xi_2(\omega)$ are two identically independent uniformly distributed random variables in $[0,1]$. A manufactured solution $u(x,t;\xi_1,\xi_2)$ with exact BO components can be constructed as
\begin{equation}
\label{example:cross}
\begin{aligned}
  u(x,t;\xi_1,\xi_2) & = 100\sin(\frac{t}{2}+\frac{\pi}{4})x^3(1-x)^3-\sqrt{3}(1.5+\sin(t))\sin(\pi x)(2\xi_1-1) \\
   & +\sqrt{3}(1.5+\cos(3t))\sin(2\pi x)(2\xi_2-1).
\end{aligned}
\end{equation}
The random forcing term $g(x,t;\omega)$ can be calculated directly and the explicit formula is given in \ref{appendix:g}. By normalizing the bases and the random coefficients, the BO expansion yields the same expression in the form of (\ref{eqn:BO-nn}):
\begin{equation}
\label{cross:bo}
\begin{aligned}
  &u_1(x,t)=-\sqrt{2}\sin(\pi x), & u_2(x,t)=\sqrt{2}\sin(2\pi x),\\
 & a_1(t)=\frac{1}{\sqrt{2}}(1.5+\sin(t)), &a_2(t)=\frac{1}{\sqrt{2}}(1.5+\cos(3t)),\\
  & Y_1(t;\xi_1,\xi_2)=\sqrt{3}(2\xi_1-1), & Y_2(t;\xi_1,\xi_2)=\sqrt{3}(2\xi_2-1).
\end{aligned}
\end{equation}


 We set $\alpha=1.5$ and use the second order GL formula to approximate the space-fractional derivatives. The initial conditions are taken directly from (\ref{cross:bo}). The number of spatial training points is $n_x=70$ and the samples are equidistantly distributed in $[0,1]$. The number of temporal training points is $n_t=70$ and the samples are drawn from a uniform distribution. For the training points in the stochastic space, we use the eighth-order Gauss-Legendre quadrature rule for both $\xi_1$ and $\xi_2$, generating
64 points in the probabilistic space. In this example, in order to improve the accuracy and numerical stability, we fix $\lambda_w=1$ and the rest weight functions are learned dynamically \cite{wang2021understanding} with initial values $\lambda_{BO}=1000$, $\lambda_{IC}=1000$, $\lambda_{BC}=1000$ and $\lambda_{g}=1$ respectively. The neural networks are trained with an Adam optimizer with a base learning rate 0.001, which we decrease exponentially every 1000 iterations by a drop rate of 0.9. The total number of iterations 400000 epochs and then L-BFGS-B optimizer is implemented.

\begin{figure}[!ht]
	\centering	
    \subfloat[]{\includegraphics[width=0.42\textwidth]{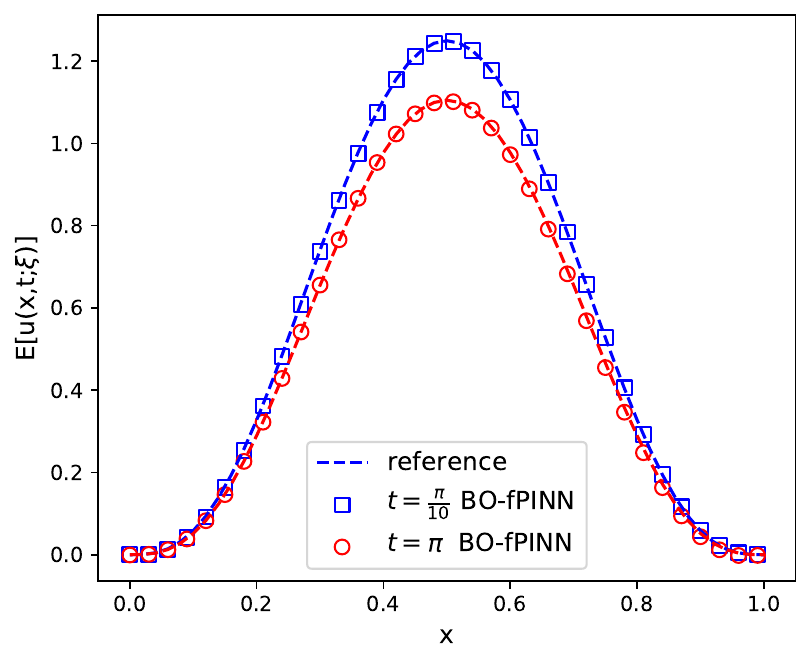}}\qquad
	\subfloat[]{\includegraphics[width=0.41\textwidth]{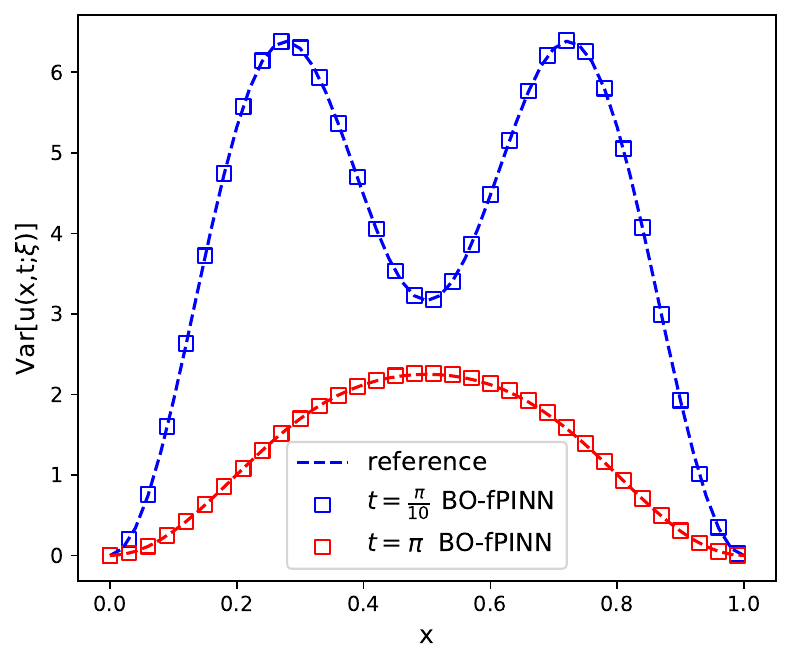}} \\
	\caption{1D stochastic fractional reaction-diffusion equation with manufactured solution. (a):~BO-fPINN solution mean at $t=\frac{\pi}{10}$ and $t=\pi$. (b):~BO-fPINN solution variance at $t=\frac{\pi}{10}$ and $t=\pi$. We use the exact solutions as reference. The scattered points denote the training points in the physical space.}\label{fig:cross-mean}
\end{figure}

\begin{figure}[!ht]
	\centering	
    \subfloat[]{\includegraphics[width=0.43\textwidth]{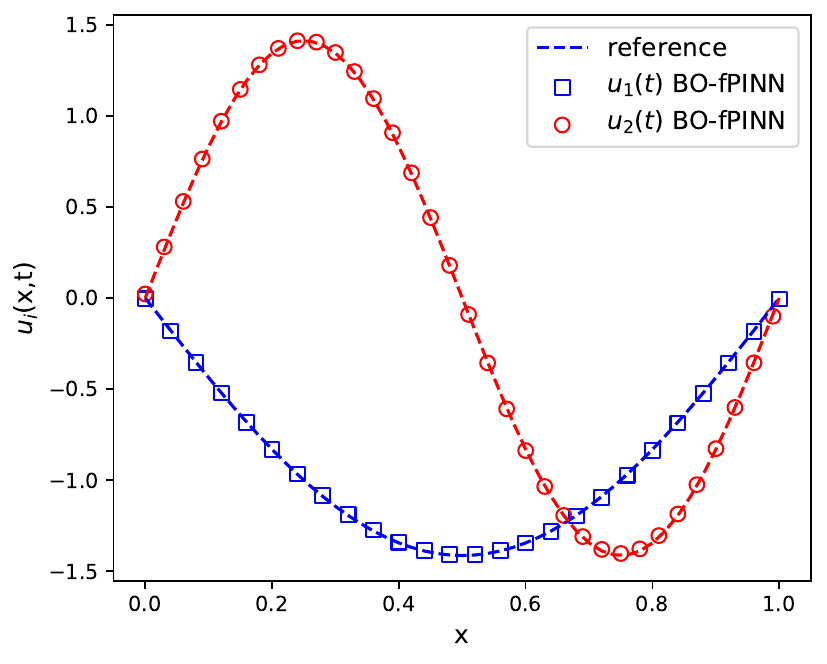}}\qquad
	\subfloat[]{\includegraphics[width=0.42\textwidth]{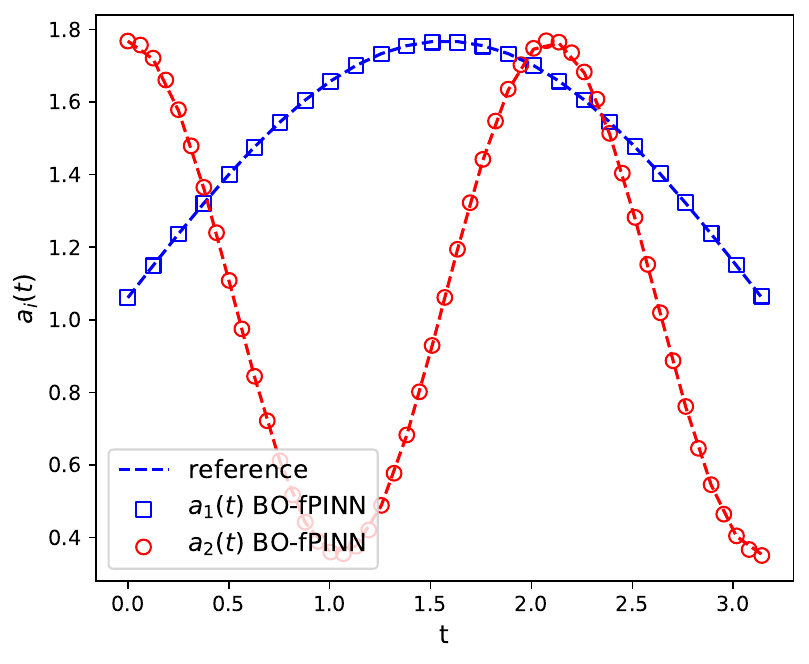}} \\	
	\caption{1D stochastic fractional reaction-diffusion equation with manufactured solution. (a):~The spatial bases $u_1$ and $u_2$ at $t=\pi$. (b):~The evolution of the scaling factors $a_i$, which correspond to the eigenvalues in the standard BO method and they have crossed in time . }\label{fig:cross-ui}
\end{figure}

In Figure \ref{fig:cross-mean}, we plot the BO-fPINN solution mean (Figure \ref{fig:cross-mean}a) and solution variance (Figure \ref{fig:cross-mean}b) versus the exact values at $t=\frac{\pi}{10}$ and $t=\pi$. In Figure \ref{fig:cross-ui}a, we compare the bases $u_i$ ($i=1,2$) at $t=\pi$ obtained from the BO-fPINN method to the exact solutions. We can see that the BO-fPINN solutions agree with the exact reference solutions very well. Figure \ref{fig:cross-ui}b shows the evolution of $a_i$ ($i=1,2$) with time, and we note that the scaling factors $a_i$ correspond to the eigenvalues in the standard BO method and they cross during the time domain. In this situation, the standard BO method
would fail while the proposed BO-fPINN method does not suffer from this issue. It is evident that the
BO-fPINN solutions agree with the exact reference solutions very well. In Table \ref{tab:error-cross}, we summarize the relative $L_2$ errors for all of the BO components at $t=\frac{\pi}{10}, \frac{\pi}{5}, \pi$, indicating the robust performance of the BO-fPINN method for problems with eigenvalue crossings.


\begin{table}[!htp]
\centering
\renewcommand\arraystretch{1.2}
\caption{1D Stochastic Fractional Reaction-Diffusion Equations with manufactured soultion. The relative $L_2$ errors of BO-fPINN solutions versus the exact solutions at $t=\frac{\pi}{10}, \frac{\pi}{5}, \pi$.
\label{tab:error-cross}}
\begin{tabular}{c|cccccc}
\hline
             &  $a_1$  &  $a_2$   &  $u_1$  &  $u_2$  & $Y_1$  &  $Y_2$ \\
\hline
$t=\frac{\pi}{10}$   & 0.014$\%$& 0.188$\%$ & 0.169$\%$&0.304$\%$&0.249$\%$&0.255$\%$  \\
\hline
$t=\frac{\pi}{5}$  & 0.092$\%$& 0.078$\%$ & 0.162$\%$&0.272$\%$&0.255$\%$&0.333$\%$  \\
\hline
$t=\pi$  & 0.426$\%$& 0.839$\%$ & 0.360$\%$&0.701$\%$&0.848$\%$&1.304$\%$  \\

\hline 
\end{tabular}
\end{table}

\subsection{1D nonlinear stochastic fractional reaction-diffusion equation with high-dimensional random inputs}
To further verify the advantage of the BO-fPINN method to predict the long-time behavior over the standard methods, we consider the following nonlinear diffusion-reaction equation with deterministic initial conditions and high-dimensional random inputs:
\begin{equation}\label{prob:1d}
  \begin{cases}
    u_t-\mu D_{|x|}^\alpha u-\epsilon f(u)=g(x,t;\omega), \quad  & x\in(-1,1),\quad t\in (0,T], \\
    u(x,0;\omega) = \sin(\pi x),\quad  &x\in(-1,1),\quad \omega\in\Omega,\\
    u(-1,t;\omega)=u(1,t;\omega)=0,&t\in(0,T],\quad \omega\in\Omega,
  \end{cases}
\end{equation}
where $f(u)=u(1-u^2)$, and $\mu$ and $\epsilon$ are time-independent diffusion and reaction coefficients, respectively. We consider four cases in this example: 1) forward problem with time independent random forcing $g(x;\omega)$; 2) forward problem with time-dependent random forcing $g(x,t;\omega)$; 3) inverse problem with time independent random forcing where the coefficients $\mu$ and $\epsilon$ are unknown but additional information for $u(x, t; \omega)$ is given; 4) transfer learning between different fractional orders.

\subsubsection{Forward problem with time independent random forcing}\label{sec:RD_notime}
We first consider the case with the source of randomness $g(x;\omega)=(1-x^2)h(x;\omega)$. The random process $h(x;\omega)$ is modeled as a Gaussian random field, i.e., $h(x;\omega) \sim \mathcal{GP}(1, C(x_1,x_2))$, where $C(x_1,x_2)$ is a squared exponential kernel with standard deviation $\sigma_h$ and correlation length $l_c$:
\begin{align*}
  C(x_1,x_2) =\sigma_h^2\exp\left(-\frac{(x_1-x_2)^2}{l_c^2}\right).
\end{align*}
For the random forcing $g(x;\omega)$, we set $\sigma_h=1$ and $l_c=0.1$, thus requiring 19 KL modes to capture at least $98\%$ of the fluctuation energy of $g(x;\omega)$. The diffusion coefficient is set to $\mu=1$ and the reaction coefficient is set to $\epsilon=1$. The neural networks used in the BO-fPINN method for this case can be found in \ref{appendix:network}. We should ponit out that $A_{nn}$ is composed of $N$ independent neural networks ($N$ is the number of BO expansion terms), each of which approximates one single scaling factor $a_i$. This is because we expect that $a_i$ may oscillate greatly in vastly different scales during the time evolution. The weighting coefficients of the loss function are taken as $\lambda_w=1$,~ $\lambda_{IC}=10$,~
$\lambda_{BC}=100$,~$\lambda_{BO}=10$,~$\lambda_{0}=0.01$. We use $n_x=71$ equidistantly distributed training points $\{x^k\}_{k=1}^{nx}$ in space, and $n_t=70$ uniformly distributed training points $\{t^k\}_{k=1}^{n_t}$ in the time domain, and $n_l=1000$ random samples $\{\xi^l\}_{l=1}^{n_{\xi}}$ in the 19-dimensional random space. The neural networks are trained with an Adam optimizer with learning rate 0.001 for 400000 epochs, followed by the L-BFGS-B optimizer.

We investigate the performance of BO-fPINN method using six BO expansion terms in (\ref{eqn:BO-nn}). To obtain the reference solution for the BO decomposition, we use the QMC-BO method (\ref{appendix:QMC-BO}) with switch time $t_s=0.01$ to solve the original BO equations due to the deterministic initial condition. The parameters for the QMC-BO method are $\Delta t=5\times 10^{-5}$, $\Delta x=0.02$, $M=1000$. A Quasi Monte Carlo method with 1000 samples is employed to solve the SFPDE to obtain the reference for the solution statistics.

\begin{figure}[!ht]
	\centering	
    \subfloat[]{\includegraphics[width=0.43\textwidth]{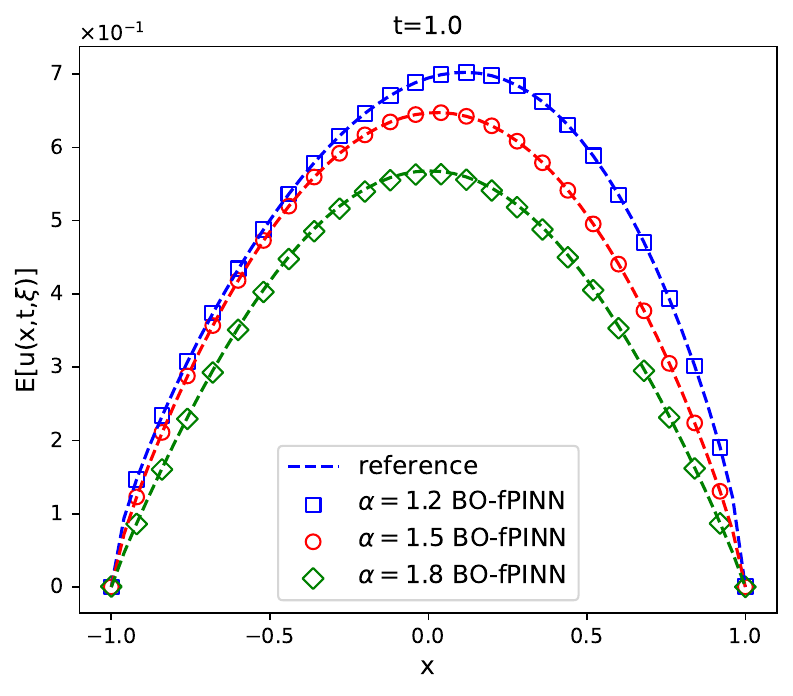}}\qquad
	\subfloat[]{\includegraphics[width=0.43\textwidth]{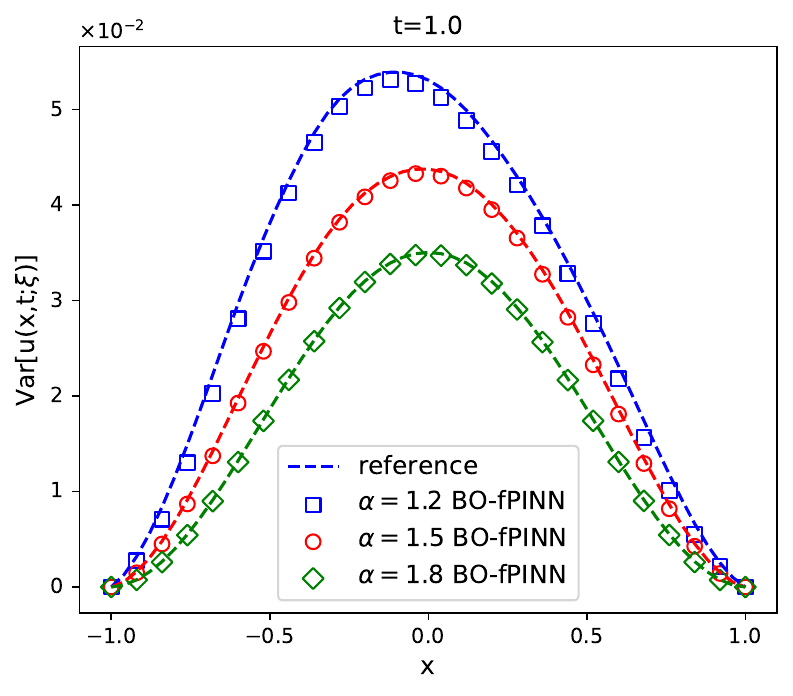}} \\
	\caption{Forward problem with time independent random forcing. (a):~BO-fPINN solution mean at $t=1.0$ with different values of the fractional order $\alpha$. (b):~BO-fPINN solution variance at $t=1.0$ with different values of the fractional order $\alpha$.
 The reference mean and variance  are calculated from the Quasi Monte Carlo simulation.}\label{fig:mean258}
\end{figure}

\begin{figure}[!ht]
	\centering	
    \subfloat[]{\includegraphics[width=0.43\textwidth]{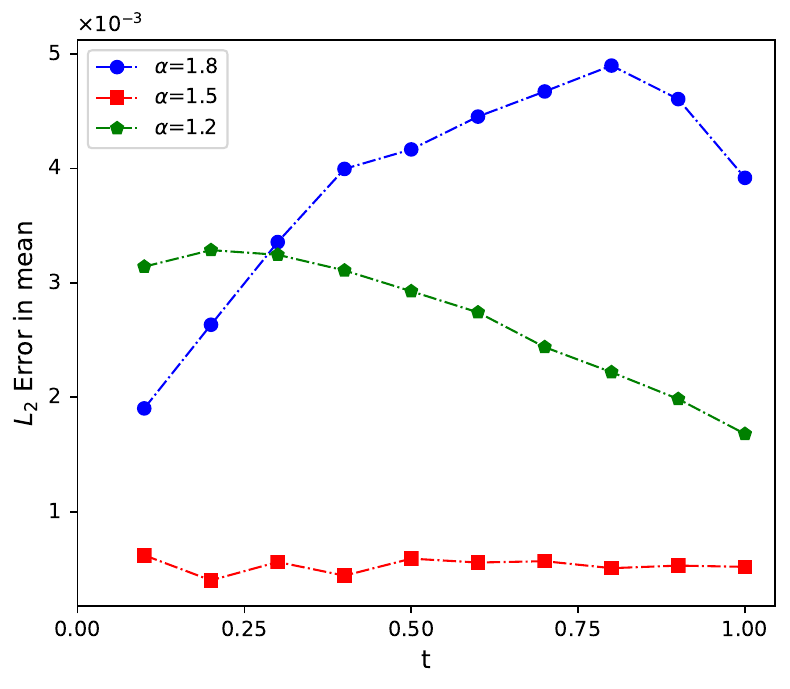}}\qquad
	\subfloat[]{\includegraphics[width=0.445\textwidth]{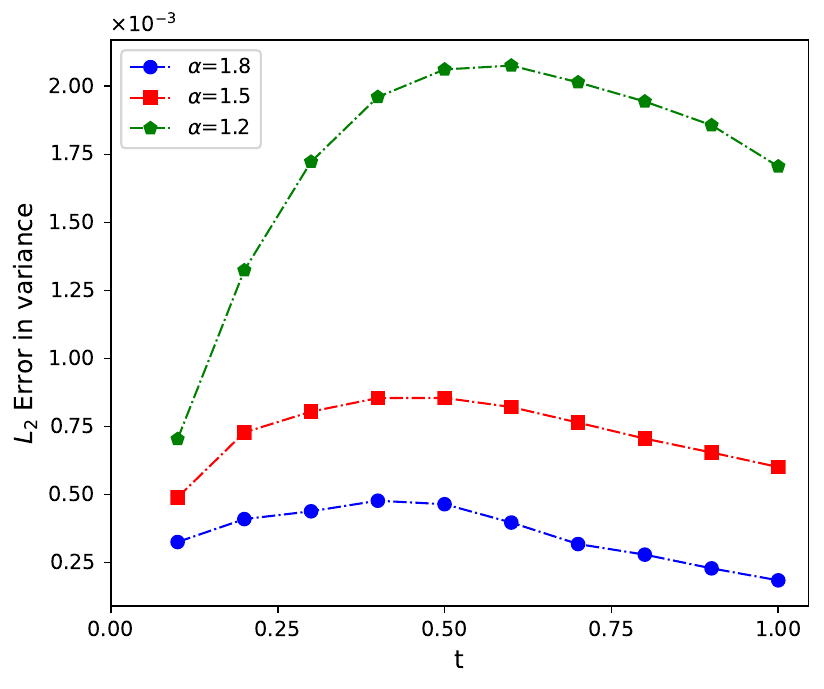}} \\
	\caption{Forward problem with time independent random forcing. (a):~The $L_2$ error in the mean obtained by BO-fPINN method with different values of the fractional order $\alpha$. (b):~The $L_2$ error in the variance obtained by BO-fPINN method with different values of the fractional order $\alpha$. We use the QMC solution as reference.
}\label{fig:l2-mean258}
\end{figure}

\begin{figure}[!ht]
	\centering	
    \subfloat[]{\includegraphics[width=0.42\textwidth]{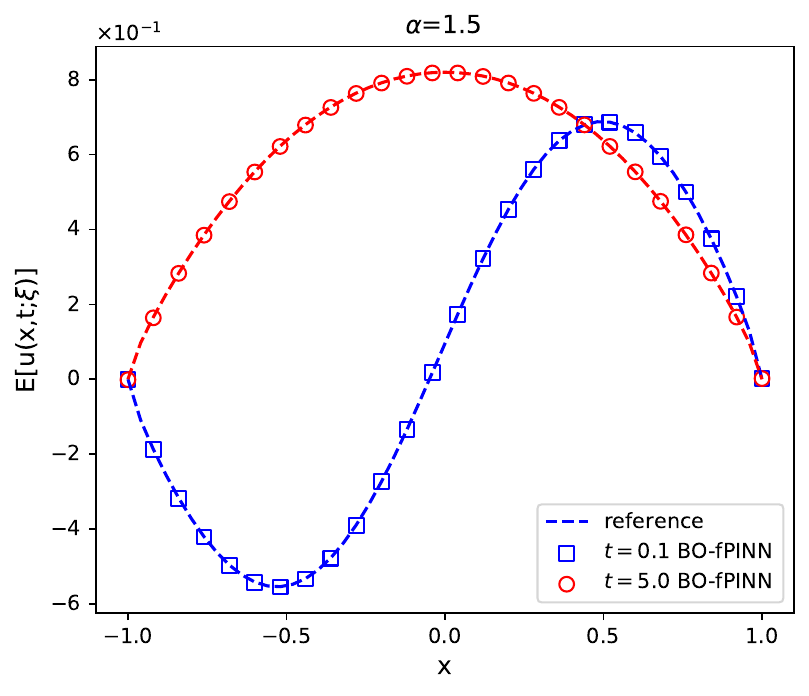}}\qquad
	\subfloat[]{\includegraphics[width=0.43\textwidth]{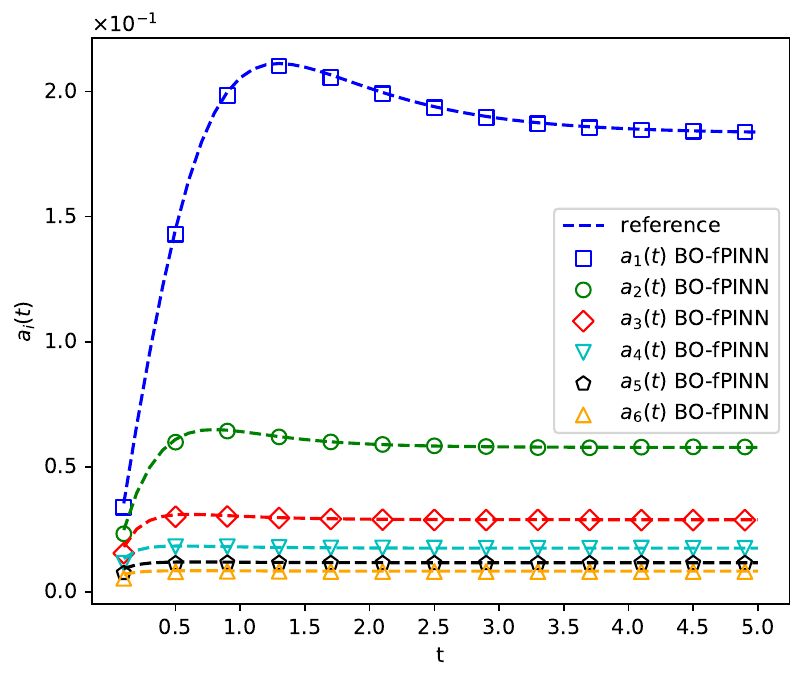}} \\	
	\caption{Forward problem with time independent random forcing. (a):~Solution mean at $t=0.1$ and $t=5.0$; the reference mean is calculated from the Quasi Monte Carlo simulation. (b):~Scaling factors $a_i$ at different time steps. The reference solutions are calculated using the QMC-BO method.}\label{fig:mean-1d-t5}
\end{figure}

\begin{figure}[!ht]
\centering{\includegraphics[height=0.36\textheight, width=0.85\textwidth]{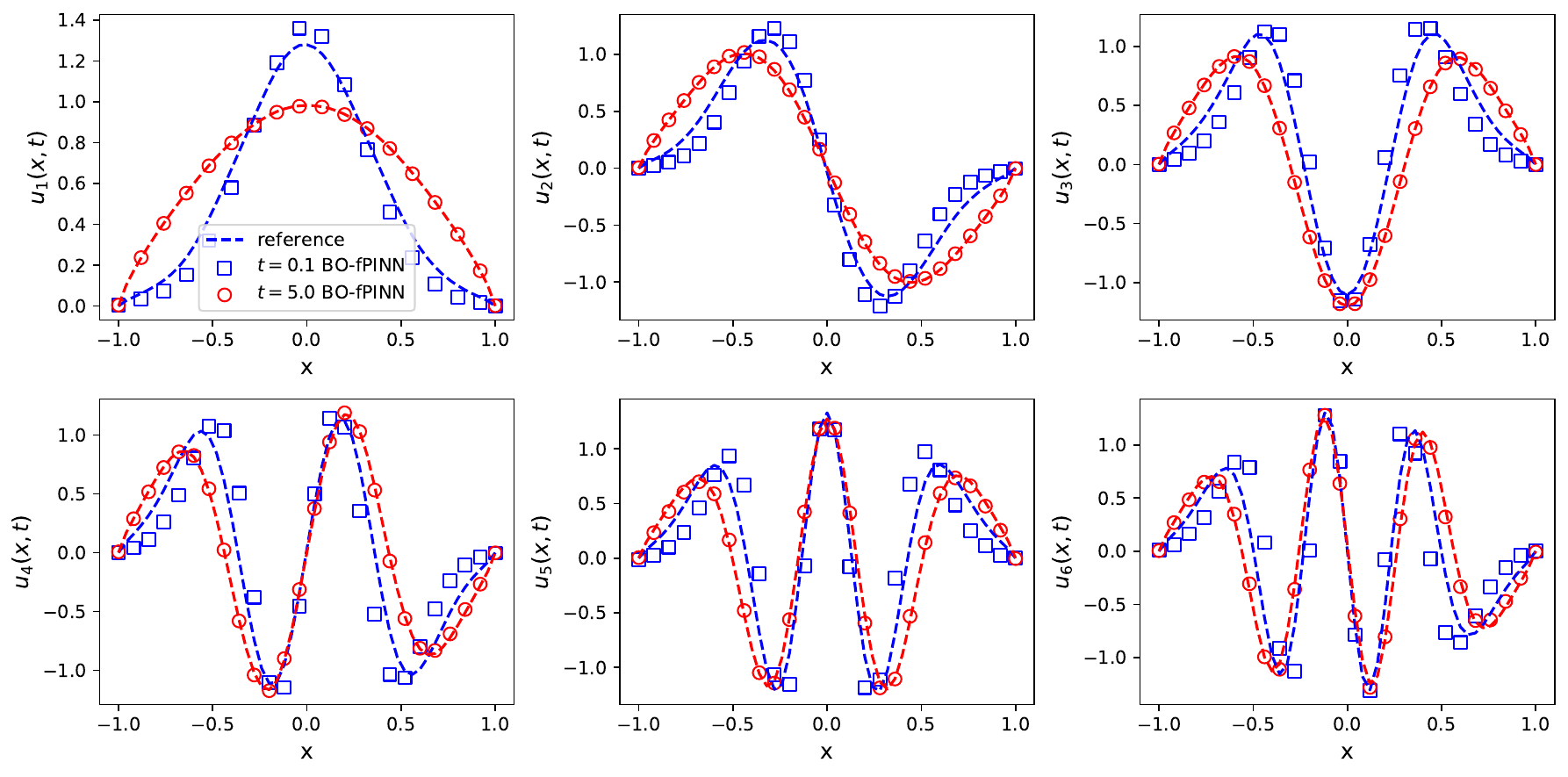}}\quad	
\caption{Forward problem with time independent random forcing. The spatial modes $u_i$ at $t=0.1$ and $t=5.0$. We use the QMC-BO solution as reference.}\label{fig:ui-1d-t5}
\end{figure}

Figures \ref{fig:mean258}a and \ref{fig:mean258}b show the BO-fPINN solution mean and variance, respectively, at $t=1.0$ with different values of the fractional order $\alpha$. Figures \ref{fig:l2-mean258}a and \ref{fig:l2-mean258}b  show the $L_2$ errors of the solution mean
and variance, respectively. We can see that the BO-fPINN solutions agree with the reference solutions very well. To further demonstrate the long-integration performance of the BO-fPINN method, we fix $\alpha=1.5$ and obtain the solution until $T=5.0$. During the training stage, we use the domain decomposition strategy proposed in \cite{zhang2020learning} to divide the time domain into five non-overlapping subdomain with equal length and then train one after the other. Figure \ref{fig:mean-1d-t5}a shows the BO-fPINN solution mean at $t=0.1$ and $t=5.0$. Figure \ref{fig:mean-1d-t5}b shows the evolution of scaling fractors $a_i$ as time grows. We compare the modal functions $u_i$ ($i=1,\cdots,6$) obtained from the BO-fPINN method with the reference solutions in Figure \ref{fig:ui-1d-t5}. These results show the efficiency of the BO representation, i.e., most of the stochasticity in this 19-dimensional SFPDE can be captured with a small number of modes. In Table \ref{tab:1d-error-t5} we report the relative $L_2$ errors for all of the BO components at the final time $T=5.0$. The proposed BO-fPINN method
generates accurate predictions at the end time (T = 5.0).
We also study the effect of the number of BO expansion modes by comparing the variances of solution calculated using six, seven and eight BO modes, as well as the performance of three different methods: BO-fPINN, gPC, and
the standard BO methods. More details about the computational results  can be found in \ref{appendix:1DRDE4.2.1}.



\begin{table}[!htp]
\centering
\renewcommand\arraystretch{1.2}
\caption{Forward problem with time independent random forcing. The relative $L_2$ errors of BO-fPINN solutions versus the reference solutions at the final time $t=5.0$.}
\label{tab:1d-error-t5}
\begin{tabular}{c|cccccc}
\hline
          &  $a_1$  &  $a_2$   &  $a_3$  &  $a_4$  & $a_5$  & $a_6$  \\

Relative $L_2$ error  & 0.037$\%$ & 0.314$\%$ & 0.554$\%$ & 0.414$\%$ & 1.273$\%$ & 1.830$\%$ \\
\hline
          &  $u_1$  &  $u_2$   &  $u_3$  &  $u_4$  & $u_5$  & $u_6$ \\

Relative $L_2$ error  & 0.0032 & 0.0132 & 0.0095 & 0.0264 & 0.0473 & 0.0595 \\
\hline
          &  $Y_1$  &  $Y_2$   &  $Y_3$  &  $Y_4$  & $Y_5$  & $Y_6$ \\

Relative $L_2$ error  & 0.0415 & 0.0522 & 0.0539 & 0.0664 & 0.0795 & 0.0890 \\
\hline 
\end{tabular}
\end{table}

\begin{figure}[!ht]
	\centering	
    \subfloat[]{\includegraphics[height=0.27\textheight, width=0.44\textwidth]{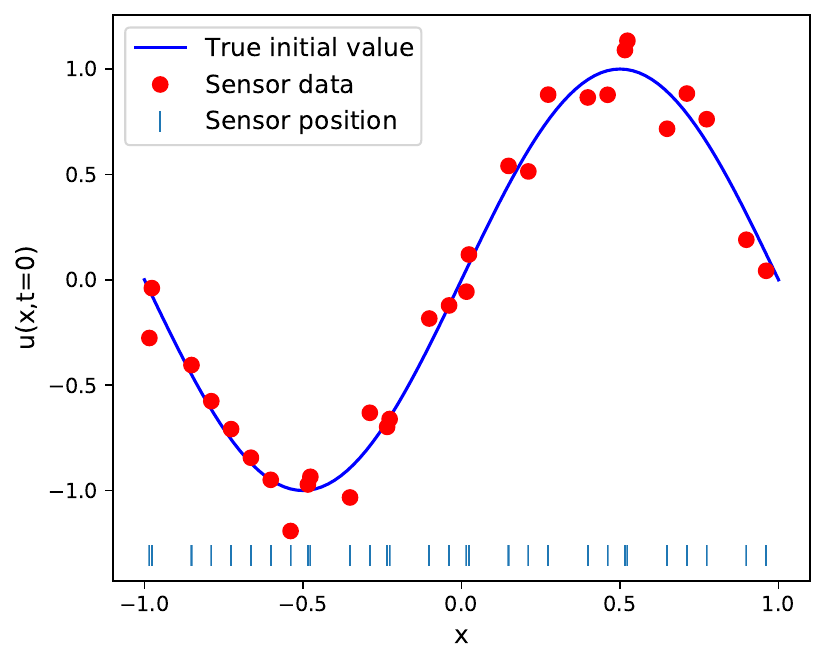}}\qquad
	\subfloat[]{\includegraphics[height=0.27\textheight, width=0.44\textwidth]{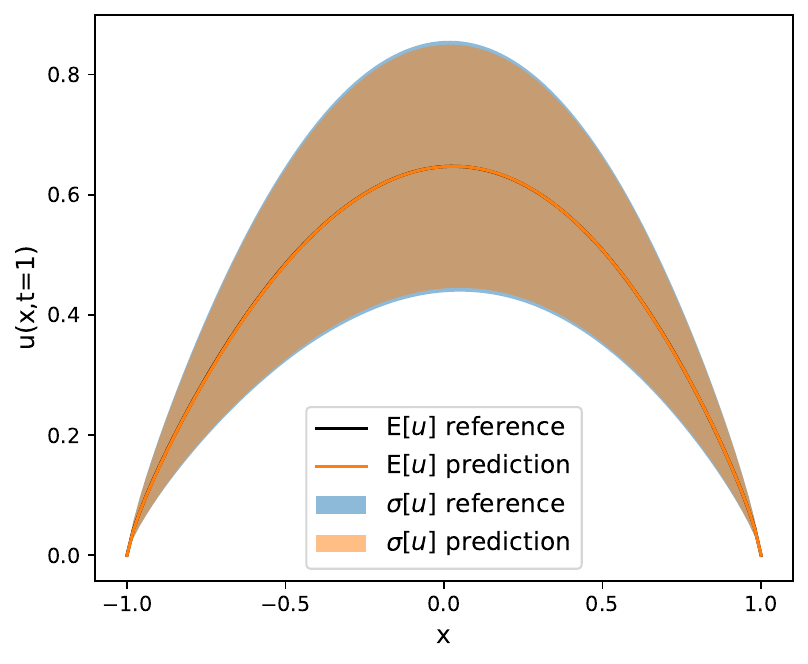}} \\	
	\caption{Forward problem with time independent random forcing (noisy data). (a):~Noisy sensor data as initial condition. (b):~Mean and standard deviation of the predicted solution $u(x,t;\omega)$ versus the reference mean and standard deviation.}\label{fig:1d-addnoise}
\end{figure}

Moreover, we solve the stochastic fractional diffusion-reaction equation with noisy sensor data as the initial condition. The noisy sensor measurements of $u(x,t=0)$ are shown in Figure \ref{fig:1d-addnoise}a, where the 30 sensors are uniformly placed in the domain and the measurements are corrupted by independent Gaussian random noise of standard deviation $0.1$. Figure \ref{fig:1d-addnoise}b presents a comparison of the BO-fPINN solution mean and standard deviation obtained from the  noisy sensor data with the reference mean and standard deviation, calculated with the Quasi Monte Carlo simulation. We can see that using noisy sensor measurements as the initial condition does not change the prediction at final time too much. The above simulation results show the superiority of the BO-fPINN method over the standard BO-method on coping with deterministic initial conditions and scattered noisy sensor measurements. Another advantage of the BO-fPINN method over the QMC-BO method is that it can solve efficiently a time-dependent nonlinear {\em inverse} stochastic problem.

\subsubsection{Inverse problem with time independent random forcing}
By using the BO-fPINN we can solve inverse problems with almost the same coding effort as solving forward problems. Here again we solve \eqref{prob:1d} with time independent random forcing given in Section 4.2.1, but we assume that we do not know the exact diffusion and reaction coefficients $\mu$, $\epsilon$, nor the  fractional order $\alpha$. Some extra information about $u(x,t;\omega)$ is provided to help infer these unknown coefficients. In this example, the extra information is the mean value of $u(x,t;\omega)$ evaluated at three locations $x=-0.5, 0, 0.5$ and at two times $t=0.1, 0.9$. We set $\sigma_h=1$ and $l_c=0.4$, and the ``hidden" values of $\mu$, $\epsilon$ and $\alpha$ are selected to be 0.5, 0.3 and 1.5, respectively. We use a BO representation with four modes, and adopt the same setup of the neural networks and training points employed in the forward problem. By including an additional term in the loss
function that calculates the MSE of the predicted $\overline{u}_{nn}(x,t)$  versus the additional measurement data, the three unknown parameters will be tuned at the training stage with the DNN parameters together. We choose the initial values of $\mu$, $\epsilon$ and $\alpha$ to be 1.0, 1.0 and 0.5$\tanh$(0.2)+1.5, respectively. The neural networks are trained with the Adam optimizer with learning rate 0.001 for 600000 epochs and then tuned with L-BFGS-B. As in the forward problem, the reference solution statistics are calculated with Quasi Monte Carlo simulation, and the reference BO components are generated by numerically solving the BO equations via the QMC-BO method.

\begin{figure}[!ht]
	\centering	
    \subfloat[]{\includegraphics[width=0.43\textwidth]{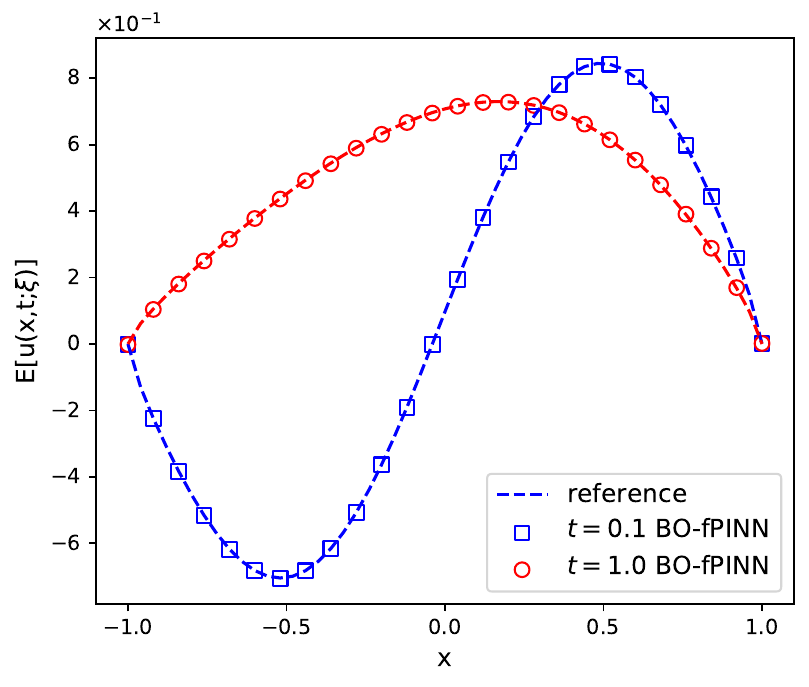}}\qquad
	\subfloat[]{\includegraphics[width=0.43\textwidth]{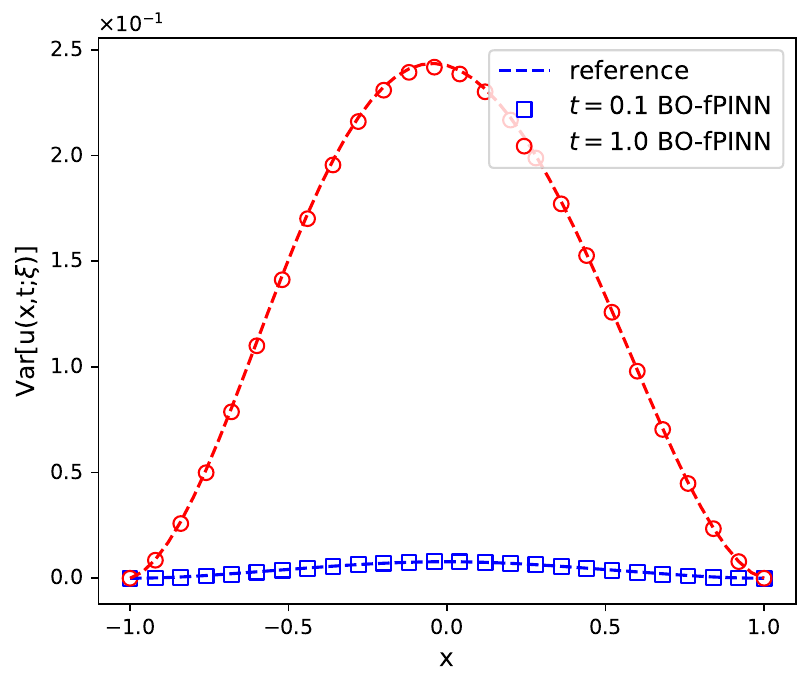}} \\	
	\caption{Inverse problem with time independent random forcing. (a):~Solution mean at $t=0.1$ and $t=1$. (b):~Solution variance at $t=0.1$ and $t=1$. The reference solutions are calculated by
solving the forward problem using Quasi Monte Carlo simulation.}\label{fig:mean-1d-inverse}
\end{figure}

\begin{figure}[!ht]
	\centering	
    \subfloat[]{\includegraphics[width=0.42\textwidth]{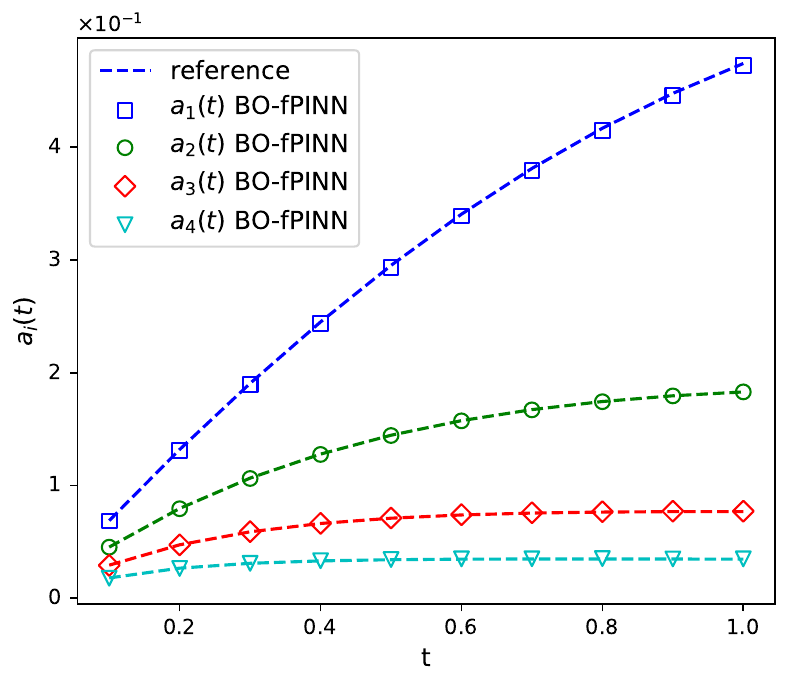}}\qquad
	\subfloat[]{\includegraphics[width=0.43\textwidth]{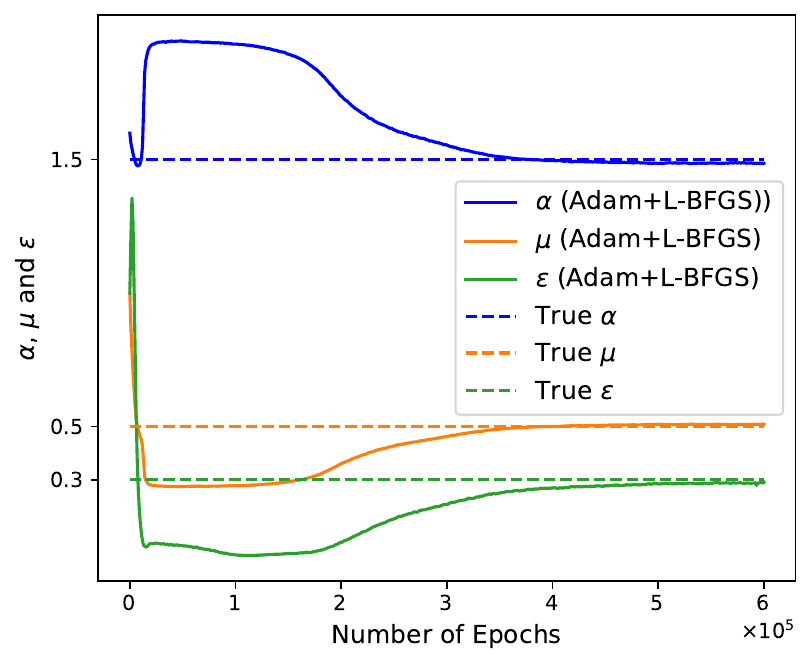}} \\	
	\caption{Inverse problem with time independent random forcing. (a):~The evolution of the scaling factors $a_i$ as the time evolving, the reference solutions are calculated by
solving the forward problem using the QMC-BO method. (b):~Parameter evolution as
the iteration of optimizer progresses. }\label{fig:ai-1d-inverse}
\end{figure}

Figure \ref{fig:mean-1d-inverse}a and \ref{fig:mean-1d-inverse}b show the predicted solution mean and variance , respectively. Figure~\ref{fig:ai-1d-inverse}a plots the scaling factors $a_i$ and Figure~\ref{fig:ai-1d-inverse}b displays the convergence history
of the predicted parameters $\mu$, $\epsilon$ and $\alpha$. We summarize the identified solutions to the inverse problems in Table \ref{tab:alpha-inverse}.  The predicted BO modes and the relative $L_2$ errors of the BO-fPINN solution versus to the reference solutions can be found in Figure \ref{fig:ui-1d-inverse} and Table \ref{tab:inverse-error}, respectively. We can see that it is evident that when compared to the reference solutions, the BO-fPINN method is still accurate for solving the inverse problem.

\begin{table}[!htp]
\centering
\renewcommand\arraystretch{1.2}
\caption{Inverse problem with time independent random forcing. Identified parameters using BO-fPINN method inverse problems with synthetic data.}
\label{tab:alpha-inverse}
\begin{tabular}{c|c}
\hline 
 True parameters                            &  Identified parameters   \\
\hline
 $\mu=0.5$, $\epsilon=0.3$, $\alpha=1.5$ & $\mu=0.508$, $\epsilon=0.291$, $\alpha=1.486$  \\
\hline 
\end{tabular}
\end{table}

\subsubsection{Forward problem with time-evolving random forcing}\label{sec:RD_time}
In this section, we still consider problem (\ref{prob:1d}) but with time-dependent stochastic  forcing
$$g(x, t;\omega)=1 +\frac{1}{10}\sum_{i=1}^{5}\gamma_{i}(t)\lambda_{i}\phi_{i}\xi_{i},$$
where $\xi_i$ are i.i.d. uniform random variables $\xi_i\sim [-1,1]$,~$\lambda_{i}$ and $\phi_{i}$ are  eigenvalue and the corresponding eigenfunction of the covariance kernel of $h(x;\omega) \sim (1-x^2)\mathcal{GP}(0, C(x_1,x_2))$, where $C(x_1,x_2)$ is a squared exponential kernel with standard deviation $1$ and correlation length $0.4$. We set the coefficients $\gamma_i(t)$ as
$\gamma(t)=[\sqrt{2(1+t)}, \sqrt{3(1+t)^2}, \sqrt{5(1+t)^3}, \sqrt{50t}, \sqrt{10t^2}]$. The diffusion coefficient is set to $\mu=0.5$ and the reaction coefficient is set to $\epsilon=0.3$. The neural networks used in the BO-fPINN method for this case can be found in the Appendix. We set $\lambda_w=1$,~$\lambda_{IC}=10$,~$\lambda_{BO}=10$,~$\lambda_{0}=0.1$. We can enforce the boundary conditions by multiplying $\bar{u}_{nn}$ and $U_{nn,i}$ with $1-x^2$, so we don't consider the $MSE_{BC}$ in the loss function for this example. We use $n_x=71$ equidistantly distributed training points $\{x_c^k\}_{k=1}^{nx}$ in space, and $n_t=70$ uniformly distributed training points $\{t_c^s\}_{k=1}^{n_t}$ in the time domain, and $n_l=1000$ random samples $\{\xi_c^l\}_{l=1}^{n_{\xi}}$ in the 5-dimensional random space. The neural networks are trained with an Adam optimizer with learning rate 0.0001 till 400000 epochs and then L-BFGS-B optimizer is implemented.

We use five BO expansion terms in (\ref{eqn:BO-nn}) for this example. To obtain the reference solution for the BO decomposition, the QMC-BO method with switch time $t_s=0.01$ is employed to solve the original BO equations. The parameters for QMC-BO method are $\Delta t=5\times 10^{-5}$, $\Delta x=0.02$, $M=1000$. A Quasi Monte Carlo method with 1000 samples is employed to solve the SFPDE to obtain the reference for the solution statistics.

\begin{figure}[!ht]
	\centering	
    \subfloat[]{\includegraphics[width=0.43\textwidth]{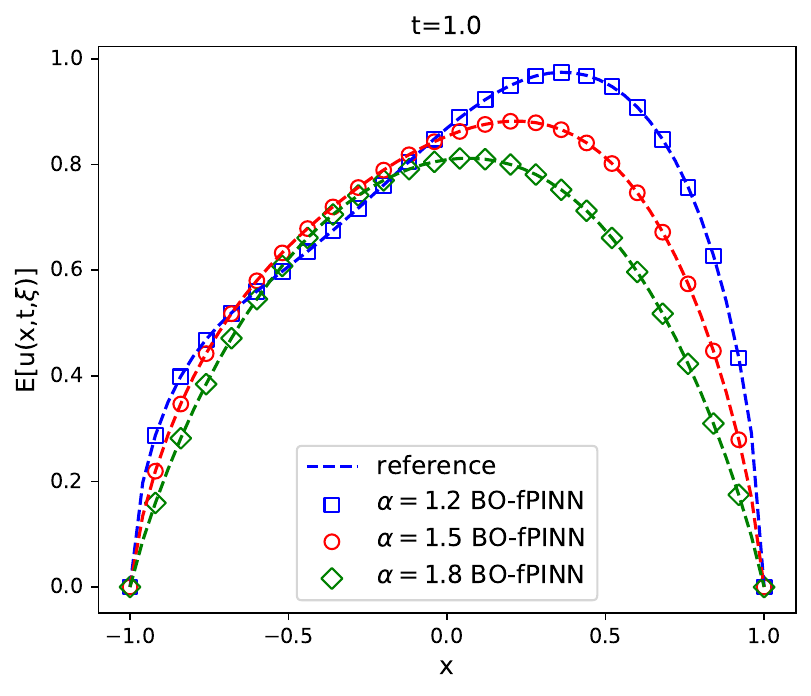}}\qquad
	\subfloat[]{\includegraphics[width=0.43\textwidth]{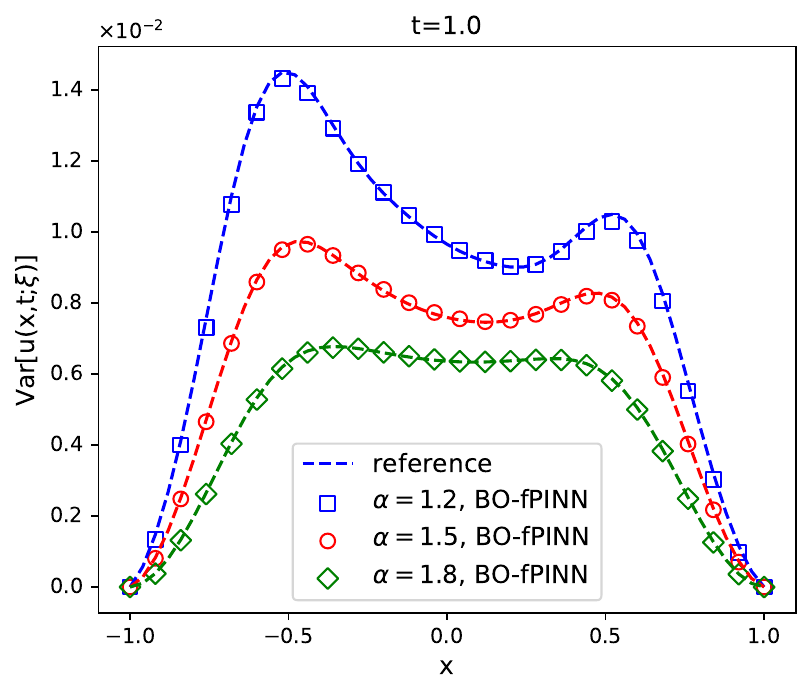}} \\
	\caption{Forward problem with time-evolving random forcing. (a):~BO-fPINN solution mean at $t=1.0$ with different values of the fractional order $\alpha$. (b):~BO-fPINN solution variance at $t=1.0$ with different values of the fractional order $\alpha$. The reference mean and variance are calculated from the Quasi Monte Carlo simulation.}\label{fig:1D-alpha121518-mean-and-variance}
\end{figure}

\begin{figure}[!ht]
	\centering	
    \subfloat[]{\includegraphics[width=0.43\textwidth]{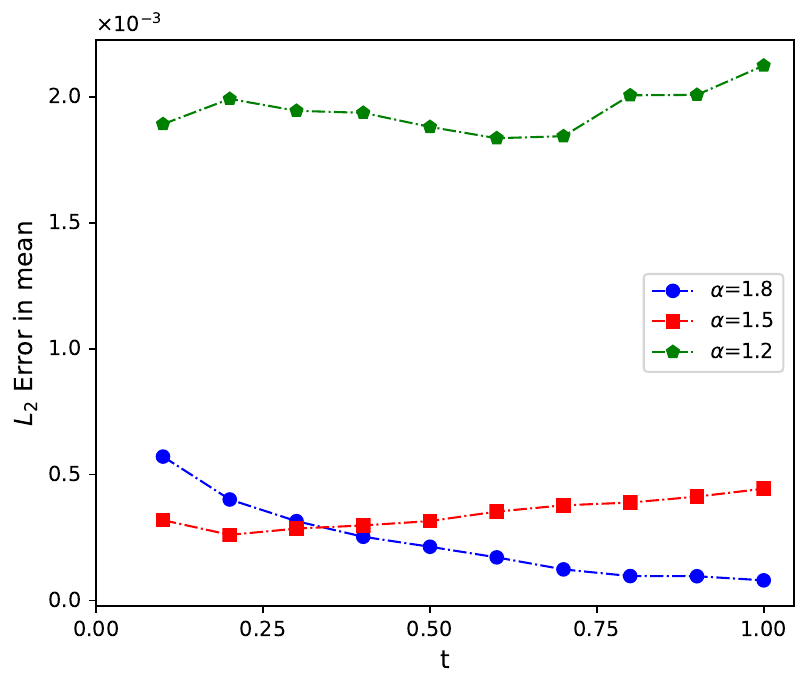}}\qquad
	\subfloat[]{\includegraphics[width=0.43\textwidth]{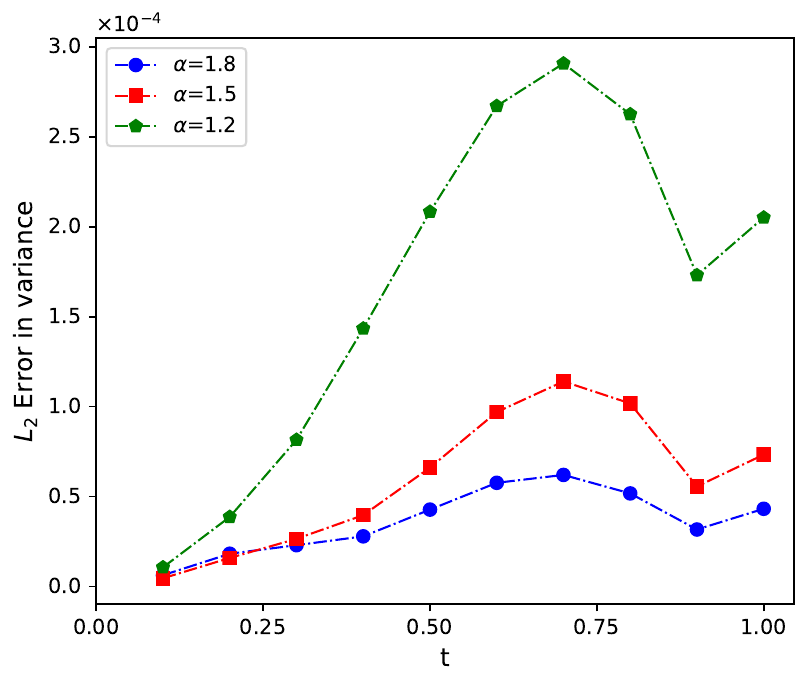}} \\
	\caption{Forward problem with time-evolving random forcing. (a):~The $L_2$ error in the mean obtained by BO-fPINN method with different values of the fractional order $\alpha$. (b):~The $L_2$ error in the variance obtained by BO-fPINN method with different values of the fractional order $\alpha$. We use the QMC solution as reference.
}\label{fig:1D-alpha121518-l2error}
\end{figure}

\begin{figure}[!ht]
	\centering	
    \subfloat[]{\includegraphics[width=0.43\textwidth]{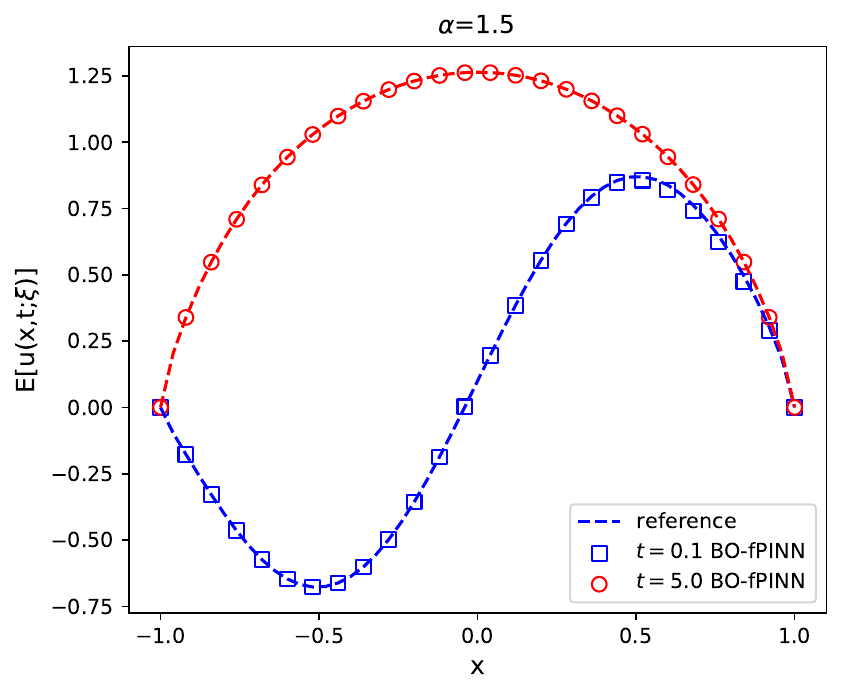}}\qquad
	\subfloat[]{\includegraphics[width=0.415\textwidth]{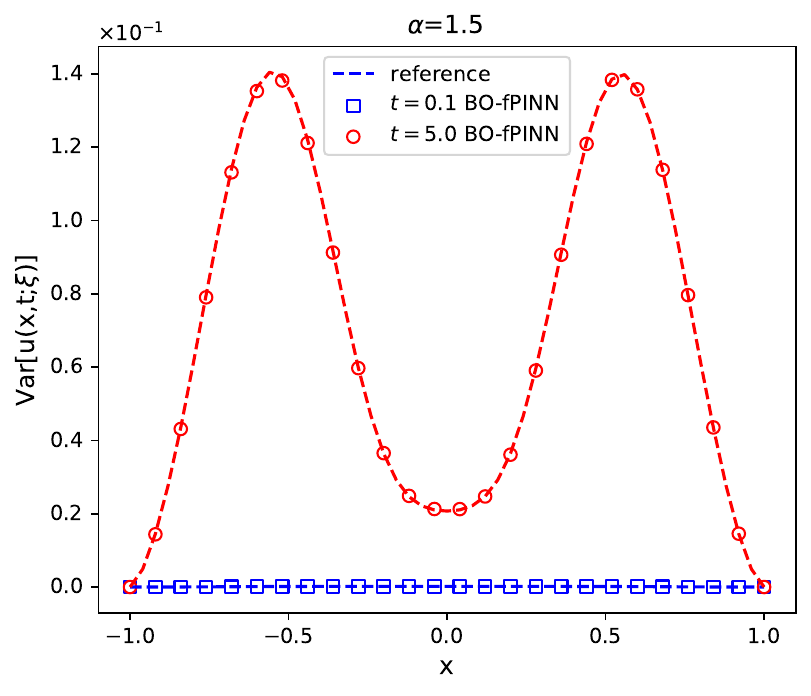}} \\	
	\caption{Forward problem with time-evolving random forcing. (a):~ Solution mean at $t=0.1$ and $t=5.0$. (b):~ Solution variance at $t=0.1$ and $t=5.0$. The reference mean and variance are calculated from the Quasi Monte Carlo simulation.}\label{fig:1D-alpha15-mean-and-variance-longtime}
\end{figure}

\begin{figure}[!ht]
\centering{\includegraphics[height=0.34\textheight, width=0.84\textwidth]{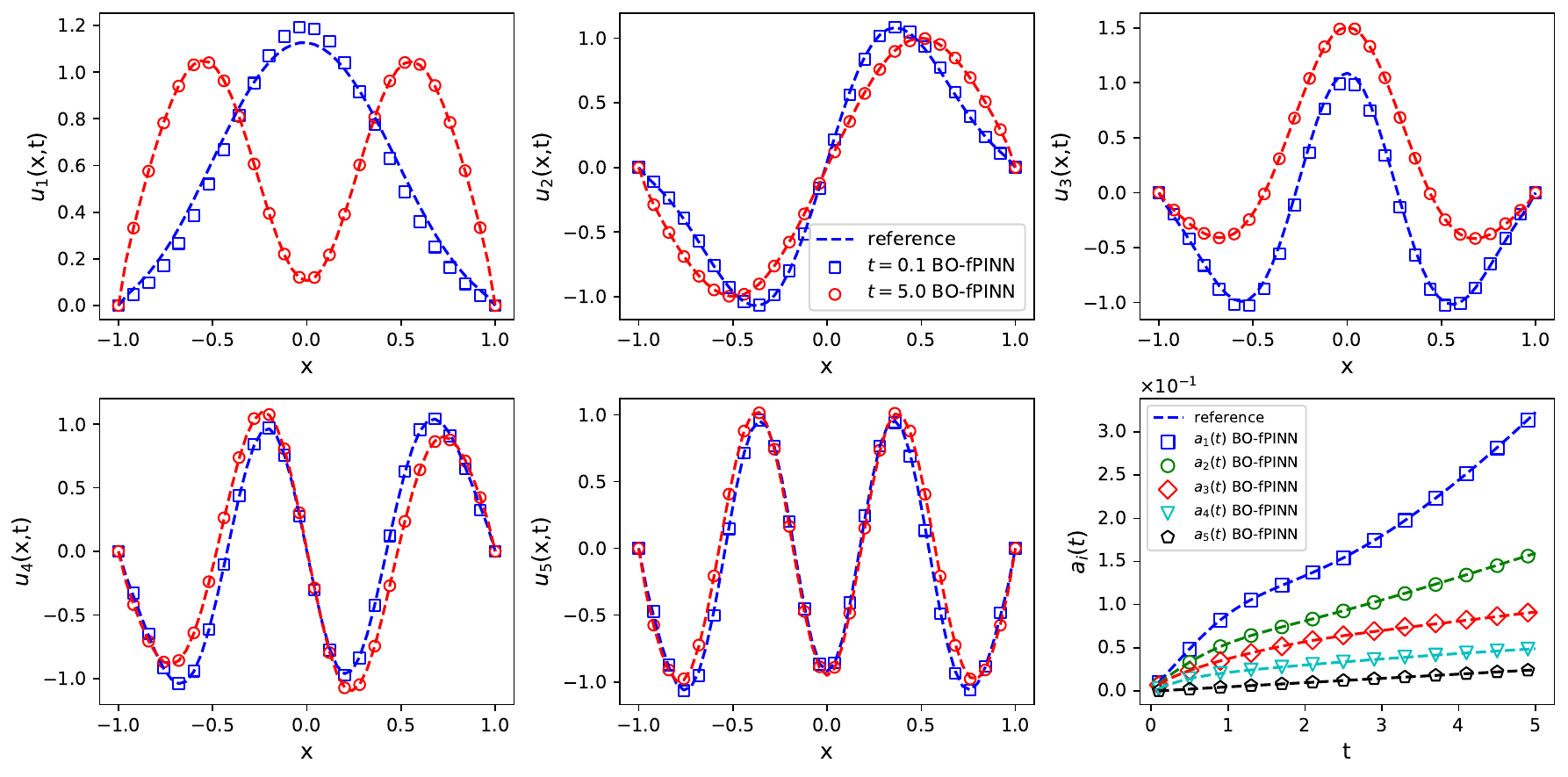}}\quad	
\caption{Forward problem with time-evolving random forcing. The spatial modes $u_i$ at $t=0.1$ and $t=5.0$. The last subplot shows the scaling factors $a_i$ at different time steps. We use the QMC-BO solution as reference.}\label{fig:1D-alpha15-ui-and-ai-longtime}
\end{figure}

Figures \ref{fig:1D-alpha121518-mean-and-variance}a and \ref{fig:1D-alpha121518-mean-and-variance}b show the BO-fPINN solution mean and variance, respectively, at $t=1.0$ with different values of the fractional order $\alpha$. Figure \ref{fig:1D-alpha121518-l2error}a and \ref{fig:1D-alpha121518-l2error}b show the L2 errors of the solution mean and variance, respectively. We can see that the BO-fPINN solutions agree closely with the reference solution. To further demonstrate the long-integration performance of the BO-fPINN method, we fix $\alpha=1.5$ and obtain the solution until $T=5.0$.  Figures \ref{fig:1D-alpha15-mean-and-variance-longtime}a and \ref{fig:1D-alpha15-mean-and-variance-longtime}b show the BO-fPINN solution mean and variance at $t=0.1$ and $t=5.0$, respectively. In Figure \ref{fig:1D-alpha15-ui-and-ai-longtime}, we compare the modal functions $u_i$ obtained from the BO-fPINN method with the reference solutions, while the last subplot of Figure \ref{fig:1D-alpha15-ui-and-ai-longtime} shows the evolution of scaling factors $a_i$, where the modes gradually pick up energy as the result of the nonlinear source term. Although the BO representation with 5 modes fits well the reference solution, we should note that more modes should be added adaptively to achieve better accuracy for the long-time integration. This will be addressed more systematically in  future work. In Table \ref{tab:1D-ai-ui-yi-rel2error-t=5} we report the relative $L_2$ errors for all of the BO components at the final time $T=5.0$. The proposed BO-fPINN method generates accurate predictions at the end time (T = 5.0).

%
%

\begin{table}[!htp]
\centering
\renewcommand\arraystretch{1.2}
\caption{Forward problem with time-evolving random forcing. The relative $L_2$ errors of BO-fPINN solutions versus the reference solutions at the final time $t=5.0$.}
\label{tab:1D-ai-ui-yi-rel2error-t=5}
\begin{tabular}{c|cccccc}
\hline
          &  $a_1$  &  $a_2$   &  $a_3$  &  $a_4$  & $a_5$  \\

Relative $L_2$ error  & 0.220$\%$ & 0.363$\%$ & 0.360$\%$ & 0.038$\%$ & 0.020$\%$ \\
\hline
          &  $u_1$  &  $u_2$   &  $u_3$  &  $u_4$  & $u_5$  \\

Relative $L_2$ error  & 0.270$\%$ & 0.577$\%$ & 0.497$\%$ & 0.670$\%$ & 0.921$\%$  \\
\hline
          &  $Y_1$  &  $Y_2$   &  $Y_3$  &  $Y_4$  & $Y_5$  \\

Relative $L_2$ error  & 0.0304 & 0.0406 & 0.0713 & 0.0479 & 0.0943 \\
\hline 
\end{tabular}
\end{table}


\subsubsection{Transfer Learning}
In general, we need to train a neural network from the beginning to obtain the neural network surrogate model when we vary some parameters in a SFPDE, such as the order of the fractional derivative. This is similar to the classical numerical method, but the neural network training has a relatively slow convergence speed due to the large number of parameters, which makes the computational  cost of retraining the network expensive. Here, we employ transfer learning to speed up the training when we solve SFPDEs with different order of the fractional derivatives~\cite{jin2021nsfnets,goswami2019transfer}.

\begin{figure}[!ht]
	\centering	
    \includegraphics[width=0.52\textwidth,height=0.20\textheight]{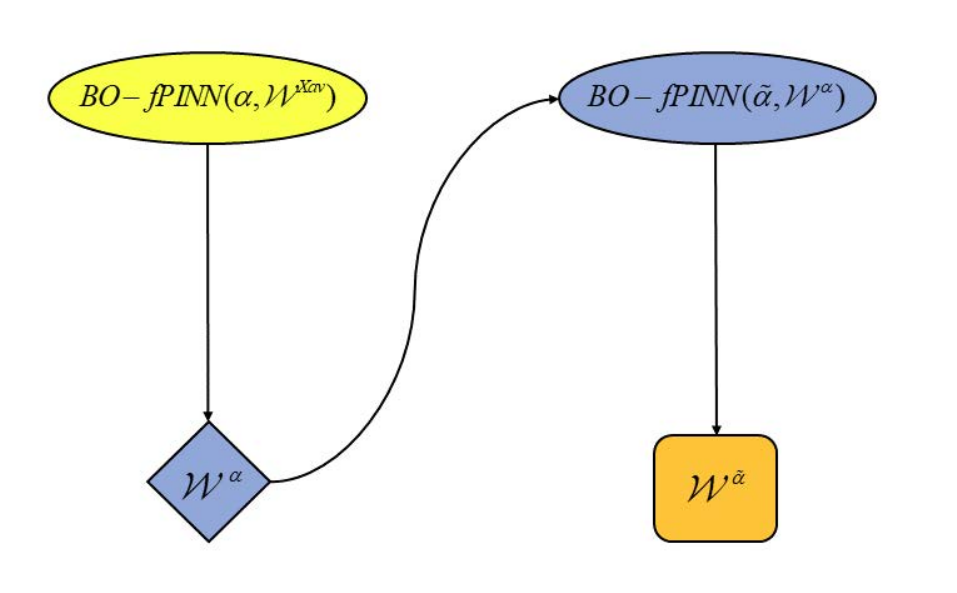}		
	\caption{A source BO-fPINN (top left), tasked to solve a SFPDE with fractional order $\alpha$ and Xavier initialization $W^{Xav}$, is trained and the weight vector $W^{\alpha}$ (bottom left) is stored. A new BO-fPINN (top right), tasked to solve a SFPDE with fractional order $\tilde{\alpha}$ and reference initialization $W^{\alpha}$,  is trained and the weight become $W^{\tilde{\alpha}}$.}\label{fig:flow-tansfer}
\end{figure}

Concretely, let $\mathcal{P}(\alpha)$ represent the variable fractional order space, then the transfer learning process is mainly divided into two steps: (1) we first train the neural network under the configuration that the variable fractional order space is $\mathcal{P}^{\alpha}$, and the network parameter model is initialized by Xavier ($W^{Xav}$) to obtain the model parameter of the source problem $W^{\alpha}$; (2) then, we can obtain the model parameter $W^{\tilde{\alpha}}$ of the target problem by training the neural network under the configuration that the variable fractional order space is $\mathcal{P}^{\tilde{\alpha}}$ and the network parameter model is initialized with $W^{\alpha}$. The schematic of transfer learning between different fractional orders is given in  Figure \ref{fig:flow-tansfer}.

We consider two different scenarios here for the forward problem: the first scenario is to solve the target problem with $\tilde{\alpha}=1.2, 1.5$ by transfer learning based on source problem with $\alpha=1.8$; the second is to solve the target problem with $\tilde{\alpha}=1.5, 1.8$ by transfer learning based on source problem with $\alpha=1.2$. The neural network structure and the number of training points used here are the same as those used in the forward problem. To demonstrate the superiority of transfer learning, we compare the computation time of the two optimization strategies: 1) Initialize the target model by $W^{\alpha}$, and fine-tune the target model only by using the L-BFGS-B via transfer learning; 2) Initialize the target model by Xavier, and optimize the target model by using the Adam and L-BFGS-B. In this section, we only show the results of the source model with $\alpha=1.8$ in transfer learning, and the results of the source model with $\alpha=1.2$ are presented in the \ref{appendix:1DRDE4.2.4}.

\begin{figure}[!ht]
\centering
\subfloat[]{\includegraphics[width=0.40\linewidth]{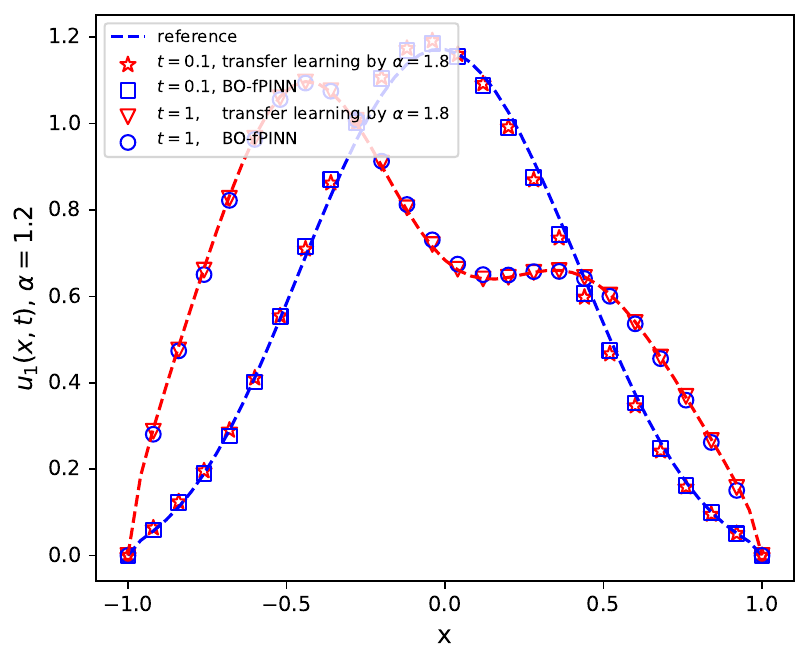}}\quad
\subfloat[]{\includegraphics[width=0.40\linewidth]{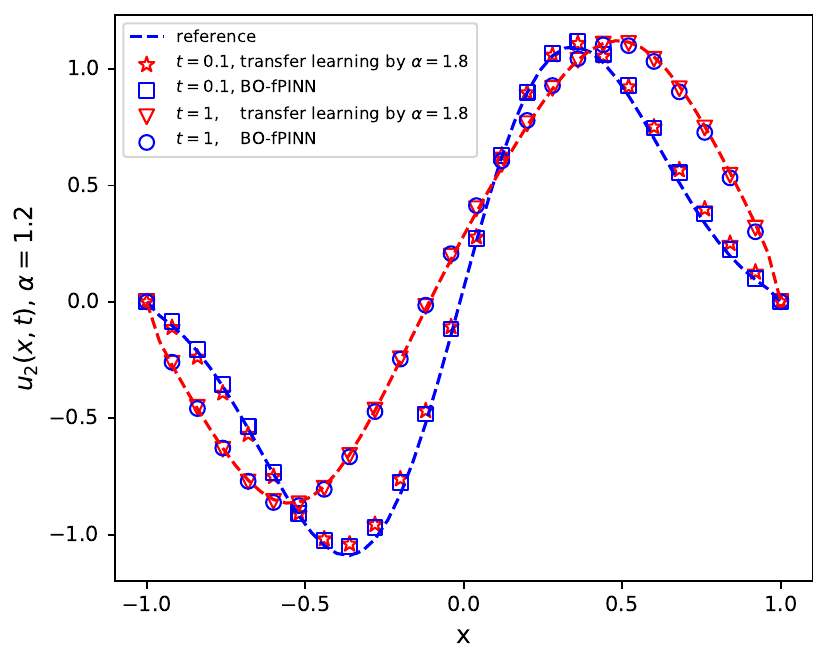}}\\
\subfloat[]{\includegraphics[width=0.40\linewidth]{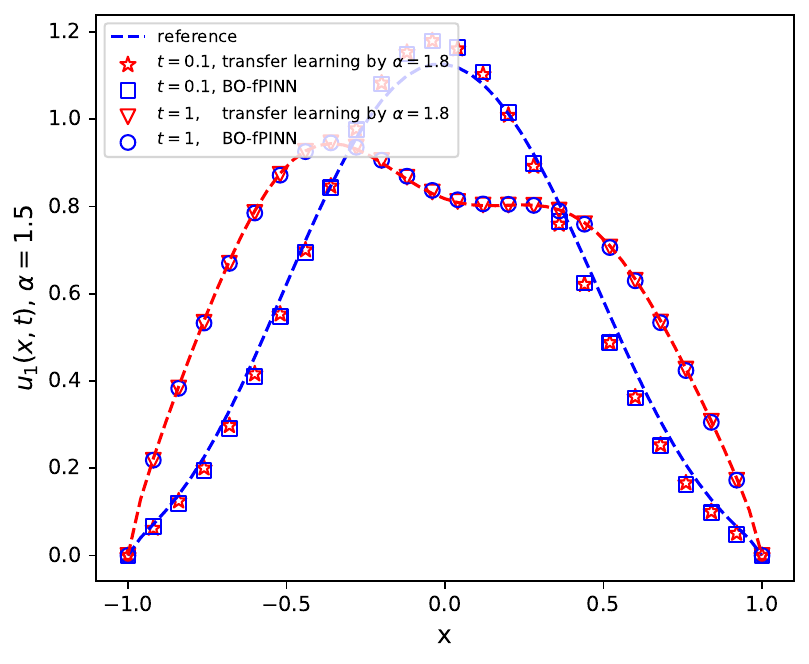}}\quad
\subfloat[]{\includegraphics[width=0.40\linewidth]{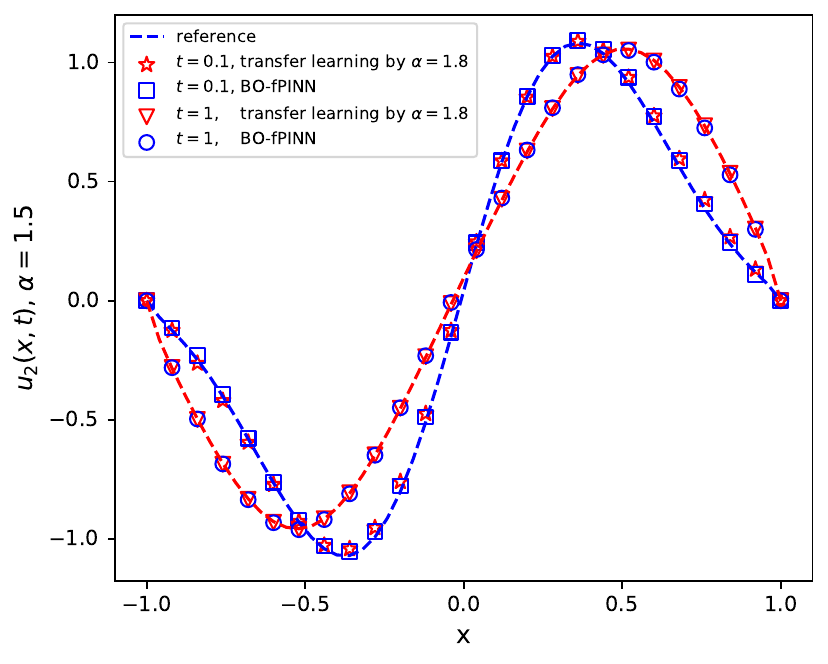}}\\
\caption{Forward problem with time-evolving random forcing (transfer learning between different fractional derivative orders). (a)-(b): The first  and second BO modes at $t=0.1, 1$ for $\alpha=1.2$ obtained by transfer learning with source model  $\alpha=1.8$, respectively. (c)-(d): The first and second BO modes at $t=0.1, 1$ for $\alpha=1.5$ obtained by transfer learning with source model  $\alpha=1.8$, respectively.The reference solutions are calculated by the QMC-BO method.}\label{fig:1D-transfer-u1-u2}
\end{figure}

Figures \ref{fig:1D-transfer-mean-and-variance}a and  \ref{fig:1D-transfer-mean-and-variance}b show the predicted solution mean and variance with $\alpha=1.2$ at $t=0.1, 1.0$ obtained by using the Xavier to initialize the target model and by using transfer learning with source model $\alpha=1.8$, respectively. We compare the modal functions $u_i$ ($i=1, 2$) obtained by directly training and transfer learning with source model $\alpha=1.8$ in Figure \ref{fig:1D-transfer-u1-u2}a-b. The trained source models used in the transfer learning here are all initialized with Xavier, trained 400000 epochs with the Adam optimizer, and then fine-tuned with the L-BFGS-B optimizer. The results for target model $\alpha=1.5$ with transfer learning by the source model $\alpha=1.8$ are reported in \ref{fig:1D-transfer-mean-and-variance}c-d and \ref{fig:1D-transfer-u1-u2}c-d. We can see that the results obtained by transfer learning with different parameters are in good agreement with the reference solution but with much less training time, which shows the effectiveness of transfer learning.

\begin{figure}[htp!]
\centering
\subfloat[]{\includegraphics[width=0.40\linewidth]{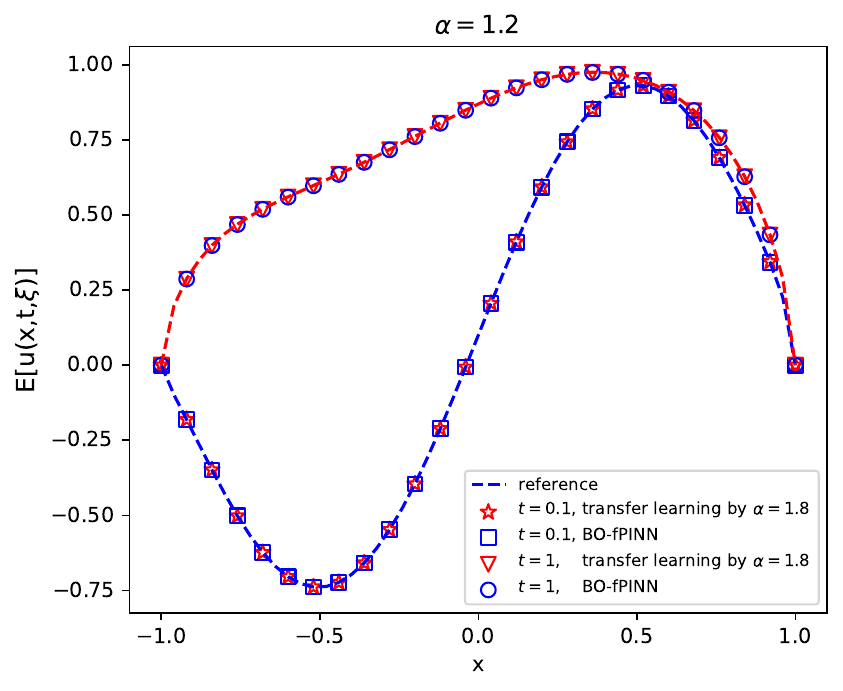}}\quad
\subfloat[]{\includegraphics[width=0.38\linewidth]{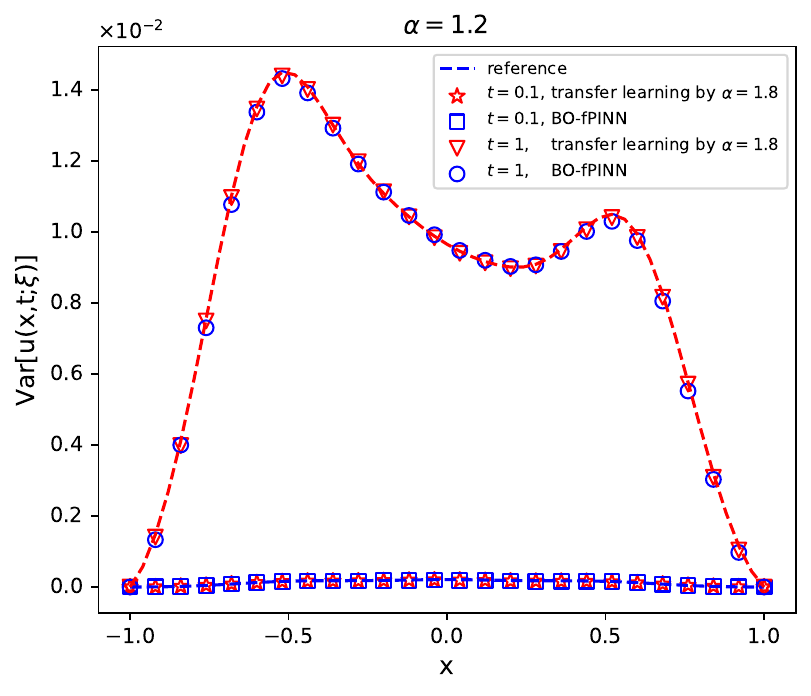}}
\\
\subfloat[]{\includegraphics[width=0.39\linewidth]{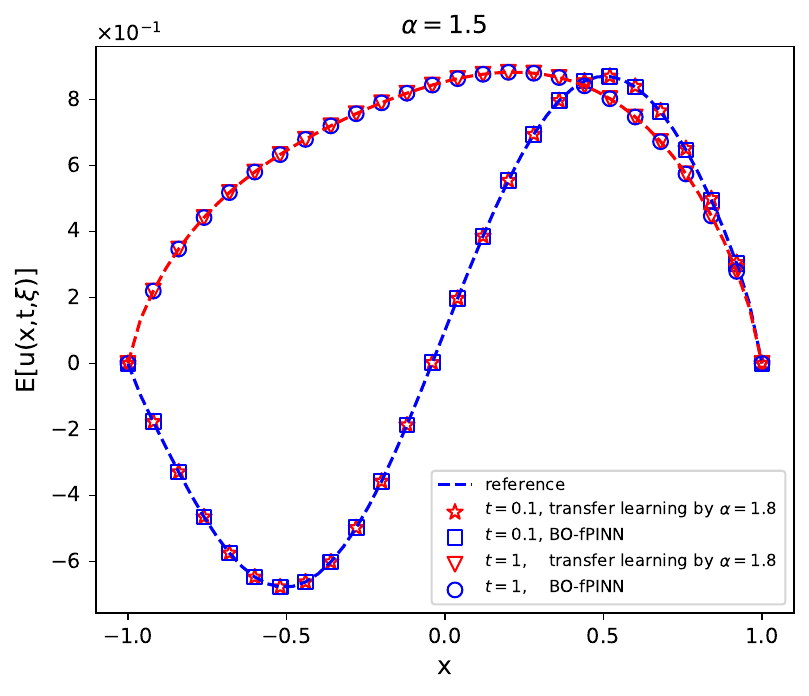}}\quad
\subfloat[]{\includegraphics[width=0.39\linewidth]{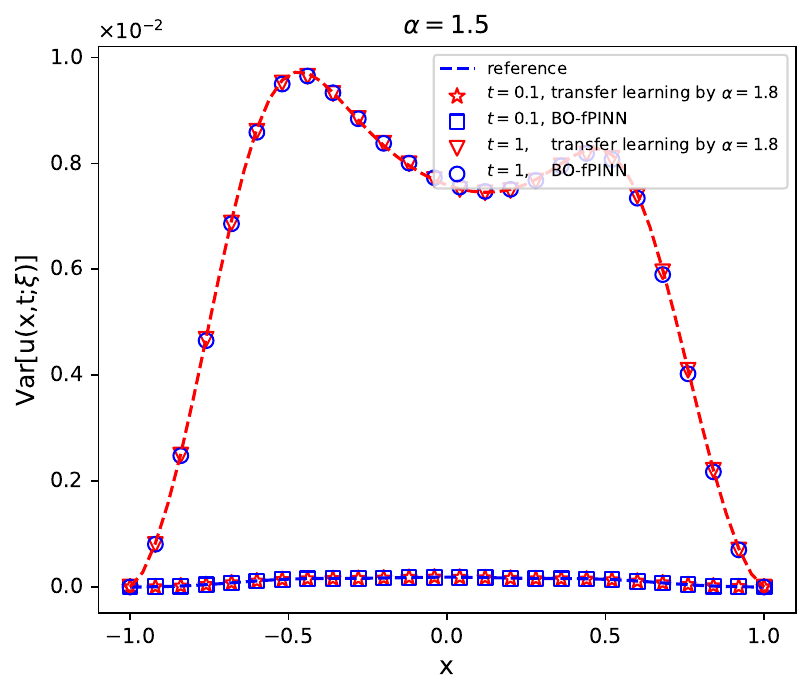}}
\caption{Forward problem with time-evolving random forcing (transfer learning between different fractional derivative orders). (a)-(b): Solution mean and variance of $\alpha=1.2$ at $t=0.1, 1$ obtained by using BO-fPINN with Xavier initialization and by transfer learning with the parameters from source model $\alpha=1.8$ as initialization, respectively. (c)-(d): Solution mean and variance of $\alpha=1.5$ at $t=0.1, 1$ by using BO-fPINN with Xavier initialization and by transfer learning with the parameters from source model $\alpha=1.8$ as initialization, respectively. The reference solutions are calculated by the QMC-BO method.}\label{fig:1D-transfer-mean-and-variance}
\end{figure}

\begin{figure}[!ht]
\centering
\subfloat[]{\includegraphics[width=0.41\linewidth]{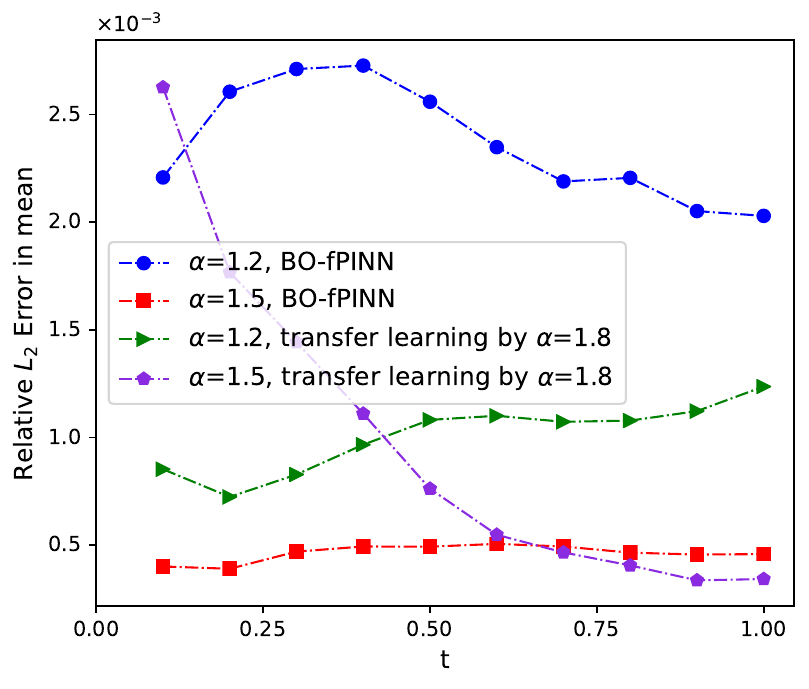}}
\quad
\subfloat[]{\includegraphics[width=0.40\linewidth]{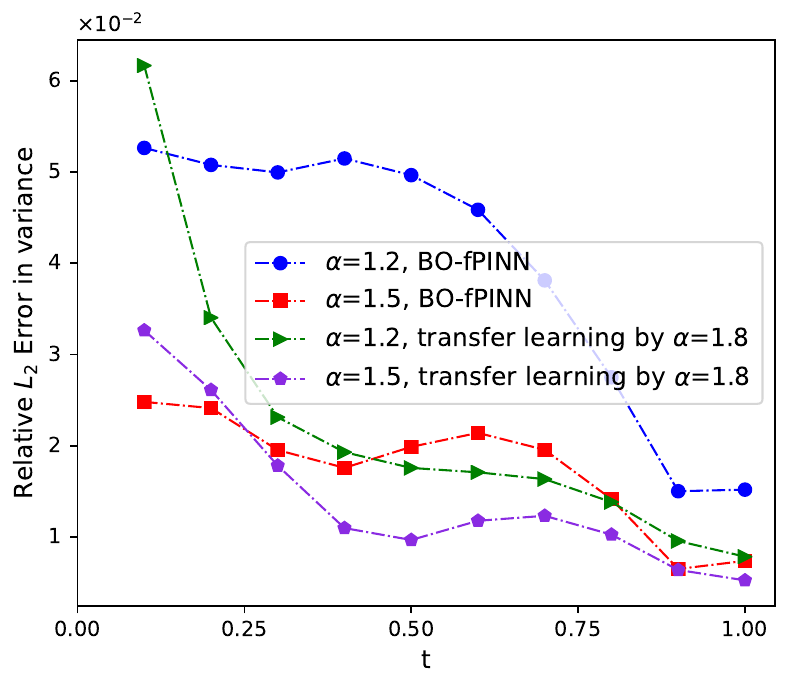}}
\caption{Forward problem with time-evolving random forcing (transfer learning between different fractional derivative orders). (a)~: The relative $L_2$ error in mean calculated by BO-fPINN and transfer learning with source model $\alpha=1.8$. (b)~: The relative $L_2$ error in variance calculated by BO-fPINN and transfer learning with source model $\alpha=1.8$. The reference solutions are calculated by the QMC-BO method. }\label{fig:1D-err-transfer}
\end{figure}

We compare the relative $L_2$ errors of the mean and variance obtained by using the source model transfer learning with parameters $\alpha=1.8$ and by initializing with Xavier the target model directly versus the reference solution in Figure \ref{fig:1D-err-transfer}. Figure \ref{fig:1D-err-transfer}a and Figure \ref{fig:1D-err-transfer}b plots the relative $L_2$ error of the predicted solution mean and variance obtained by transfer learning with source model $\alpha=1.8$. It is evident that when compared with the reference solution, the BO-fPINN with transfer learning still obtains accurate results. However, transfer learning runs almost 60 times faster than using the Adam and L-BFGS-B optimizers to train the model from scratch, as reported in  Table \ref{tab:transfer}.

\begin{table}[!htp]
\centering
\renewcommand\arraystretch{1.2}
\caption{Forward problem with time-evolving random forcing (transfer learning between different fractional derivative orders). Number of iterations and time of two optimization strategies. The training is performed on a RTX 2090 GPU.}
\label{tab:transfer}
\begin{tabular}{c|cc}
\hline 
 Optimization strategy                                                          & Number of iterations     & relative time  \\
\hline
Adam(40w)+L-BFGS-B($W^{Xav}$) & 308677   & 58.2  \\
\hline 
Finetuning with L-BFGS-B($W^{\alpha=1.8}$)    & 4762     & 1    \\
\hline
\end{tabular}
\end{table}

\subsection{2D stochastic fractional reaction-diffusion equation}\label{sec:RD_2d}

We consider the following 2D problem
\begin{equation}
  \begin{cases}
    u_t-D_{|x_1|}^\alpha u - D_{|x_2|}^\alpha -K(1+\sigma\xi)f(u)=g(x;\omega),  \quad & x\in \mathcal{D}=(0,1)^2,\quad t\in (0,T], \\
    u(x,0;\omega) = \sin(2\pi x_1)\sin(2\pi x_2),\quad &x\in \mathcal{D}, \quad\omega\in \Omega,\\
    u|_{\partial D}=0 \quad& x \in \partial\mathcal{D},\quad \omega\in \Omega, \quad t\in [0,T],
  \end{cases}
\end{equation}
where $\xi \sim U(-\sqrt{3},\sqrt{3})$, $f(u)=u(1-u^2)$,  $K=1$, $\sigma=1$, the randomness comes from both $\xi$ and the forcing term, $g(x;\omega)$,  is given by
\begin{align}
  g(x;\omega) = \mu_{g}(x) + 3\sum_{i=0}^{P_1}\sum_{j=i}^{P_2}\sqrt{v_iv_j}\phi_i(x_1)\phi_j(x_2)\xi_{ij}(\omega), \quad \xi_{ij}\sim U(-\sqrt{3}, \sqrt{3}) 
\end{align}
where  $\mu_{g}(x)=1$, $v_0:=\frac{1}{2}$, $\phi_0(x):=1$, and
\begin{align*}
  v_i:=\frac{1}{2}\exp\left(-\pi i^2l^2\right),\quad \quad \phi_i(x_k):=\sqrt{2}\cos(i\pi x_k), \quad k=1,2,i\geq 1,
\end{align*}
here we set $l=\frac{1}{3}$ and $P_1=P_2=4$.
We consider the case of $\alpha=1.8$. The neural networks used in the BO-fPINN method for this case can be found in \ref{appendix:network}. We note that the number of independent neural networks $A_{nn}$ is determined by the number of BO expansion terms. The weighting coefficients of the loss function are taken as $\lambda_w=1$,~$\lambda_{IC}=5$,~$\lambda_{BC}=50$,~$\lambda_{BO}=5$,~$\lambda_{0}=0$. We use $n_{x_1}=n_{x_2}=31$ equidistantly distributed training points $\{x_1^k\}_{k=1}^{nx_1}$ and $\{x_2^k\}_{k=1}^{nx_2}$ in spatial space, and $n_l=1000$ random samples $\{\xi^l\}_{l=1}^{n_{\xi}}$ in the random space. We sample 20 and 40 points from the uniform distribution of $[0,0.1]$ and $[0.1,1]$, respectively, to form the training point $\{t^s\}_{k=1}^{n_t}$ in the time domain. The neural networks are trained with an Adam optimizer with learning rate 0.0001 till 400000 epochs and then L-BFGS-B optimizer is implemented.

\begin{figure}[!htp]
\centering
\subfloat[]{\includegraphics[width=0.32\linewidth]{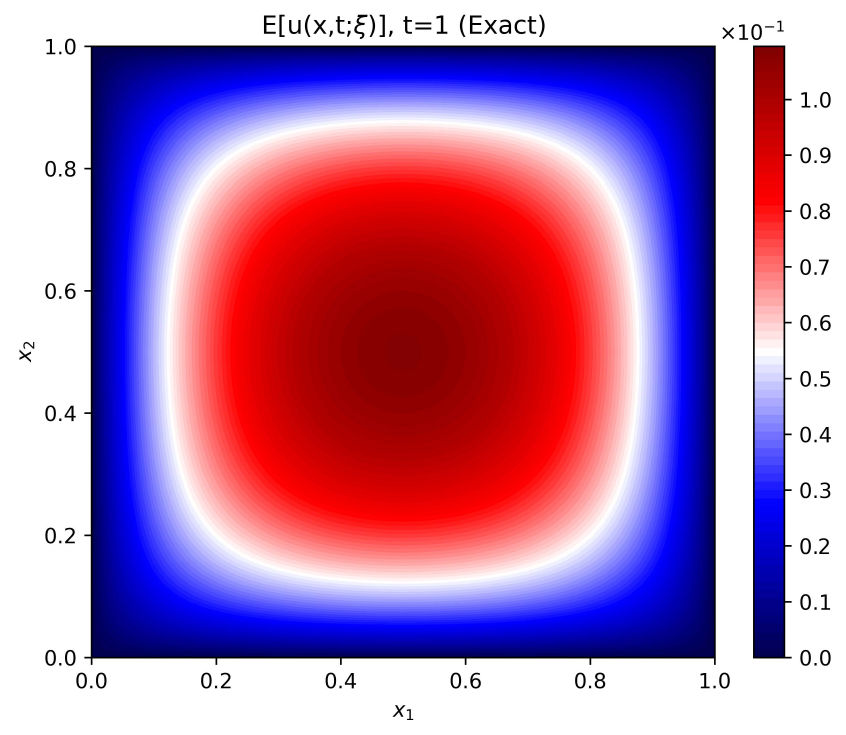}}\hfill
\subfloat[]{\includegraphics[width=0.32\linewidth]{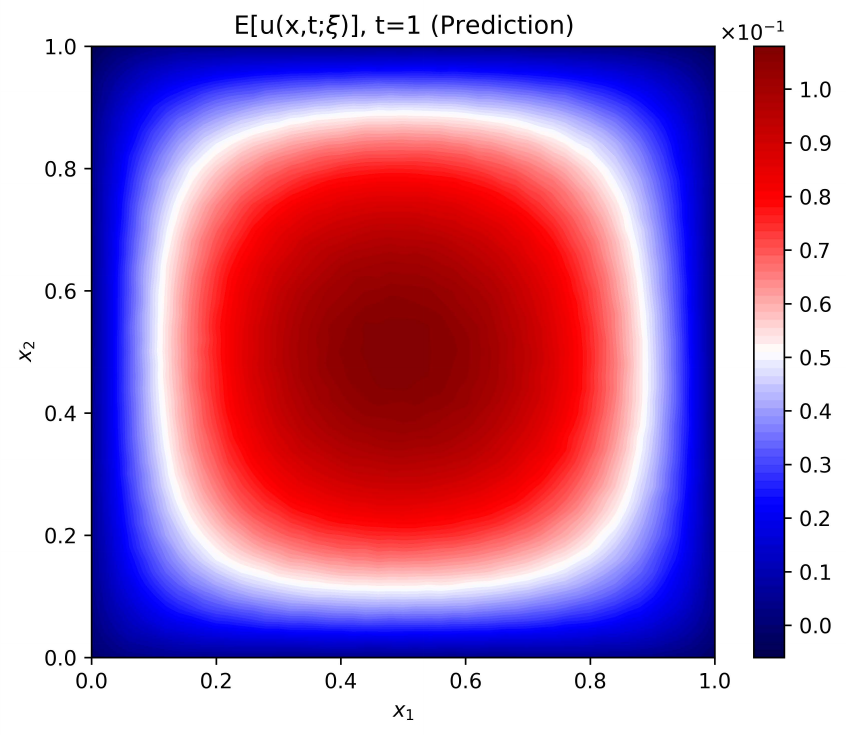}}\hfill
\subfloat[]{\includegraphics[width=0.32\linewidth]{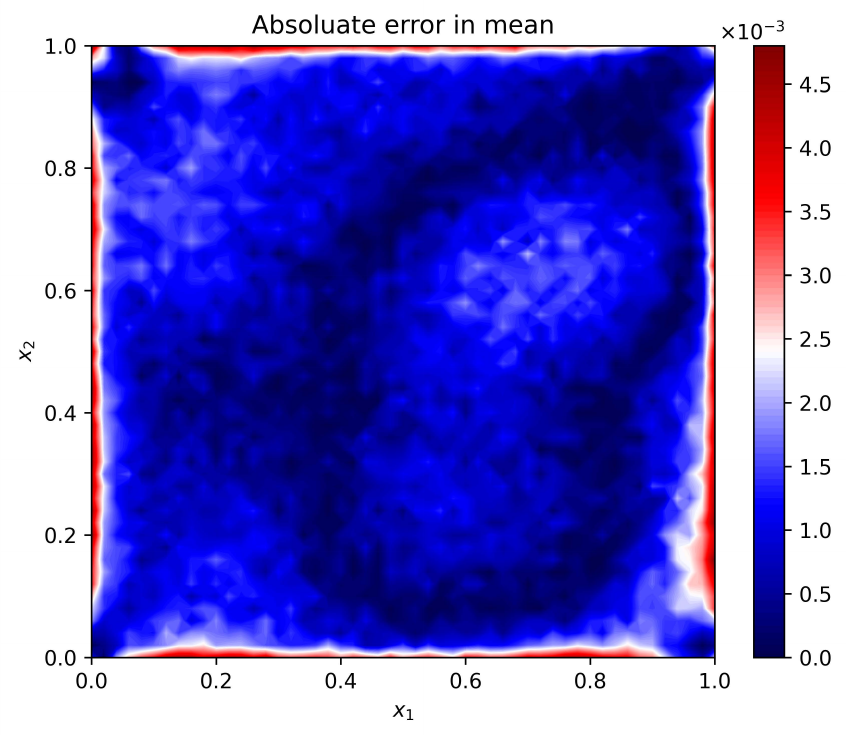}}
\\
\subfloat[]{\includegraphics[width=0.32\linewidth]{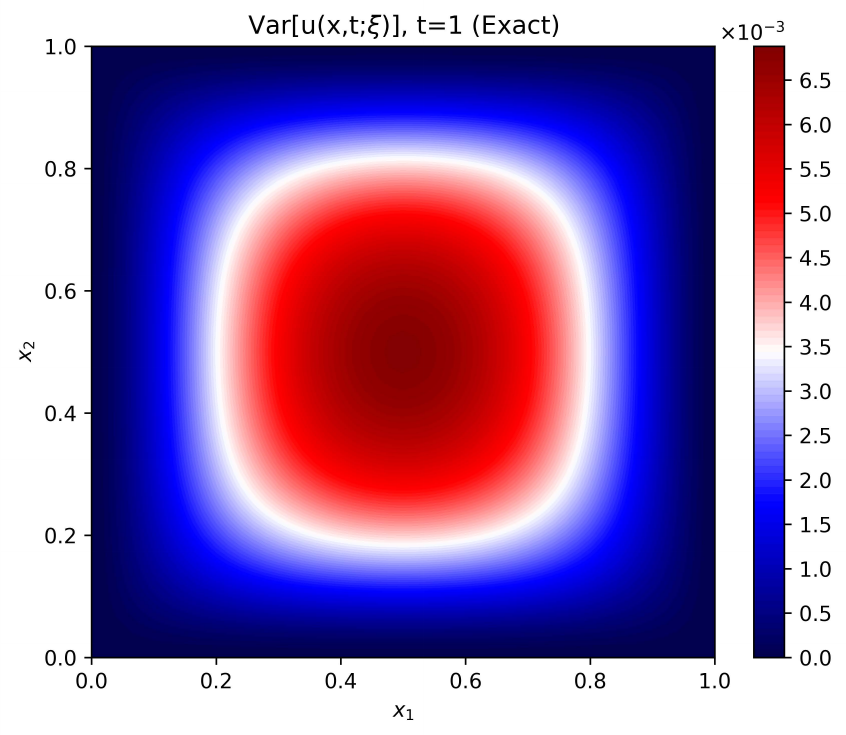}}\hfill
\subfloat[]{\includegraphics[width=0.32\linewidth]{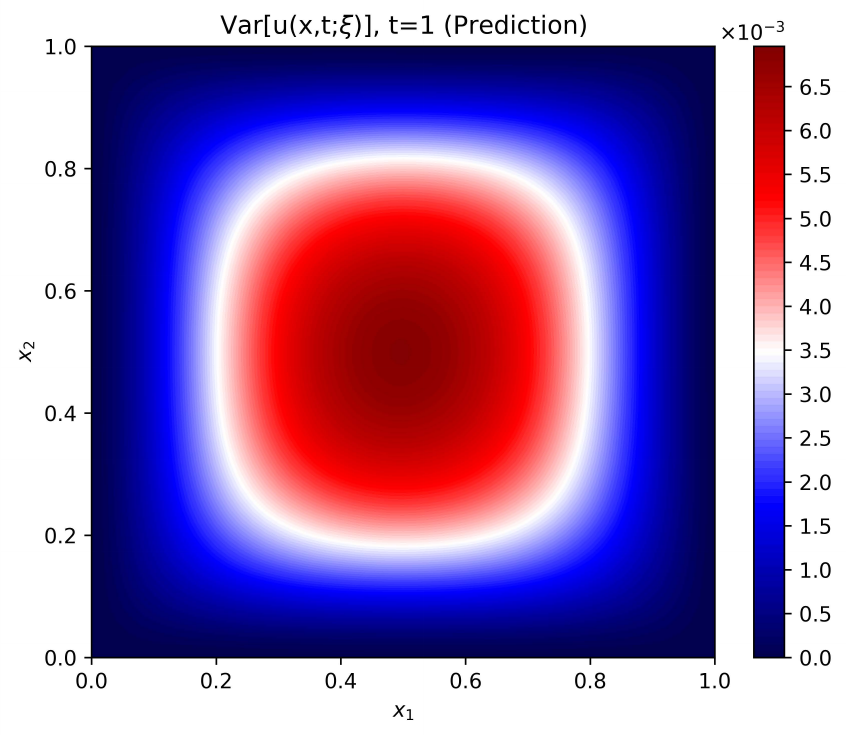}}\hfill
\subfloat[]{\includegraphics[width=0.32\linewidth]{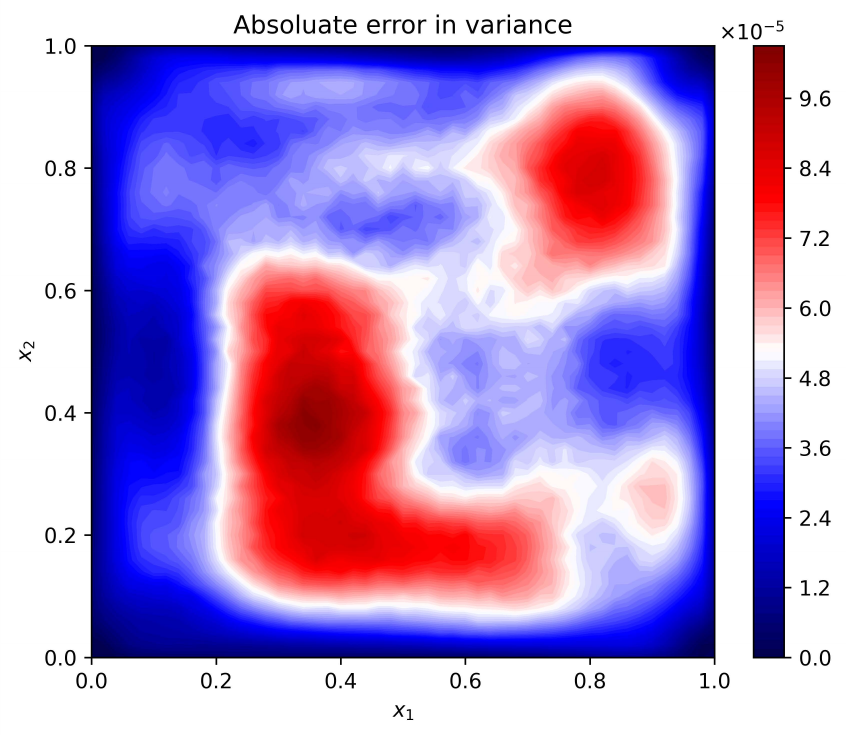}}

\caption{2D stochastic fractional reaction-diffusion equation, $\alpha=1.8$. (a): BO-fPINN method predictive mean. (b): Reference solution mean calculated by the QMC method. (c): Absolute error between reference mean and the BO-fPINN predictive mean. (e): BO-fPINN method predictive variance. (b): Reference solution variance calculated by QMC method. (c): Absolute error between reference variance and the BO-fPINN predictive variance.
}\label{fig:2D-mean-and-variance}
\end{figure}

\begin{figure}[!htp]
	\centering	
    \subfloat[]{\includegraphics[width=0.42\textwidth]{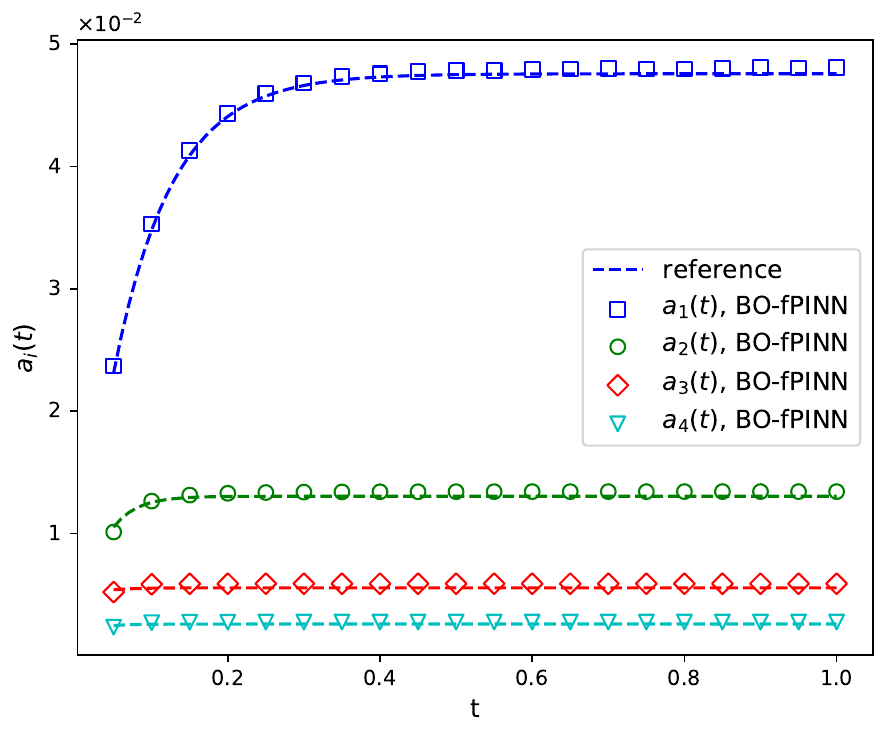}}\qquad
	\subfloat[]{\includegraphics[width=0.42\textwidth]{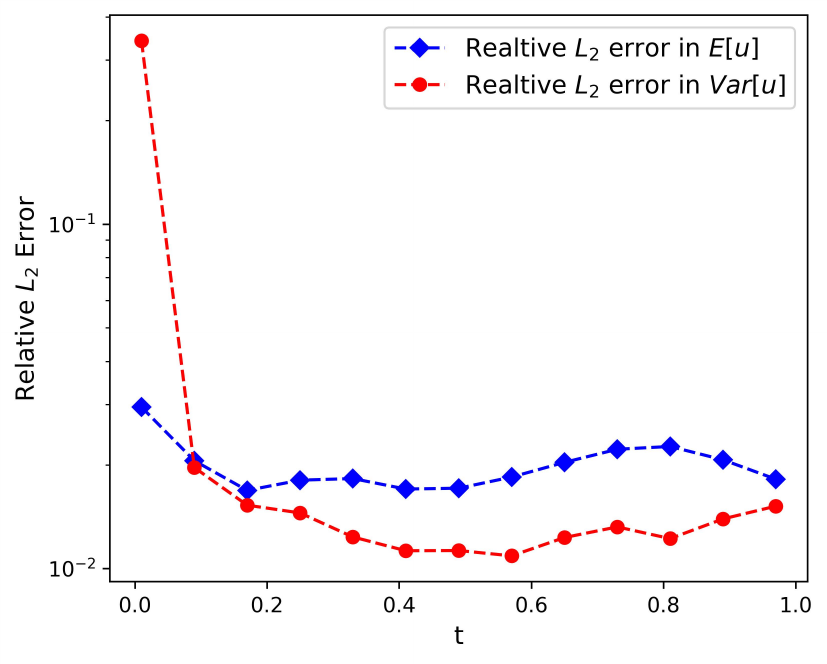}} \\	
	\caption{2D stochastic fractional reaction-diffusion equation, $\alpha=1.8$. (a):~ The evolution of the scaling factors $a_i$ as time evolves. (b):~ The relative $L_2$ error in the mean and variance obtanied by the BO-fPINN method. The reference mean and variance are calculated from the QMC method.}\label{fig:2D-ai}
\end{figure}

\begin{figure}[!ht]
\centering
\subfloat[]{\includegraphics[width=0.82\linewidth]{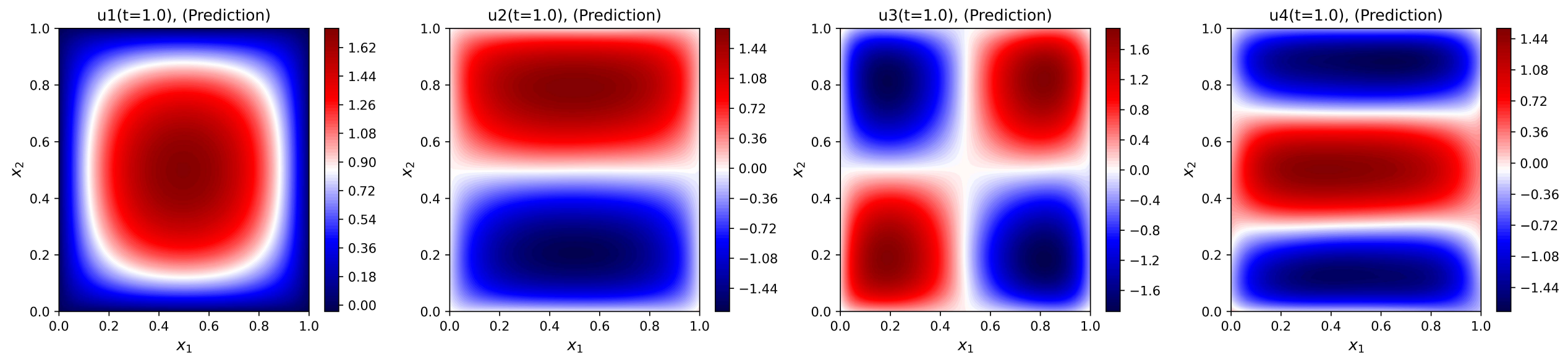}}\\
\subfloat[]{\includegraphics[width=0.82\linewidth]{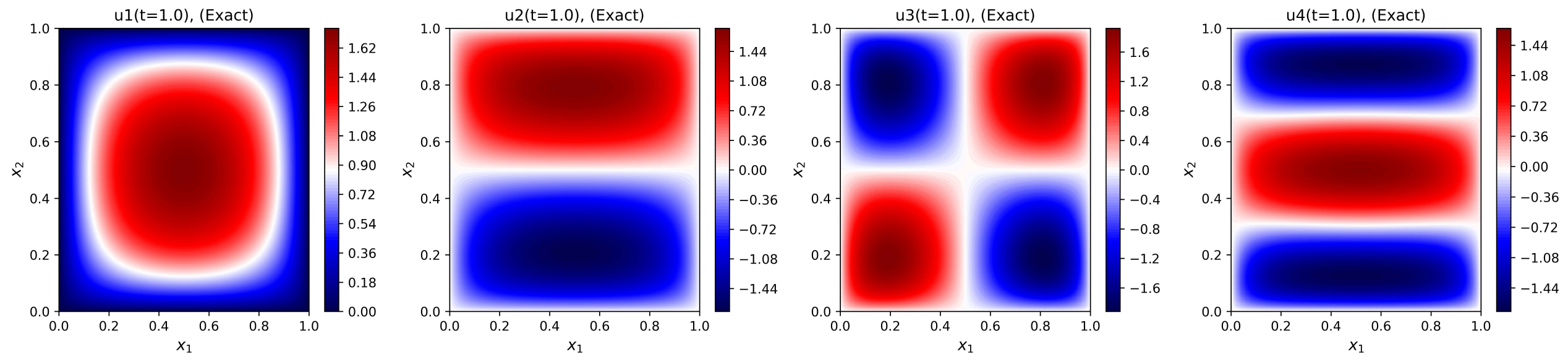}}\\
\subfloat[]{\includegraphics[width=0.82\linewidth]{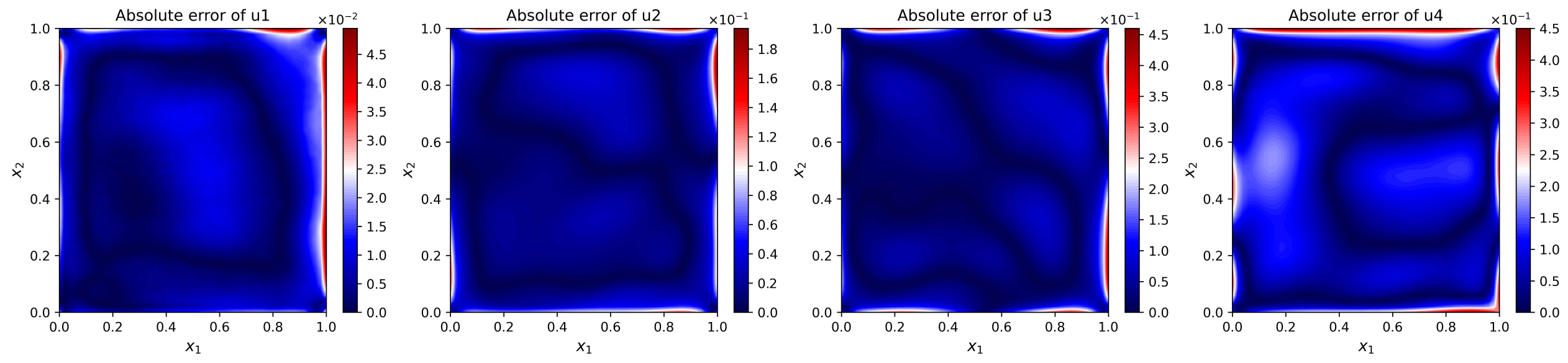}}
\caption{2D stochastic fractional reaction-diffusion equation, $\alpha=1.8$.  (a):~ The BO-fPINN spatial basis $u_i$ at $t=1$. (b):~ The reference spatial modes $u_i$ at $t=1$. (c):~ The absolute error in the spatial modes $u_i$ obtained by BO-fPINN method at $t=1$. The reference solutions are calculated by the QMC method.
}\label{fig:2D-ui-14}
\end{figure}

To obtain the reference solution, a Quasi Monte Carlo method with 1000 samples is employed to solve the SFPDE. For each sample point, the deterministic fractional PDE is numerically solved by discretizing the spatial fractional derivative using GL formula while CN-ADI scheme is implemented for temporal discretization. The predicted solution mean and variance as well as the absolute error are presented in Figure \ref{fig:2D-mean-and-variance}. Figure \ref{fig:2D-ai}(a) shows the evolution of the scaling factors $a_i$, and Figure \ref{fig:2D-ai}(b) plots the  the relative $L_2$ error of the mean and variance. Figure \ref{fig:2D-ui-14} shows the predicted BO modes and the absolute errors at $t=1$. We observe good accuracy of the BO-fPINN method.

\section{Summary}
We have presented the BO-fPINN method that combines the bi-orthogonal decomposition for stochastic processes with the physics-informed neural network for fractional PDEs. We demonstrated the performance of BO-fPINN method via several fractional partial differential equations that include a 1D nonlinear stochastic fractional reaction-diffusion equation with manufactured solution, high-dimensional random inputs/noisy initial conditions, and 2D stochastic fractional reaction-diffusion equation. The 1D nonlinear stochastic fractional reaction-diffusion equation with high-dimensional random inputs includes two cases of both time-independent random forcing and time-evolving random forcing. Compared to traditional numerical methods, BO-fPINN has an advantage that it can solve the inverse problem efficiently with the identical formulation and code as for the forward problems. We also applied transfer learning for the fractional order to accelerate the training of BO-fPINN method. Taken together, BO-fPINN seems to be an effective method for general time-dependent SFPDEs, especially for long-time integration for which existing traditional solvers are limited or are computationally  prohibitive. A current limitation is that we still employ auxiliary grids for computing the fractional derivatives numerically but in future work we plan to replace this approach with efficient sampling methods~\cite{guo2022monte}.

\section*{Acknowledgments}
Ling Guo is supported by the NSF of China (12071301, 92270115) and the Shanghai Municipal Science and Technology Commission (20JC1412500).
Fanhai Zeng is supported by the National Key R\&D Program of China (2021YFA1000202, 2021YFA1000200),
the NSF of China (12171283, 12120101001),
the startup fund from Shandong University (11140082063130),
and the Science Foundation Program for Distinguished Young Scholars of Shandong (Overseas) (2022HWYQ-045).

\FloatBarrier


\begin{thebibliography}{10}
\expandafter\ifx\csname url\endcsname\relax
  \def\url#1{\texttt{#1}}\fi
\expandafter\ifx\csname urlprefix\endcsname\relax\def\urlprefix{URL }\fi
\expandafter\ifx\csname href\endcsname\relax
  \def\href#1#2{#2} \def\path#1{#1}\fi

\bibitem{benson2000application}
D.~A. Benson, S.~W. Wheatcraft, M.~M. Meerschaert, Application of a fractional
  advection-dispersion equation, Water Resources Research 36~(6) (2000)
  1403--1412.

\bibitem{mainardi2010fractional}
F.~Mainardi, Fractional calculus and waves in linear viscoelasticity: An
  introduction to mathematical models, World Scientific, 2010.

\bibitem{song2016fractional}
F.~Song, C.~Xu, G.~E. Karniadakis, A fractional phase-field model for two-phase
  flows with tunable sharpness: Algorithms and simulations, Computer Methods in
  Applied Mechanics and Engineering 100~(305) (2016) 376--404.

\bibitem{yu2016fractional}
Y.~Yu, P.~Perdikaris, G.~E. Karniadakis, Fractional modeling of viscoelasticity
  in 3d cerebral arteries and aneurysms, Journal of Computational Physics 323
  (2016) 219--242.

\bibitem{lischke2020fractional}
A.~Lischke, G.~Pang, M.~Gulian, F.~Song, C.~Glusa, X.~Zheng, Z.~Mao, W.~Cai,
  M.~M. Meerschaert, M.~Ainsworth, et~al., What is the fractional {Laplacian}?
  a comparative review with new results, Journal of Computational Physics 404
  (2020) 109009.

\bibitem{tao2015recent}
T.~Tang, T.~Zhou, Recent developments in high order numerical methods for
  uncertainty quantification, Scientia Sinica Mathematica 45~(7) (2015)
  891--928.

\bibitem{liu2018quasi}
W.~Liu, M.~R\"{o}ckner, J.~L. da~Silva, Quasi-linear (stochastic) partial
  differential equations with time-fractional derivatives, SIAM Journal on
  Mathematical Analysis 50~(3) (2018) 2588--2607.

\bibitem{liu2018fourier}
F.~Liu, Y.~Yan, M.~Khan, Fourier spectral methods for stochastic space
  fractional partial differential equations driven by special additive noises,
  EudoxusPress.

\bibitem{yokoyama2016regularity}
S.~Yokoyama, Regularity for the solution of a stochastic partial differential
  equation with the fractional {Laplacian}, in: Mathematical Fluid Dynamics,
  Present and Future, Springer, 2016, pp. 597--613.

\bibitem{kuo2012quasi}
F.~Y. Kuo, C.~Schwab, I.~H. Sloan, Quasi-monte carlo finite element methods for
  a class of elliptic partial differential equations with random coefficients,
  SIAM Journal on Numerical Analysis 50~(6) (2012) 3351--3374.

\bibitem{kuo2016application}
F.~Y. Kuo, D.~Nuyens, Application of quasi-monte carlo methods to elliptic pdes
  with random diffusion coefficients: a survey of analysis and implementation,
  Foundations of Computational Mathematics 16~(6) (2016) 1631--1696.

\bibitem{xiu2002wiener}
D.~Xiu, G.~E. Karniadakis, The wiener--askey polynomial chaos for stochastic
  differential equations, SIAM Journal on Scientific Computing 24~(2) (2002)
  619--644.

\bibitem{xiu2010numerical}
D.~Xiu, Numerical methods for stochastic computations, in: Numerical Methods
  for Stochastic Computations, Princeton University Press, 2010.

\bibitem{guo2020constructing}
L.~Guo, A.~Narayan, T.~Zhou, Constructing least-squares polynomial
  approximations, SIAM Review 62~(2) (2020) 483--508.

\bibitem{sapsis2009dynamically}
T.~P. Sapsis, P.~F. Lermusiaux, Dynamically orthogonal field equations for
  continuous stochastic dynamical systems, Physica D: Nonlinear Phenomena
  238~(23-24) (2009) 2347--2360.

\bibitem{sapsis2011dynamically}
T.~P. Sapsis, Dynamically orthogonal field equations for stochastic fluid flows
  and particle dynamics, Ph.D. thesis, Massachusetts Institute of Technology
  (2011).

\bibitem{sapsis2012dynamical}
T.~P. Sapsis, P.~F. Lermusiaux, Dynamical criteria for the evolution of the
  stochastic dimensionality in flows with uncertainty, Physica D: Nonlinear
  Phenomena 241~(1) (2012) 60--76.

\bibitem{cheng2013dynamically1}
M.~Cheng, T.~Y. Hou, Z.~Zhang, A dynamically bi-orthogonal method for
  time-dependent stochastic partial differential equations i: Derivation and
  algorithms, Journal of Computational Physics 242 (2013) 843--868.

\bibitem{cheng2013dynamically2}
M.~Cheng, T.~Y. Hou, Z.~Zhang, A dynamically bi-orthogonal method for
  time-dependent stochastic partial differential equations ii: Adaptivity and
  generalizations, Journal of Computational Physics 242 (2013) 753--776.

\bibitem{musharbash2015error}
E.~Musharbash, F.~Nobile, T.~Zhou, Error analysis of the dynamically orthogonal
  approximation of time dependent random pdes, SIAM Journal on Scientific
  Computing 37~(2) (2015) A776--A810.

\bibitem{choi2014equivalence}
M.~Choi, T.~P. Sapsis, G.~E. Karniadakis, On the equivalence of dynamically
  orthogonal and bi-orthogonal methods: Theory and numerical simulations,
  Journal of Computational Physics 270 (2014) 1--20.

\bibitem{babaee2017robust}
H.~Babaee, M.~Choi, T.~P. Sapsis, G.~E. Karniadakis, A robust
  bi-orthogonal/dynamically-orthogonal method using the covariance
  pseudo-inverse with application to stochastic flow problems, Journal of
  Computational Physics 344 (2017) 303--319.

\bibitem{graepel2003solving}
T.~Graepel, Solving noisy linear operator equations by {Gaussian} processes:
  Application to ordinary and partial differential equations, in: ICML, Vol.~3,
  2003, pp. 234--241.

\bibitem{sarkka2011linear}
S.~S{\"a}rkk{\"a}, Linear operators and stochastic partial differential
  equations in {Gaussian} process regression, in: International Conference on
  Artificial Neural Networks, Springer, 2011, pp. 151--158.

\bibitem{bilionis2016probabilistic}
I.~Bilionis, Probabilistic solvers for partial differential equations, ArXiv
  e-prints (2016) ArXiv--1607.

\bibitem{raissi2018numerical}
M.~Raissi, P.~Perdikaris, G.~E. Karniadakis, Numerical {Gaussian} processes for
  time-dependent and nonlinear partial differential equations, SIAM Journal on
  Scientific Computing 40~(1) (2018) A172--A198.

\bibitem{lagaris1998artificial}
I.~E. Lagaris, A.~Likas, D.~I. Fotiadis, Artificial neural networks for solving
  ordinary and partial differential equations, IEEE Transactions on Neural
  Networks 9~(5) (1998) 987--1000.

\bibitem{lagaris2000neural}
I.~E. Lagaris, A.~C. Likas, D.~G. Papageorgiou, Neural-network methods for
  boundary value problems with irregular boundaries, IEEE Transactions on
  Neural Networks 11~(5) (2000) 1041--1049.

\bibitem{khoo2021solving}
Y.~Khoo, J.~Lu, L.~Ying, Solving parametric pde problems with artificial neural
  networks, European Journal of Applied Mathematics 32~(3) (2021) 421--435.

\bibitem{raissi2017physics}
M.~Raissi, P.~Perdikaris, G.~E. Karniadakis, Physics informed deep learning
  (part i): Data-driven solutions of nonlinear partial differential equations,
  ArXiv Preprint ArXiv:1711.10561.

\bibitem{chen1993approximations}
T.~Chen, H.~Chen, Approximations of continuous functionals by neural networks
  with application to dynamic systems, IEEE Transactions on Neural networks
  4~(6) (1993) 910--918.

\bibitem{chen1995universal}
T.~Chen, H.~Chen, Universal approximation to nonlinear operators by neural
  networks with arbitrary activation functions and its application to dynamical
  systems, IEEE Transactions on Neural Networks 6~(4) (1995) 911--917.

\bibitem{Raissi2017PhysicsID}
M.~Raissi, P.~Perdikaris, G.~E. Karniadakis, Physics informed deep learning
  (part ii): Data-driven discovery of nonlinear partial differential equations,
  ArXiv Abs/1711.10566.

\bibitem{raissi2019deep}
M.~Raissi, Z.~Wang, M.~S. Triantafyllou, G.~E. Karniadakis, Deep learning of
  vortex-induced vibrations, Journal of Fluid Mechanics 861 (2019) 119--137.

\bibitem{tartakovsky2020physics}
A.~M. Tartakovsky, C.~O. Marrero, P.~Perdikaris, G.~D. Tartakovsky,
  D.~Barajas-Solano, Physics-informed deep neural networks for learning
  parameters and constitutive relationships in subsurface flow problems, Water
  Resources Research 56~(5) (2020) e2019WR026731.

\bibitem{kissas2020machine}
G.~Kissas, Y.~Yang, E.~Hwuang, W.~R. Witschey, J.~A. Detre, P.~Perdikaris,
  Machine learning in cardiovascular flows modeling: Predicting arterial blood
  pressure from non-invasive 4d flow mri data using physics-informed neural
  networks, Computer Methods in Applied Mechanics and Engineering 358 (2020)
  112623.

\bibitem{gao2022failure}
Z.~Gao, L.~Yan, T.~Zhou, Failure-informed adaptive sampling for pinns, arXiv
  preprint arXiv:2210.00279.

\bibitem{gao2023failure}
Z.~Gao, T.~Tang, L.~Yan, T.~Zhou, Failure-informed adaptive sampling for pinns,
  part ii: combining with re-sampling and subset simulation, arXiv preprint
  arXiv:2302.01529.

\bibitem{pang2019fpinns}
G.~Pang, L.~Lu, G.~E. Karniadakis, fpinns: Fractional physics-informed neural
  networks, SIAM Journal on Scientific Computing 41~(4) (2019) A2603--A2626.

\bibitem{han2017deep}
J.~Han, A.~Jentzen, et~al., Deep learning-based numerical methods for
  high-dimensional parabolic partial differential equations and backward
  stochastic differential equations, Communications in Mathematics and
  Statistics 5~(4) (2017) 349--380.

\bibitem{raissi2018forward}
M.~Raissi, Forward-backward stochastic neural networks: Deep learning of
  high-dimensional partial differential equations, ArXiv Preprint
  ArXiv:1804.07010.

\bibitem{zhu2018bayesian}
Y.~Zhu, N.~Zabaras, Bayesian deep convolutional encoder--decoder networks for
  surrogate modeling and uncertainty quantification, Journal of Computational
  Physics 366 (2018) 415--447.

\bibitem{tartakovsky2018learning}
A.~M. Tartakovsky, C.~O. Marrero, P.~Perdikaris, G.~D. Tartakovsky,
  D.~Barajas-Solano, Learning parameters and constitutive relationships with
  physics informed deep neural networks, arXiv preprint arXiv:1808.03398.

\bibitem{yan2019adaptive}
L.~Yan, T.~Zhou, An adaptive surrogate modeling based on deep neural networks
  for large-scale bayesian inverse problems, arXiv preprint arXiv:1911.08926.

\bibitem{zhang2019quantifying}
D.~Zhang, L.~Lu, L.~Guo, G.~E. Karniadakis, Quantifying total uncertainty in
  physics-informed neural networks for solving forward and inverse stochastic
  problems, Journal of Computational Physics 397 (2019) 108850.

\bibitem{zhang2020learning}
D.~Zhang, L.~Guo, G.~E. Karniadakis, Learning in modal space: Solving
  time-dependent stochastic pdes using physics-informed neural networks, SIAM
  Journal on Scientific Computing 42~(2) (2020) A639--A665.

\bibitem{zhao2021spectral}
Y.~Zhao, Z.~Mao, L.~Guo, Y.~Tang, G.~E. Karniadakis, A spectral method for
  stochastic fractional pdes using dynamically-orthogonal/bi-orthogonal
  decomposition, Journal of Computational Physics (2022) 111213.

\bibitem{wang2021understanding}
S.~Wang, Y.~Teng, P.~Perdikaris, Understanding and mitigating gradient flow
  pathologies in physics-informed neural networks, SIAM Journal on Scientific
  Computing 43~(5) (2021) A3055--A3081.

\bibitem{jin2021nsfnets}
X.~Jin, S.~Cai, H.~Li, G.~E. Karniadakis, Nsfnets (navier-stokes flow nets):
  Physics-informed neural networks for the incompressible navier-stokes
  equations, Journal of Computational Physics 426 (2021) 109951.

\bibitem{goswami2019transfer}
S.~Goswami, C.~Anitescu, S.~Chakraborty, T.~Rabczuk, Transfer learning enhanced
  physics informed neural network for phase-field modeling of fracture,
  Theoretical and Applied Fracture Mechanics 106 (2020) 102447.

\bibitem{guo2022monte}
L.~Guo, H.~Wu, X.~Yu, T.~Zhou, Monte carlo fpinns: Deep learning method for
  forward and inverse problems involving high dimensional fractional partial
  differential equations, Computer Methods in Applied Mechanics and Engineering
  400 (2022) 115523.

\bibitem{li2019theory}
C.~Li, M.~Cai, Theory and numerical approximations of fractional integrals and
  derivatives, SIAM, 2019.

\bibitem{zhao2015series}
L.~Zhao, W.~Deng, A series of high-order quasi-compact schemes for space
  fractional diffusion equations based on the superconvergent approximations
  for fractional derivatives, Numerical Methods for Partial Differential
  Equations 31~(5) (2015) 1345--1381.

\end{thebibliography}

\newpage
\appendix
\setcounter{remark}{0}
\setcounter{theorem}{0}
\setcounter{definition}{0}
\setcounter{figure}{0}
\setcounter{table}{0}
\renewcommand{\thetheorem}{\Alph{section}.\arabic{theorem}}
\renewcommand{\thedefinition}{\Alph{section}.\arabic{definition}}
\renewcommand{\theremark}{\Alph{section}.\arabic{remark}}
\renewcommand{\thefigure}{\Alph{section}.\arabic{figure}}
\renewcommand{\thetable}{\Alph{section}.\arabic{table}}

\section{Definition of fractional derivatives and the GL approximation}\label{appendix:fractional-derivative-def}
\begin{definition}[Riemann-Liouville fractional derivative \cite{li2019theory}] 
The $\alpha$-order left and right Riemann-Liouville derivatives of the function $u(x)$,  $x\in [a,b]$, are defined by
\begin{align}\label{eqn:RL_rg}
  \sideset{_{RL}}{_{a,x}^{\alpha}}{\mathop{D}}u(x) &=\frac{1}{\Gamma(n-\alpha)}\left [\frac{d^n}{dx^n}\int_{a}^{x}(x-s)^{n-\alpha -1}u(s)ds\right],\\
  \sideset{_{RL}}{_{x,b}^{\alpha}}{\mathop{D}}u(x) &=\frac{1}{\Gamma(n-\alpha)}\left [\frac{d^n}{dx^n}\int_{x}^{b}(s-x)^{n-\alpha -1}u(s)ds\right],
\end{align}
where~$n-1<\alpha\leq n$, and $n$ is a positive integer.
\end{definition}

\begin{definition}[Riesz derivative \cite{li2019theory}]\label{def:Riesz}
The $\alpha$-order Riesz derivative of the function $u(x)$ on $[a,b]$, is defined by
\begin{align}
  D_{|x|}^{\alpha}u(x)=-\frac{1}{2\cos(\frac{\pi \alpha}{2})}(\sideset{_{RL}}{_{a,x}^{\alpha}}{\mathop{D}}u(x)
+\sideset{_{RL}}{_{x,b}^{\alpha}}{\mathop{D}}u(x)),
\end{align}
where~$n-1<\alpha\leq n$, $n$ is a positive integer, $\sideset{_ {RL}}{_{a,x}^{\alpha}}{\mathop{D}}$ and $\sideset{_{RL}}{_{x,b}^{\alpha}}{\mathop{D}}$ are left and right Riemann-Liouville derivatives.
\end{definition}

\subsection{Gr\"{u}nwald-Letnikov (GL) formula for fractional derivatives}\label{appendix:GL}
Based on uniform grid $x^j=a+j\Delta x$ for $j=0,1,2,\cdots, N$, the first- and second-order shifted GL finite difference operator for  approximating $D_{|x|}^{\alpha}u(x)$ can be defined as \cite{zhao2015series}:
\begin{align}\label{eqn:RL_rg-2}
D_{|x|}^{\alpha}u(x)&\approx  c_{\alpha}\delta_{\Delta x,1}^{\alpha}u(x^j),\\
D_{|x|}^{\alpha}u(x)&\approx  -\frac{1}{2\cos(\frac{\pi\alpha}{2})}
\left[\frac{\alpha}{2}\delta_{\Delta x,1}^{\alpha}u(x^j)+\left(1-\frac{\alpha}{2}\right)\delta_{\Delta x,0}^{\alpha}u(x^j)\right],
\quad x = x^j,
\end{align}
where  
\begin{align*}
  \delta_{\Delta x,p}^{\alpha}u(x^j):=\frac{1}{\Delta x^{\alpha}}\left[\sum_{k=0}^{j}w_k^{\alpha}u(x^j-(k-p)\Delta x)+\sum_{k=0}^{N-j}w_k^{\alpha}u(x^j+(k-p)\Delta x)\right],
\end{align*}
with  $w_0^{\alpha}=1$,$w_k^{\alpha}=(1-\frac{\alpha+1}{k})w_{k-1}^{\alpha}$, for $k\geq 1$.

\section{Generalized Karhunen-Lo\`eve (KL) expansion}\label{appendix:KL}
For a time-dependent random field $u(x,t;\omega)$, the generalized KL expansion at a given time $t$ is
\begin{equation}\label{eqn:kl}
    u(x,t;\omega) = \overbar{u}(x,t) + \sum_{i=1}^{\infty}\sqrt{\lambda_i} \phi_i(x,t) \xi_{i}(t;\omega), \quad \omega \in \Omega,
\end{equation}
where $\overbar{u}$ is the mean, $\xi_{i}(t;\omega)$ ($i=1,2,3,...$) are zero-mean independent random variables, $\lambda_i$ and $\phi_i$ are the $i^\text{th}$ largest eigenvalue and the corresponding eigenfunction of the covariance kernel, i,e., they solve the following eigenproblem:
\begin{equation}
    \int_{\D} C_{u(x,t)u(y,t)}\phi_i(y,t) \dd y = \lambda_i \phi_i(x, t).
\end{equation}
Here $C_{u(x,t)u(y,t)} = \E \left[ (u(x,t;\omega)-\overbar{u}(x,t))^{\mathrm{T}} (u(y,t;\omega)-\overbar{u}(y,t)) \right]$ is the covariance kernel of $u$~\cite{sapsis2009dynamically}.

\section{Numerical results for fPINN method}\label{appendix:pinn-gl}
\begin{figure}[!htp]
\centering	
\subfloat[]{\includegraphics[width=0.42\textwidth]{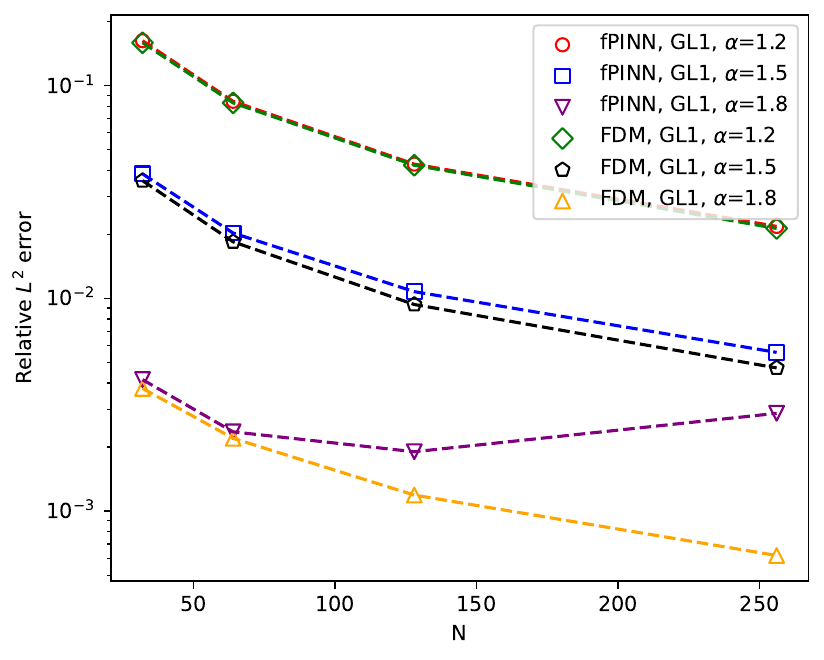}}\qquad
\subfloat[]{\includegraphics[width=0.42\textwidth]     {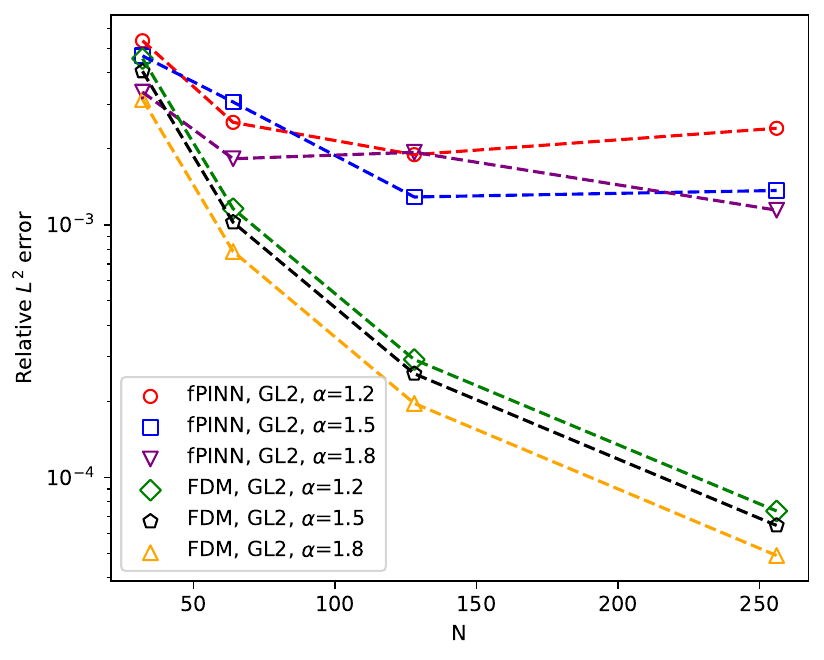}} \\
\caption{Comparison diagram of fPINN method and finite difference method for solving nonlinear fractional partial differential equations. (a):~ The relative $L^{2}$ error between the exact solution and the solutions of fPINN method and finite difference method for different values of the fractional order $\alpha$ when we use the first order GL formula. (b):~ the relative $L^{2}$ error between the exact solution and the solutions of fPINN method and finite difference method for different values of the fractional order $\alpha$ when we use the second order GL formula.}\label{fig:pinn-fdm-gl}
\end{figure}

To test the convergence of the fPINN method for deterministic fractional PDEs, we consider the following 1D nonlinear time-dependent fractional PDE:
\begin{equation}\label{prob:nonlinear-fpde}
  \begin{cases}
    u_t-D_{|x|}^\alpha u-u(1-u^2) = g(x, t), \quad x\in (0,1),\quad t\in [0,T], \\
    u(x,0) = 100(1-x)^3x^3, \quad x\in (0,1)\\
    u(0,t)=u(1,t) = 0.
  \end{cases}
\end{equation}
We consider a fabricated solution $u=100e^{-t}(1-x)^3x^3$. The corresponding forcing term is
\begin{align*}
  g(x,t) = & -200e^{-t}(1-x)^3x^3 - 10^{6}e^{-3t}(1-x)^3x^3 \\
  & + \frac{1}{2cos\left(\frac{\pi \alpha}{2}\right)}\left[\frac{\Gamma(4)}{\Gamma(4-\alpha)}\left(x^{3-\alpha} +(1-x)^{3-\alpha}\right) - 3\frac{\Gamma(5)}{\Gamma(5-\alpha)}\left(x^{4-\alpha} +(1-x)^{4-\alpha}\right) \right.\\
   &\left. + 3\frac{\Gamma(6)}{\Gamma(6-\alpha)}\left(x^{5-\alpha} +(1-x)^{5-\alpha}\right) - \frac{\Gamma(7)}{\Gamma(7-\alpha)}\left(x^{6-\alpha} +(1-x)^{6-\alpha}\right) \right].\\
\end{align*}
Following \cite{pang2019fpinns}, we use the first-order/second-order shifted GL formula to compute the fractional derivative of the DNN outputs, the DNN surrogate has 4 hidden layers with 20 neurons in each layer and the activation function is chosen as the tanh function. For the FDM scheme, we set $\delta t=0.001$, $\Delta x= \frac{1}{N}$, where $N=[32,64,128,512]$. In Figure \ref{fig:pinn-fdm-gl}, we compare the solutions of fPINNs against FDM for different fractional derivative orders.
We observe that for a smaller number of training points
the fPINN yields solution accuracy very close to that of the FDM. But the error of the fPINN saturates as the number of trianing points increasing, which is due to the optimization error of the nueral netwroks.

\section{Classical hybrid QMC-BO method for time-dependent SFPDEs }\label{appendix:QMC-BO}

Based on the expansion (\ref{eqn:BOexpansion}) and the BO constraints, we can derive a set of independent and explicit equations for all the unknowns quantities. Here we state the BO evolution equations for the SFPDE \eqref{eqn:SFPDE} without proof. For more details of the derivation, the authors can see \cite{cheng2013dynamically2,zhang2020learning, zhao2021spectral} and references therein.

Define the matrices $S$ and $M$ whose entries are
\begin{equation}
S_{ij} = \left\langle u_i, \frac{\partial u_j}{\partial t} \right\rangle, \quad M_{ij} = \mathbb{E} \bigg[ Y_i \frac{\partial Y_j}{\partial t}  \bigg].
\end{equation}
Then by taking the time derivative of $\left\langle \cdot,\cdot\right\rangle $, we have
\begin{equation}
\begin{aligned}
&S_{ij} + S_{ji} = \left\langle u_i, \frac{\partial u_j}{\partial t} \right\rangle + \left\langle \frac{\partial u_i}{\partial t}, u_j \right\rangle   &for ~i \neq j, \\
&S_{ij} = \frac{1}{2} \frac{d \lambda_{i}(t)}{dt}   &for ~i = j,\\
&M_{ij} + M_{ji} = \mathbb{E} \bigg[ Y_i \frac{\partial Y_j}{\partial t}  \bigg] + \mathbb{E} \bigg[ \frac{\partial Y_i}{\partial t} Y_{j}  \bigg] = 0.
\end{aligned}
\end{equation}

\begin{theorem}[see~\cite{cheng2013dynamically2}]
We assume that the bases and stochastic coefficients satisfy the BO condition. Then, the original SFPDE \eqref{eqn:SFPDE} is reduced to the following system of equations:
\begin{align}\label{eqn:bo-eqns3}
\begin{cases}
  &\frac{\partial \overline{u}(x,t)}{\partial t}=\mathbb{E}[\mathcal{N}_x^{\alpha,\beta}[u(x,t;\omega)]+g],\\
  & \lambda_i\frac{dY_i(t;\omega)}{dt}=-\sum_{j=1}^{N}S_{ij}Y_j + \langle\mathcal{N}_x^{\alpha,\beta}[u]+g-\mathbb{E}[\mathcal{N}_x^{\alpha,\beta}[u]+g], u_i(x,t)\rangle, \\
  &\frac{\partial u_i(x, t)}{\partial t}=-\sum_{j=1}^{N}M_{ij}u_j + \mathbb{E}[\mathcal{N}_x^{\alpha,\beta}[u]Y_i+gY_i]. 
\end{cases}
\end{align}
Moreover, if $\lambda_i \ne \lambda_j$ for $i \ne j,\ i,j,=1,2,\dots,N$, the $N$-by-$N$ matrices $S$ and $M$ have closed form expression:
\begin{equation*}
    M_{ij}=\begin{cases}
             \frac{G_{ij}+G_{ji}}{-\lambda_i+\lambda_j}, & \mbox{if } i\neq j, \\
             0, & \mbox{if} i=j,
           \end{cases}\quad
     S_{ij}=\begin{cases}
             G_{ij}+\lambda_i M_{ij}, & \mbox{if } i\neq j, \\
             G_{ii}, & \mbox{if} i=j,
           \end{cases}
\end{equation*}
where the matrix $G_{ij}=\langle \mathbb{E}[\mathcal{N}_x^{\alpha}[u]Y_j+gY_j],u_i\rangle$.
\end{theorem}
The initial condition is generated from the KL expansion of $u_0(x;\omega)$
\begin{align}\label{eqn:bo_t0}
\begin{cases}
\overline{u}(x,t_0)=\mathbb{E}[u_0(x;\omega)],\\
  u_i(x,t_0)=\sqrt{\lambda_i}v_i(x,t_0),\\
  Y_i(t_0;\omega)=\frac{1}{\sqrt{\lambda_i}}\langle u(x,t_0;\omega)-\overline{u}(x,t_0),v_i\rangle.
\end{cases}
\end{align}

For problems with random initial condition, we can numerically solve the above BO equations with the shifted GL formulas in space and a 3rd level Adams¨CBashforth scheme in time. However, for the stochastic fractional PDEs with deterministic initial condition, we need to start with a Quasi Monte Carlo method or generalized Polynomail Chaos until $t_s$ and then switch to solving the BO equations. This procedure is called the hybrid QMC-BO method and is summarized in Algorithm \ref{alg:QMC-BO} \cite{choi2014equivalence,babaee2017robust}.

\begin{algorithm}
\setlength{\baselineskip}{20pt}
\caption{~QMC-BO~ for solving time-dependent SFPDEs}\label{alg:QMC-BO}
\textbf{Step 1:} Run QMC up to $t=t_s$ from $t=0$.\\
\textbf{Step 2:} At $t=t_s$, use the KL decomposition for the solution:$$u(x,t_s;\omega)=\overline{u}(x,t_s)+\sum_{i=1}^{N}Y_i(t_s;\omega)u_i(x,t_s).$$
\textbf{Step 3:} From the KL decomposition, we can initialize:  $\overline{u}(x,t_s)$,~ $Y_i(t_s;\omega)$,~$u_i(x,t_s)$, $i=1,\cdots,N.$\\
\textbf{Step 4:} Switch over to the BO method up to time $t=T$.
\end{algorithm}

\section{Additional results of Section 4}\label{appendix:test}

\subsection{Neural network architectures used in the numerical examples}\label{appendix:network}
We list the NN architectures used in section 4 in Table \ref{tab:NN}.
\begin{table}[!htp]
\centering
\renewcommand\arraystretch{1.2}
\caption{The width and depth of the hidden layer of ~$\overline{u}_{nn}(x,t),~A_{nn}(t),~U_{nn}(x,t),~Y_{nn}(t;\xi)$ in different examples in the article. In the table, the number before the $\times$ represents the width of the hidden layer, the latter represents the depth of the hidden layer.}
\label{tab:NN}
\begin{tabular}{c|cccccc}
\hline
             &  $\overline{u}_{nn}(x,t)$  &  $A_{nn}(x,t)$  &  $U_{nn}(x,t)$  &  $Y_{nn}(x;\xi)$  \\
\hline
Section \ref{sec:RD}   & 4 $\times$ 64 & 4 $\times$ 64 & 4 $\times$ 64 & 4 $\times$ 32  \\
\hline
Section \ref{sec:RD_notime}  & 3 $\times$ 32 & 3 $\times$ 4 & 3 $\times$ 64 & 3 $\times$ 64  \\
\hline
Section \ref{sec:RD_time}  & 3 $\times$ 32 & 3 $\times$ 4 & 3 $\times$ 64 & 3 $\times$ 32  \\
\hline
Section \ref{sec:RD_2d}  & 4 $\times$ 32 & 3 $\times$ 4 & 4 $\times$ 32 & 4 $\times$ 32  \\

\hline 
\end{tabular}
\end{table}

\subsection{Additional results of Section 4.1}\label{appendix:g}
The corresponding frocing term is given as
\vspace{-8pt}
\begin{align*}
  g(x,t;\omega)\approx &  50\cos(\frac{t}{2}+\frac{\pi}{4})x^3(1-x)^3 -\sqrt{3}\cos(t)\sin(\pi x)(2\xi_1-1)\\
  &-3\sqrt{3}\sin(3t)\sin(2\pi x)(2\xi_2-1)+ \frac{50}{\cos(\frac{\pi\alpha}{2})}\sin(\frac{t}{2}+\frac{\pi}{4})\left [\frac{\Gamma(4)}{\Gamma(4-\alpha)}\left (x^{3-\alpha}+(1-x)^{3-\alpha}\right )\right.\\
   & -3\frac{\Gamma(5)}{\Gamma(5-\alpha)}\left (x^{4-\alpha}+(1-x)^{4-\alpha}\right ) +3\frac{\Gamma(6)}{\Gamma(6-\alpha)}\left (x^{5-\alpha}+(1-x)^{5-\alpha}\right ) \\
   &\left. -\frac{\Gamma(7)}{\Gamma(7-\alpha)}\left (x^{6-\alpha}+(1-x)^{6-\alpha}\right ) \right]\\
   &-\frac{\sqrt{3}}{2\cos(\frac{\pi\alpha}{2})}(1.5+\sin(t))(2\xi_1-1)\left[-\pi^2x^{2-\alpha}
   \sum_{m=1}^{M}\frac{(-1)^{m-1}}{\Gamma(4m-2-\alpha)}(\pi x)^{2m-1} \right.\\
   &\left.+ \frac{\pi x^{1-\alpha}}{\Gamma(2-\alpha)} +\sum_{m=1}^{M}\frac{(-1)^{m-1}}{\Gamma(2m-\alpha)}(\pi (1-x))^{2m-1}(1-x)^{-\alpha} \right]\\
   &+\frac{\sqrt{3}}{2\cos(\frac{\pi\alpha}{2})}(1.5+\cos(3t))(2\xi_2-1)\left[-4\pi^2x^{2-\alpha}
   \sum_{m=1}^{M}\frac{(-1)^{m-1}}{\Gamma(4m-2-\alpha)}(2\pi x)^{2m-1} \right.\\
   &\left.+ \frac{2\pi x^{1-\alpha}}{\Gamma(2-\alpha)} -\sum_{m=1}^{M}\frac{(-1)^{m-1}}{\Gamma(2m-\alpha)}(2\pi (1-x))^{2m-1}(1-x)^{-\alpha}\right ]\\
   &-u(1-u^2).
\end{align*}
Where $M$ is the number of items truncated by Taylor expansion, we choose $M=50$ in our simulations.

\subsection{Additional results of Section 4.2.1}\label{appendix:1DRDE4.2.1}
\begin{figure}[!ht]
	\centering
    \subfloat[]{\includegraphics[width=0.42\textwidth]{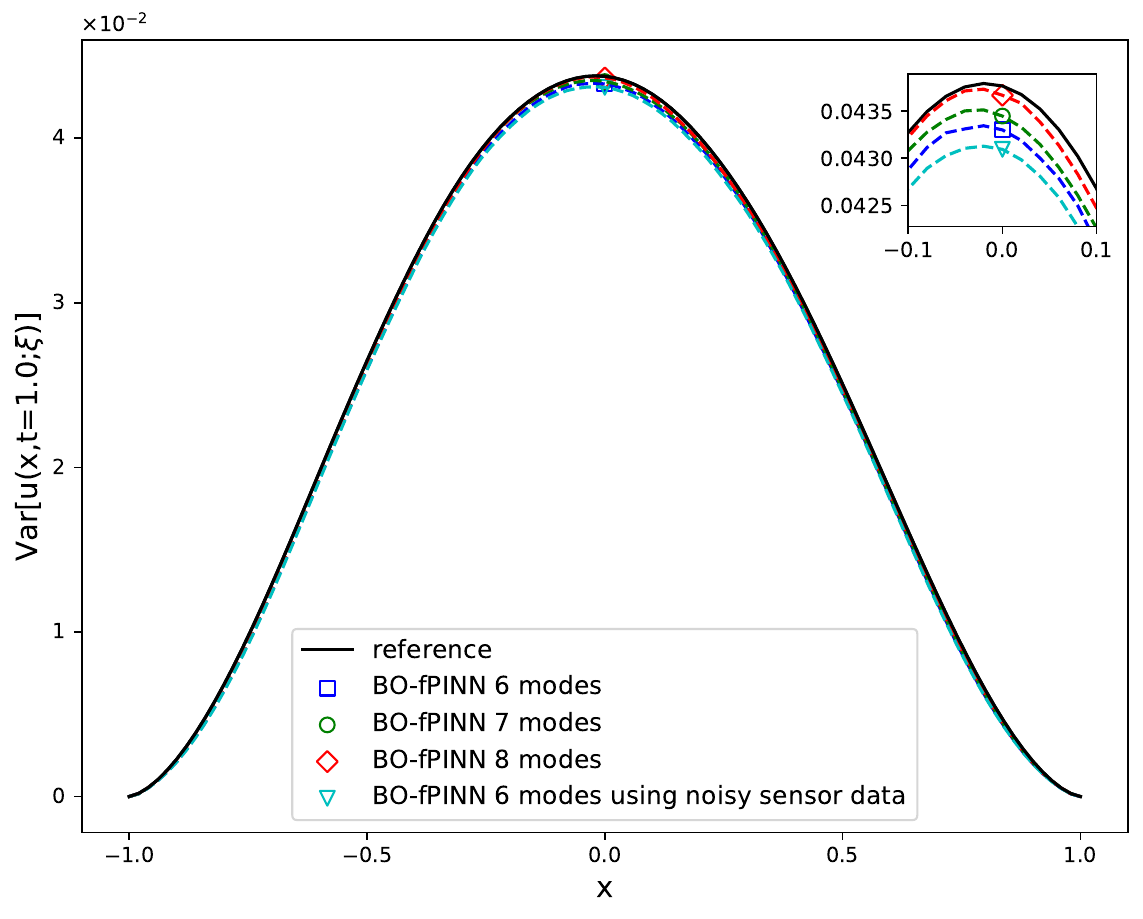}}\qquad
	\subfloat[]{\includegraphics[width=0.43\textwidth]{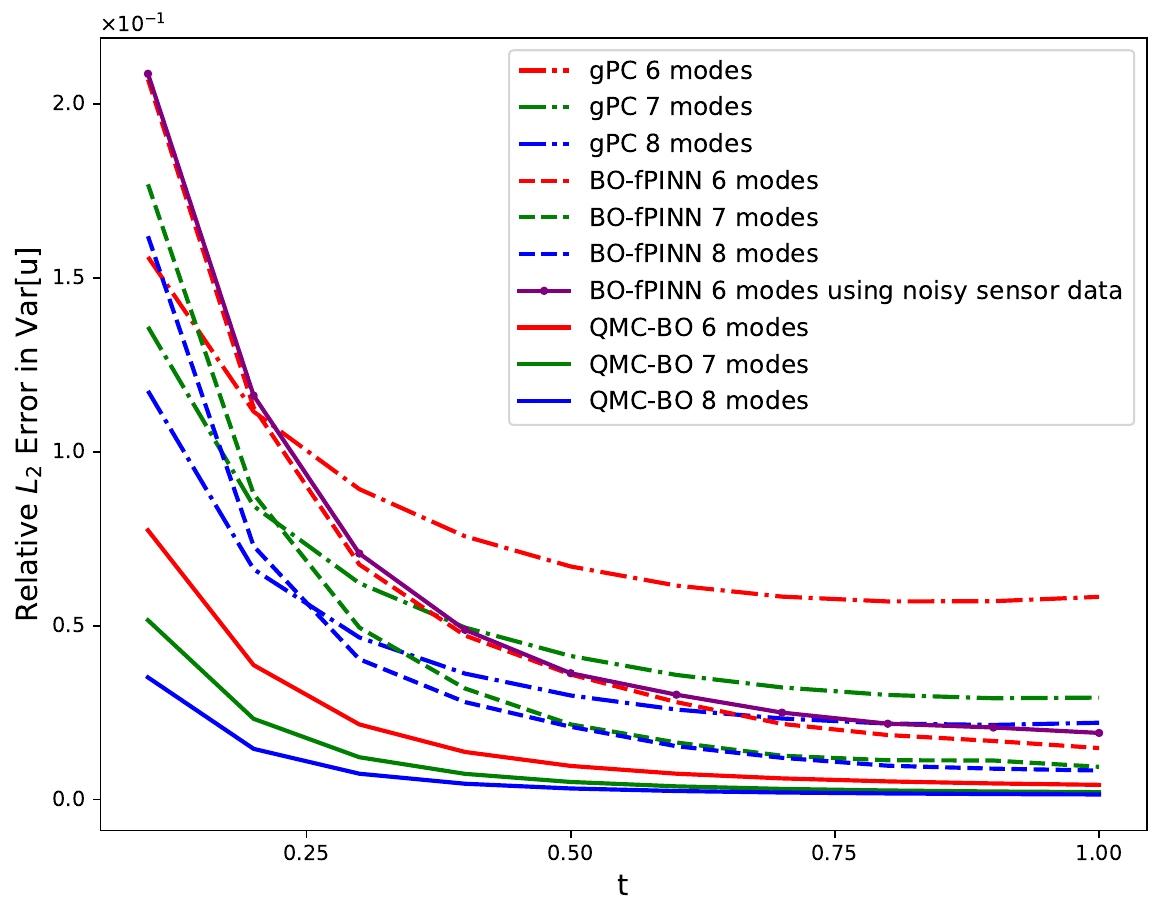}} \\	
	\caption{Forward problem with time independent random forcing. (a):~Variance of the BO-FPINN solution calculated using 6, 7, and 8 BO modes; the reference variance is calculated from the Monte Carlo simulation. (b):~The comparision of the $L_2$ erros of the solution variance calculated by the BO-FPINN method, the standard numerical BO method, and the gPC method. The gPC method generates the largest error since it fails to capture the dynamic evolution of stochastic basis for nonlinear problems. Similar accuracy is obtained with BO-fPINN method when using noisy sensor data as the initial condition. }\label{fig:1d-var-error}
\end{figure}

we analyze the effect of the number of BO expansion modes by comparing the variances of solution calculated using six, seven, and eight BO modes. Figure \ref{fig:1d-var-error}a plots the predicted variance of BO-fPINN method at time t = 1.0 versus the reference solution. Figure \ref{fig:1d-var-error}b compares the relative $L_2$ error of the solution variance at $t=1.0$ obtained via three different methods: BO-fPINN, gPC and QMC-BO. We can see that the gPC method generates the largest error as it fails to capture the evolution of the system's stochastic structure due to the non-linearity, therefore, to achieve the same accuracy, gPC method has to include a larger number of modes than the BO method. Meanwhile, we can observe that better accuracy can be achieved when more modes are included in the BO expansion. However, the BO-fPINN method is less accurate than the standard numerical BO method due to dominant optimization errors of DNNs. We obtain similar accuracy when using noisy sensor data as the initial condition and thus the BO-fPINN method is robust for noisy data. These findings are consistent with the findings in \cite{zhang2020learning}.

\subsection{Additional results of Section 4.2.2}\label{appendix:1DRDE4.2.2}
\begin{figure}[!ht]
  \centering
  \includegraphics[height=0.48\textheight, width=0.80\textwidth]{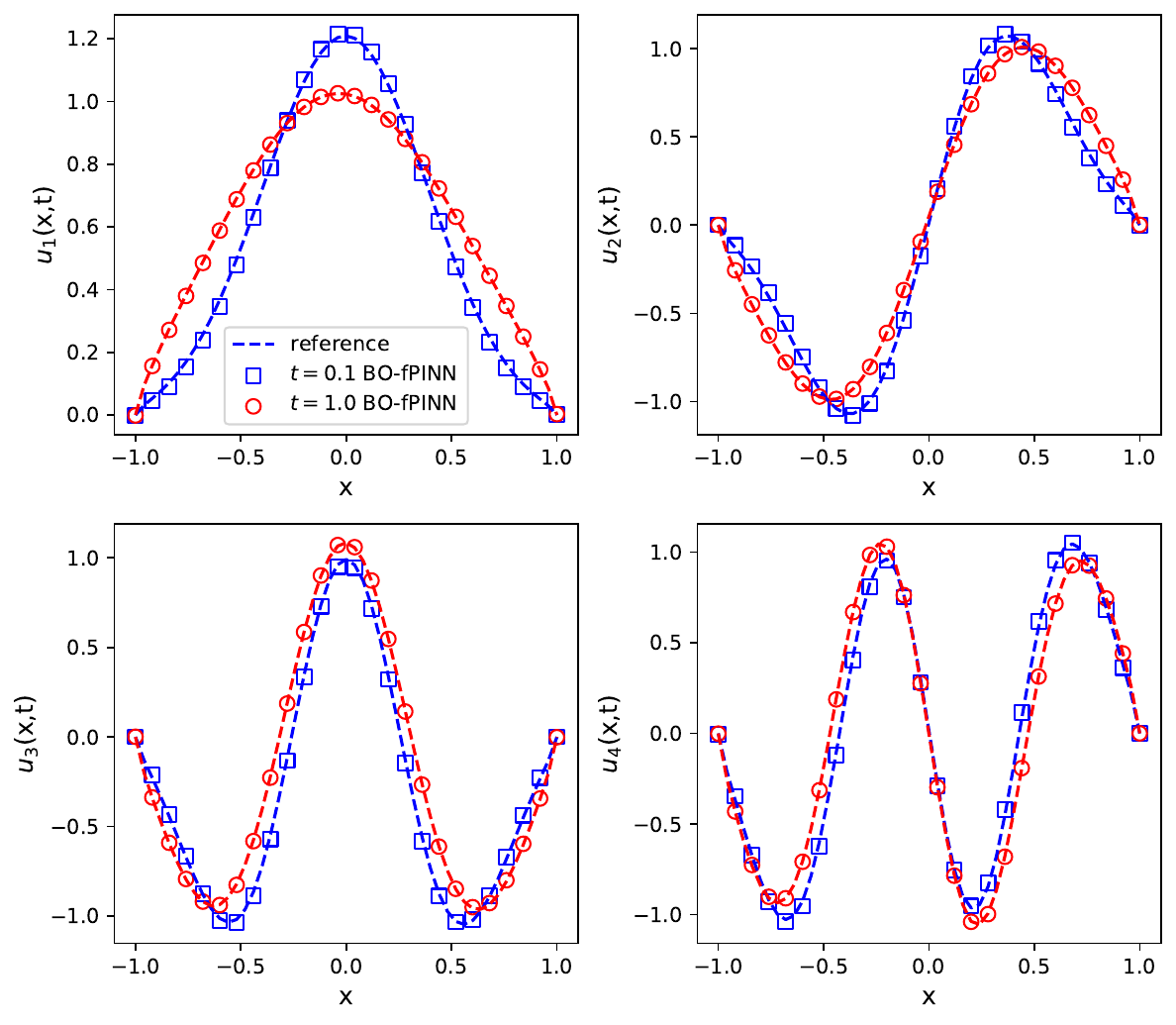}
  \caption{Inverse problem with time independent random forcing. The spatial modes at $t=0.1$ and $t=1$. We use QMC-BO solutions as reference.}
  \label{fig:ui-1d-inverse}
\end{figure}

Figure \ref{fig:ui-1d-inverse} shows the modal functions $u_i$ obtained by BO-fPINN method at $t=0.1,1$ for inverse problem with time independent random forcing. The reference solutions are calculated by solving the forward problem using the QMC-BO method. We can see that the predicted solution by BO-fPINN and the reference solution are almost completely coincident. In Table \ref{tab:inverse-error}, we summarize the relative $L_2$ errors for all of the BO components at time $t=1.0$. We conclude that the BO-fPINN method can efficiently solve the inverse problem with time independent random forcing.

\begin{table}[!htp]
\centering
\caption{Inverse problem with time independent random forcing. The relative $L_2$ errors of BO-fPINN solutions versus the reference solutions at the final time $t=1.0$. The reference solutions are calculated by solving the forward problem using the QMC-BO method. }
\label{tab:inverse-error}
\begin{tabular}{c|cccccc}
\hline
          &  $a_1$  &  $a_2$   &  $a_3$  &  $a_4$   \\

Relative $L_2$ error  & 0.436$\%$ & 0.022$\%$ & 0.471$\%$ & 0.629$\%$ \\
\hline
          &  $u_1$  &  $u_2$   &  $u_3$  &  $u_4$   \\

Relative $L_2$ error  & 0.293$\%$ & 0.343$\%$ & 0.603$\%$ & 0.542$\%$  \\
\hline
          &  $Y_1$  &  $Y_2$   &  $Y_3$  &  $Y_4$    \\

Relative $L_2$ error  & 0.0098 & 0.0140 & 0.0287 &  0.0434 \\
\hline
\end{tabular}
\end{table}

\subsection{Additional results of Section 4.2.3}\label{appendix:1DRDE4.2.3}
\begin{figure}[!ht]
	\centering
    \subfloat[]{\includegraphics[height=0.27\textheight, width=0.44\textwidth]{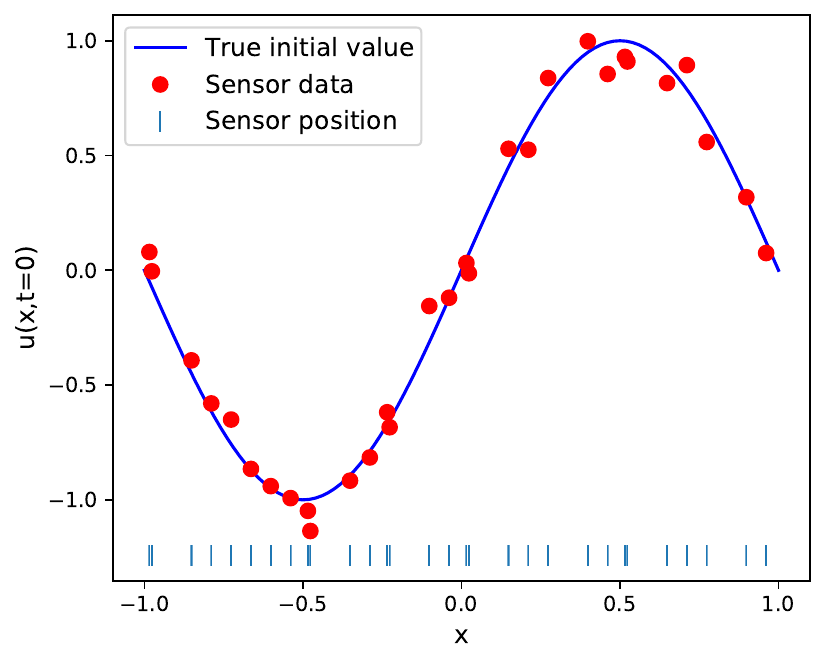}}\qquad
    \subfloat[]{\includegraphics[height=0.27\textheight, width=0.44\textwidth]{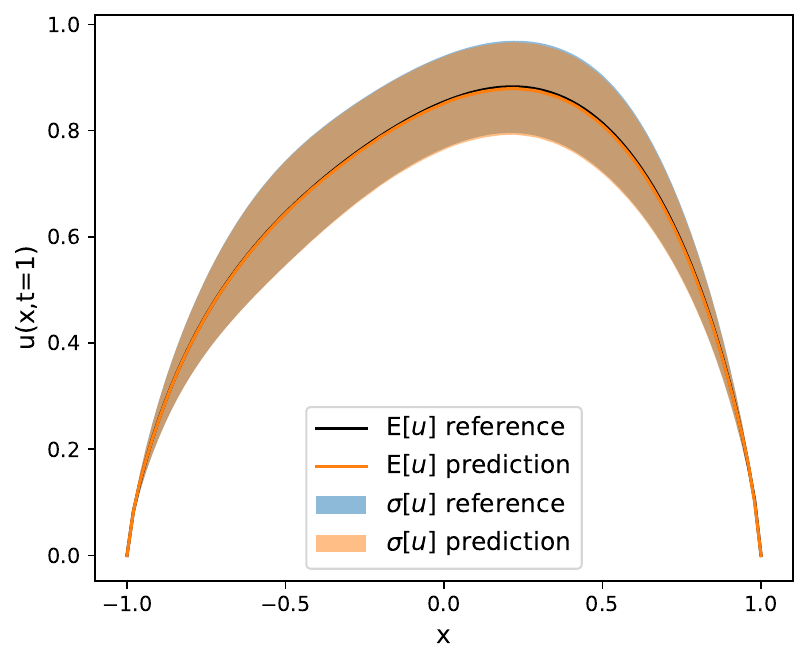}}
	\caption{Forward problem with time-evolving random forcing (noisy data). (a):~ Noisy sensor data as initial condition. (b):~ Mean and deviation of the predicted solution $u(x,t;\omega)$ versus the reference mean and standard deviation.}\label{fig:1D-addnoise-variance}
\end{figure}

\begin{figure}[!ht]
	\centering	
        \subfloat[]{\includegraphics[height=0.205\textheight, width=0.30\textwidth]{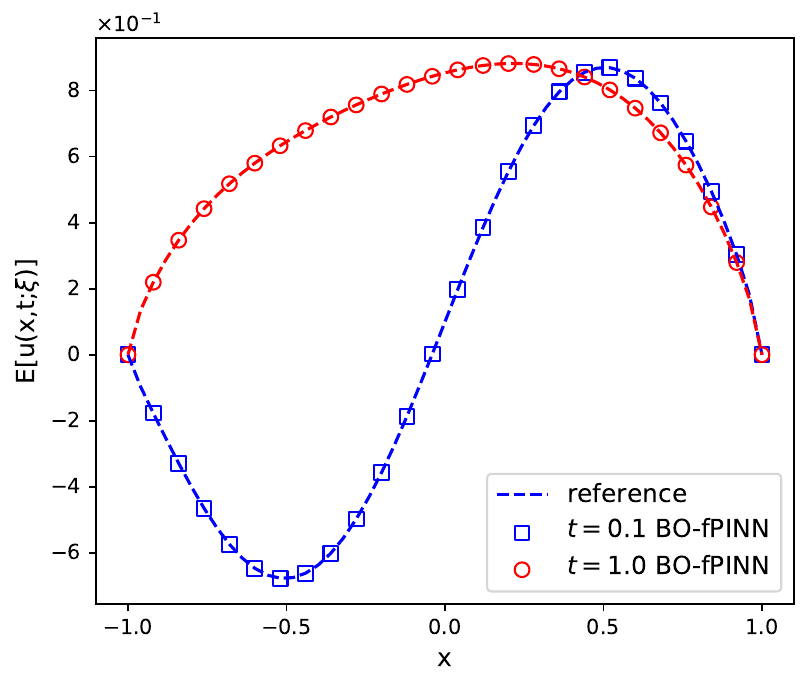}}\qquad
	\subfloat[]{\includegraphics[height=0.205\textheight, width=0.30\textwidth]{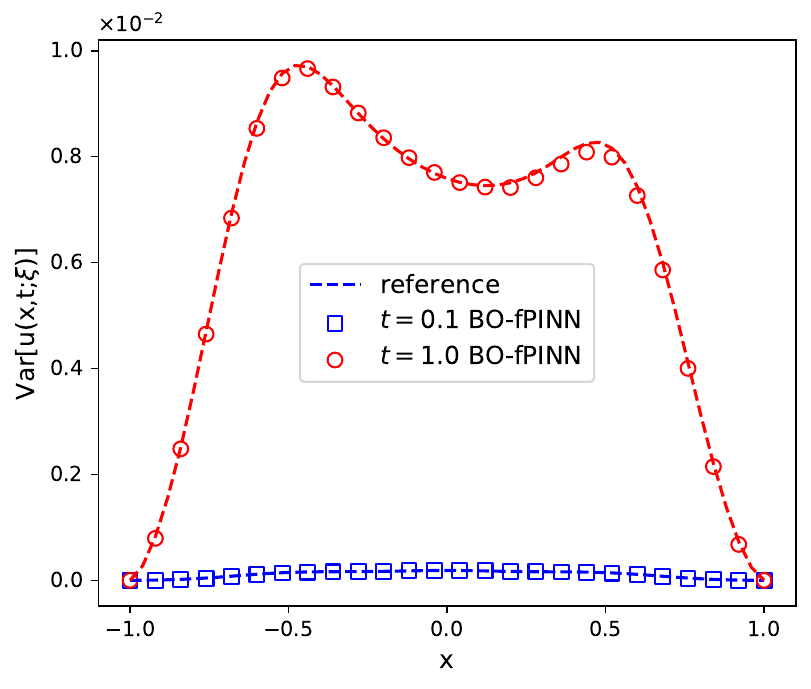}}\qquad
        \subfloat[]{\includegraphics[height=0.20\textheight, width=0.30\textwidth]{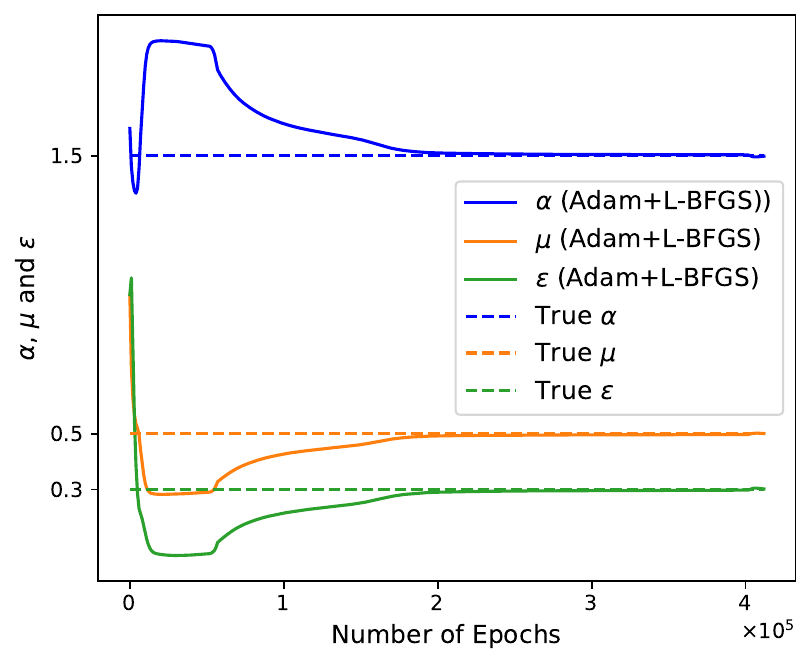}} \\	
	\caption{Inverse problem with time-evolving random forcing. ((a):~Solution mean at $t=0.1$ and $t=1$. (b):~Solution variance at $t=0.1$ and $t=1$. (c):~Parameter evolution as the iteration of optimizer progresses. The reference solutions are calculated by solving the forward problem using Quasi Monte Carlo simulation.}\label{fig:1D-inverse-mean-and-variance}
\end{figure}

\begin{figure}[!ht]
  \centering
  \includegraphics[height=0.36\textheight, width=0.86\textwidth]{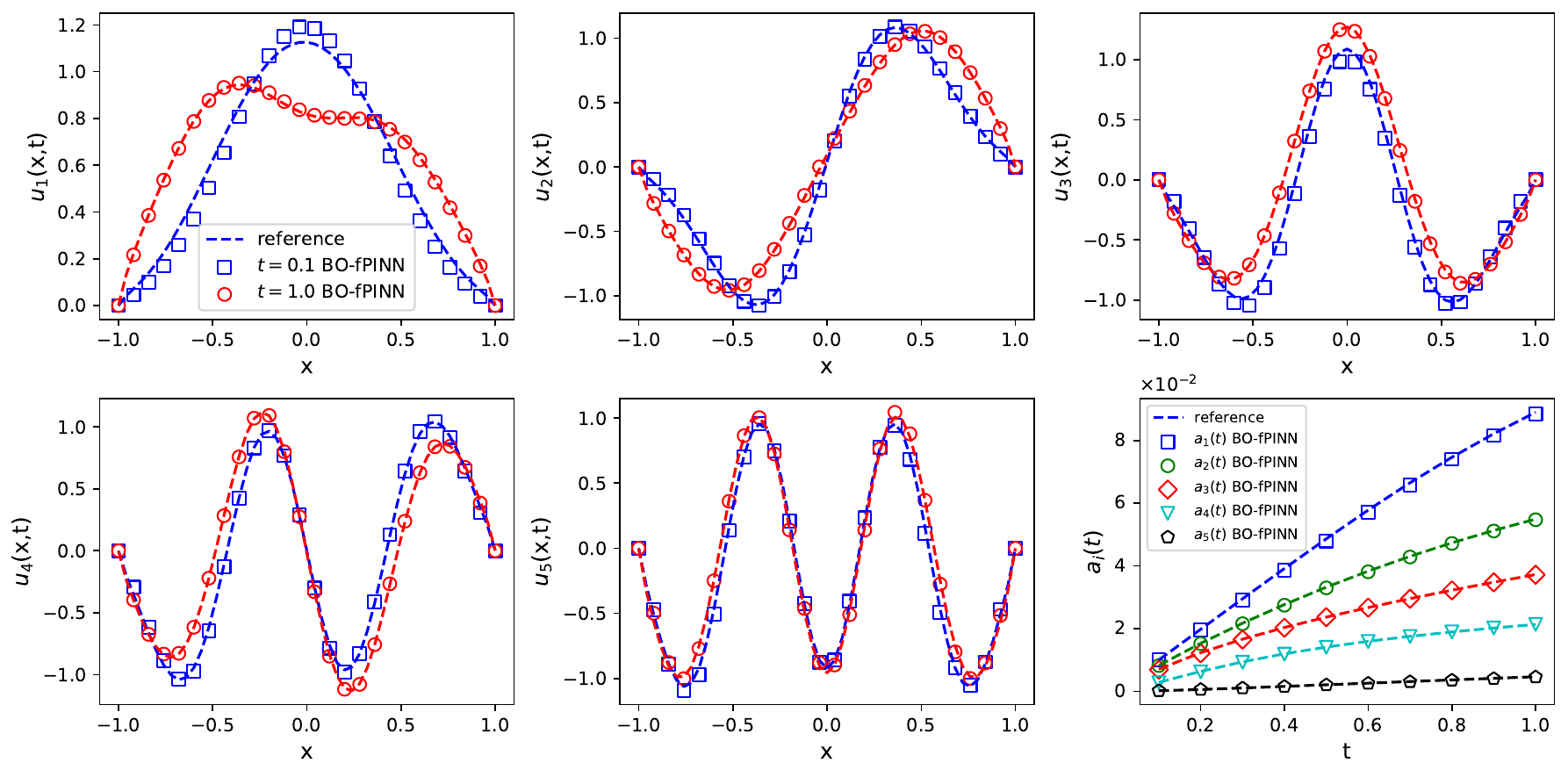}
  \caption{Inverse problem with time-evolving random forcing. The spatial modes at $t=0.1$ and $t=1$, the last subplot shows the scaling factors $a_i$ at different time steps, the reference solutions calculated using the QMC-BO method.}
  \label{fig:1D-inverse-ui-ai}
\end{figure}

For the 1D nonlinear stochastic fractional reaction-diffusion equation with time-evolving random forcing, we also consider the case with noisy sensor data as the initial condition, as well as the inverse problem.
The noisy sensor measurements of $u(x,t=0)$ is shown in Figure \ref{fig:1D-addnoise-variance}a, where the 30 sensors are uniformly placed in the domain and the measurements are corrupted by independent Gaussian random noise of standard deviation $0.1$. Figure \ref{fig:1D-addnoise-variance}b reports a comparison of the BO-fPINN solution mean and standard deviation obtained from the  noisy sensor data, and the reference mean and standard deviation calculated with the Quasi Monte Carlo simulation. The inverse problem considered here has the same setup as the one in section 4.2.2, the results are shown in Figure (\ref{fig:1D-inverse-mean-and-variance})-(\ref{fig:1D-inverse-ui-ai}) and Table \ref{tab:1D-inverse-error}.

\begin{table}[!htp]
\centering
\caption{Inverse problem with time-evolving random forcing. The relative $L_2$ errors of BO-fPINN solutions versus the reference solutions at the final time $t=1.0$. The reference solutions are calculated by solving the forward problem using the QMC-BO method.}
\label{tab:1D-inverse-error}
\begin{tabular}{c|cccccc}
\hline
          &  $a_1$  &  $a_2$   &  $a_3$  &  $a_4$  &  $a_5$    \\

Relative $L_2$ error  & 0.619$\%$ & 0.342$\%$ & 0.146$\%$ & 0.227$\%$ & 2.078$\%$ \\
\hline
          &  $u_1$  &  $u_2$   &  $u_3$  &  $u_4$  &  $u_5$   \\

Relative $L_2$ error  & 0.982$\%$ & 1.068$\%$ & 1.257$\%$ & 0.959$\%$ & 4.234$\%$  \\
\hline
          &  $Y_1$  &  $Y_2$   &  $Y_3$  &  $Y_4$  &  $Y_5$    \\

Relative $L_2$ error  & 0.0386 & 0.0334 & 0.0374 &  0.0314 &  0.0592 \\
\hline
\end{tabular}
\end{table}

\begin{table}[!htp]
\centering
\renewcommand\arraystretch{1.2}
\caption{Inverse problem with time-evolving random forcing. Identified parameters using the BO-fPINN method inverse problems with synthetic data.}
\label{tab:1D-parameters-inverse}
\begin{tabular}{c|c}
\hline 
 True parameters                            &  Indentified parameters   \\
\hline
 $\mu=0.5$, $\epsilon=0.3$, $\alpha=1.5$ & $\mu=0.501$, $\epsilon=0.302$, $\alpha=1.498$  \\
\hline 
\end{tabular}
\end{table}

\subsection{Additional results of Section 4.2.4}\label{appendix:1DRDE4.2.4}
Here we report the results for solving the target problem with  $\alpha = 1.5, 1.8$ by transfer learning method based on source problem with $\alpha = 1.2$.
\begin{figure}[htp!]
\centering
\subfloat[]{\includegraphics[width=0.40\linewidth]{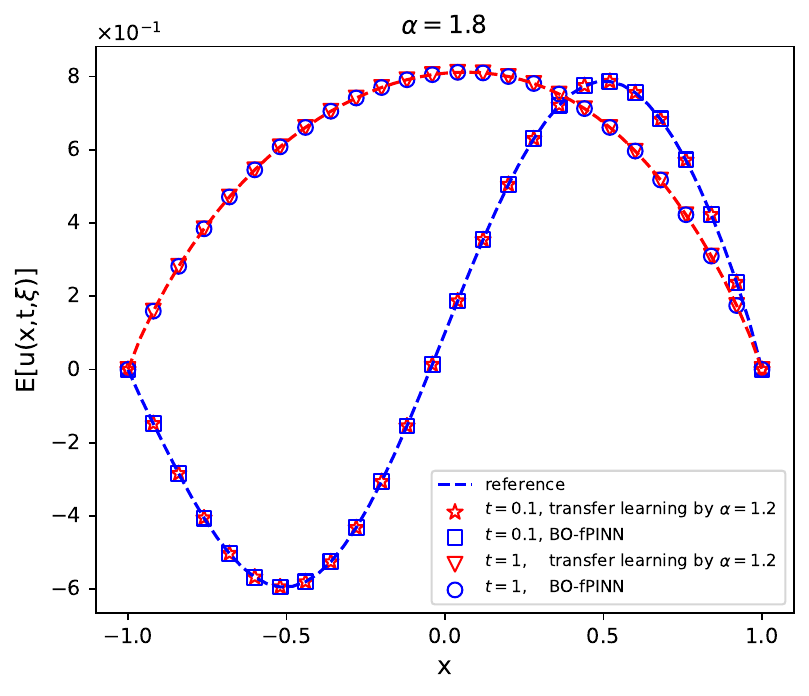}}\quad
\subfloat[]{\includegraphics[width=0.39\linewidth]{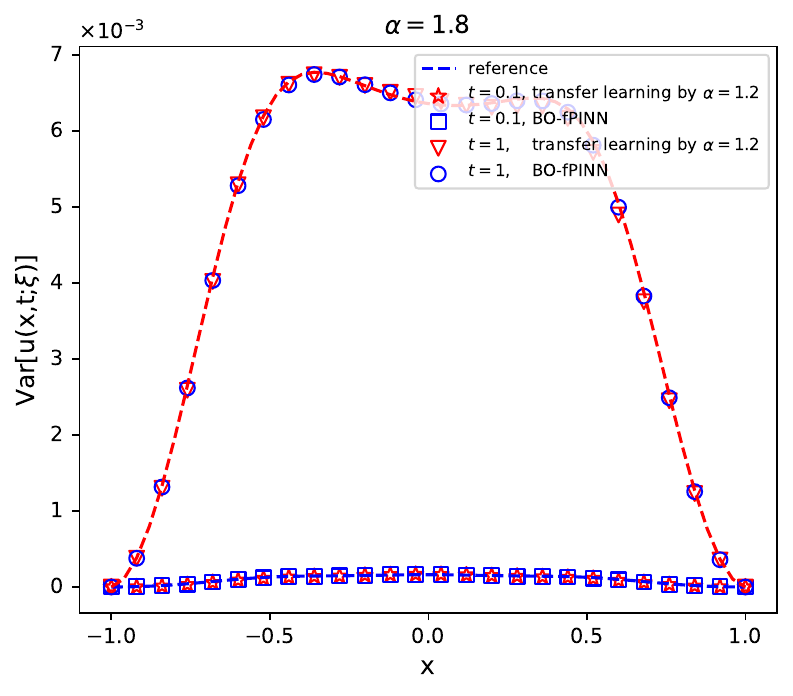}}
\caption{Forward problem with time-evolving random forcing (transfer learning between different fractional derivative orders). (a): Solution mean of $\alpha=1.8$ at $t=0.1, 1$ obtained by using BO-fPINN with Xavier initialization and by transfer learning with the parameters from source model $\alpha=1.2$ as initialization. (b): Solution variance of $\alpha=1.8$ at $t=0.1, 1$ by using BO-fPINN with Xavier initialization and by transfer learning with the parameters from source model $\alpha=1.2$ as initialization. The reference solutions are calculated by the QMC-BO method.}\label{fig:1D-transfer-mean-and-variance-byalpha12}
\end{figure}

\begin{figure}[!ht]
\centering
\subfloat[]{\includegraphics[width=0.40\linewidth]{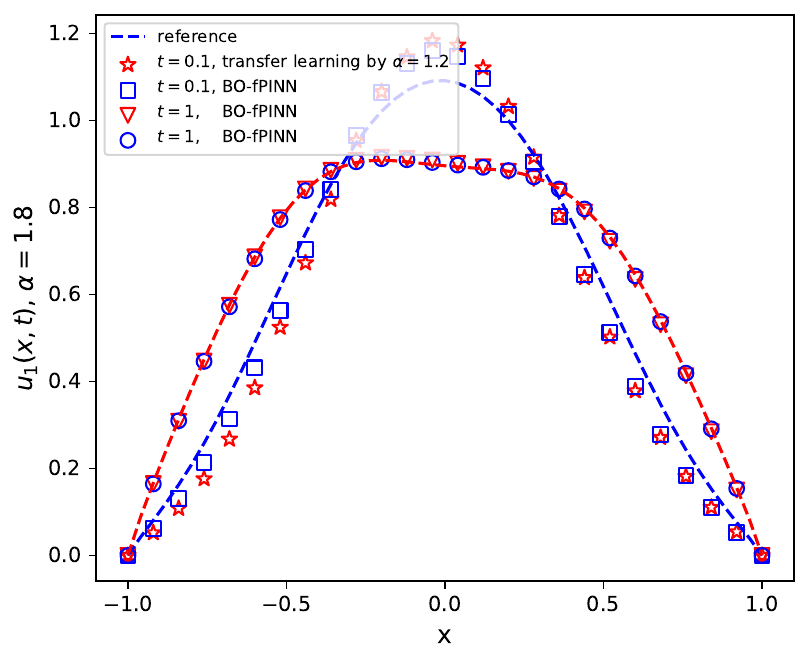}}\quad
\subfloat[]{\includegraphics[width=0.41\linewidth]{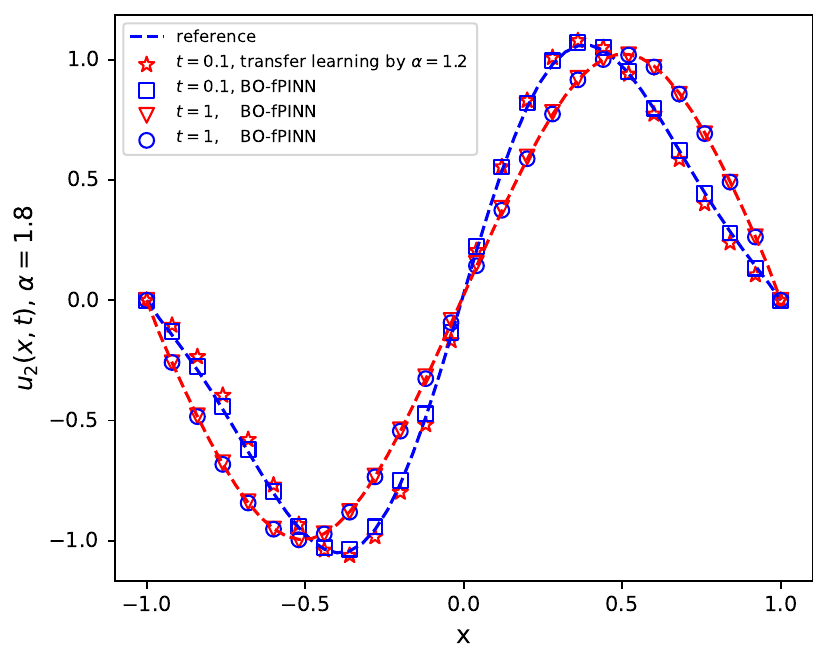}}
\caption{Forward problem with time-evolving random forcing (transfer learning between different fractional derivative orders). (a): The first BO modes at $t=0.1, 1$ for $\alpha=1.8$ obtained by transfer learning with source model $\alpha=1.2$. (b): The second BO modes at $t=0.1, 1$ for $\alpha=1.8$ obtained by transfer learning with source model $\alpha=1.2$. The reference solutions are calculated by the QMC-BO method.}\label{fig:1D-transfer-u1-u2-byalpha12}
\end{figure}

\begin{figure}[!ht]
\centering
\subfloat[]{\includegraphics[width=0.41\linewidth]{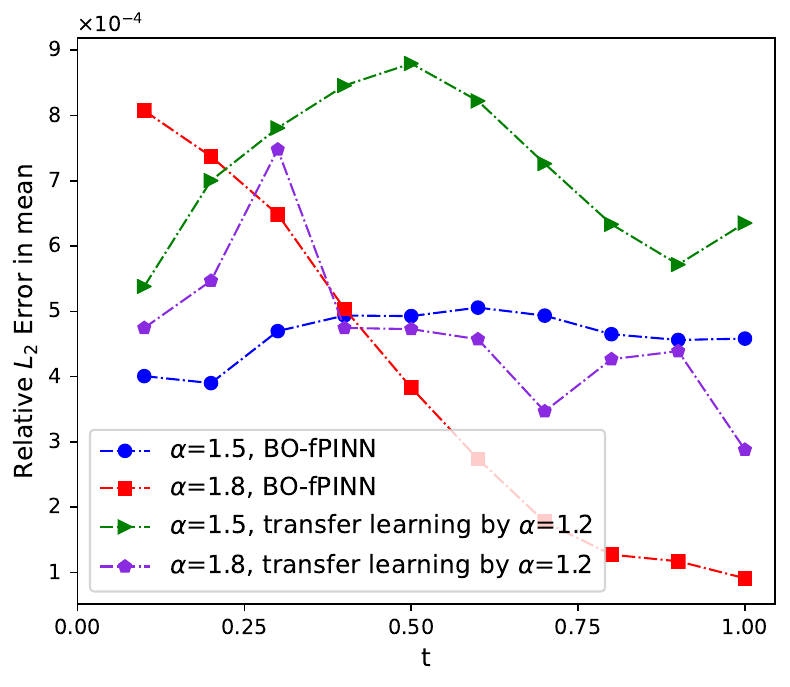}}\quad
\subfloat[]{\includegraphics[width=0.41\linewidth]{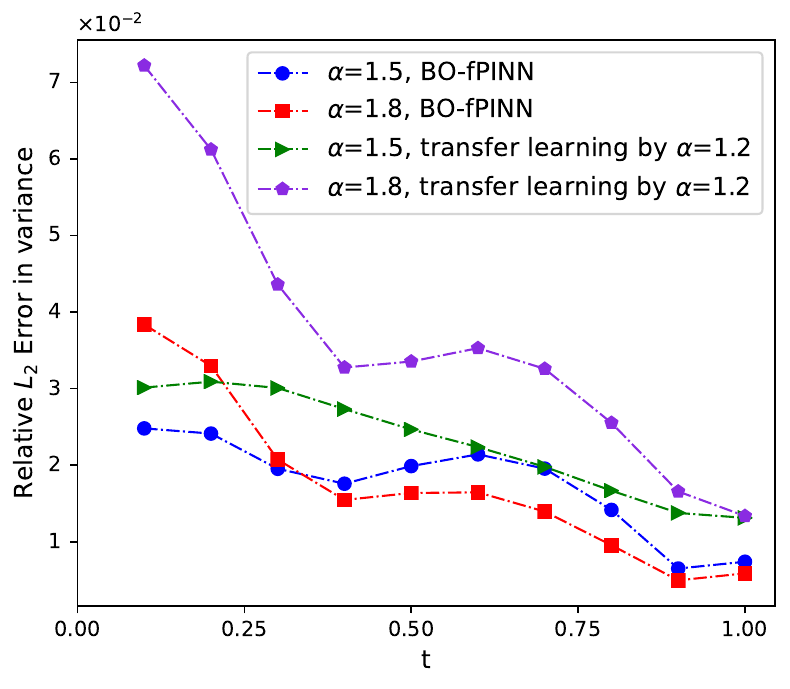}}\\
\quad
\caption{Forward problem with time-evolving random forcing (transfer learning between different fractional derivative orders). (a)~: The relative $L_2$ error in mean calculated by BO-fPINN and transfer learning with source model $\alpha=1.2$. (b)~: The relative $L_2$ error in variance calculated by BO-fPINN and transfer learning with source model $\alpha=1.2$. The reference solutions are calculated by the QMC-BO method.}\label{fig:1D-err-transfer-byalpha12}
\end{figure}

\end{document}